\documentclass[twoside,11pt,a4paper]{article}
\usepackage{times}
\usepackage{jmlr}
\usepackage[utf8]{inputenc}\usepackage{macros}

\definecolor{darkblue}{rgb}{0.0,0.0,0.7}
\RequirePackage[linktocpage=true,
hyperfootnotes=false,
colorlinks = true,citecolor = darkblue,urlcolor = darkblue, linkcolor=magenta, backref=page,
]{hyperref}
\usepackage{natbib}

\usepackage[backref=page]{hyperref}

\renewcommand*{\backrefalt}[4]{    \ifcase #1     \or (Cited on page:~#2.)     \else (Cited on pages:~#2.)    \fi    }

\definecolor{linkequation}{rgb}{0.0, 0.0, 1.0}
\usepackage{xcolor}
\newcommand*{\SavedEqref}{}
\let\SavedEqref\eqref
\renewcommand*{\eqref}[1]{  \begingroup
    \hypersetup{
      linkcolor=linkequation,
      linkbordercolor=linkequation,
    }    \SavedEqref{#1}  \endgroup
}

\newcommand*{\SavedRef}{}
\let\SavedEqref\ref
\renewcommand*{\ref}[1]{  \begingroup
    \hypersetup{
      linkcolor=linkequation,
      linkbordercolor=linkequation,
    }    \SavedRef{#1}  \endgroup
}

\usepackage[colorinlistoftodos]{todonotes}
\usepackage[symbol*]{footmisc}
\DefineFNsymbolsTM{myfnsymbols}{              }\usepackage{multicol}
\setfnsymbol{myfnsymbols}
\usepackage{tabularx}
\allowdisplaybreaks

\usepackage{authblk}

\begin{document}

\title{Bounds in Wasserstein Distance for Locally Stationary Functional Time Series}

\author[1,2]{Jan Nino G. Tinio}
\author[1]{Mokhtar Z. Alaya}
\author[1]{Salim Bouzebda}
\affil[1]{\footnotesize Université de Technologie de Compiègne,  \authorcr
Laboratoire de Mathématiques Appliquées de Compiègne, CS 60 319 - 60 203 Compiègne Cedex
}
\affil[2]{\footnotesize Department of Mathematics, Caraga State University, Butuan City, Philippines\linebreak\linebreak\texttt{e-mails}: \texttt{jan-nino.tinio@utc.fr}, \quad \texttt{alayaelm@utc.fr}, \quad \texttt{salim.bouzebda@utc.fr}}

\maketitle

\begin{abstract}
Functional time series (FTS) extend traditional methodologies to accommodate data observed as functions/curves. A significant challenge in FTS consists of accurately capturing the time-dependence structure, especially with the presence of time-varying covariates. When analyzing time series with time-varying statistical properties, locally stationary time series (LSTS) provide a robust framework that allows smooth changes in mean and variance over time.  This work investigates Nadaraya-Watson (NW) estimation procedure for the conditional distribution of locally stationary functional time series (LSFTS), where the covariates reside in a semi-metric space endowed with a semi-metric. Under small ball probability and mixing condition, we establish convergence rates of NW estimator for LSFTS with respect to Wasserstein distance. The finite-sample performances of the model and the estimation method are illustrated through extensive numerical experiments both on functional simulated and real data.
\end{abstract}


\section{Introduction} \label{sec:introduction}
In recent years, data collected across various fields increasingly exhibit functional or curve-like characteristics, commonly referred to as functional data whose values come from infinite-dimensional space. This advancement is driven by the proliferation of data collected on progressively refined temporal and spatial grids, for instance, in meteorology, medicine, satellite imagery, economics and finance, environmental science, and many others \citep{soukarieh2024weak, Bouzebda20234, MR4384114}. Specifically, \cite{Chenetal2016} demonstrated intriguing applications in biometrics, \cite{AuevanDelft2020} explored the relevance of functional data to environmental science, and \cite{Bugnietal2009, Aueetal2015} focused on applications in econometrics. The statistical modeling of such data, conceptualized as random functions/curves, has given rise to various intricate theoretical and numerical research inquiries. Functional data are investigated for multiple reasons: they assist in developing data representation strategies that highlight significant features, analyzing mean behaviors and deviations, and constructing models when temporal dependency is present in the data \citep{Ramsay&Silverman2002, Ferraty&Vieu2006}. 

The statistical challenges in studying these data have drawn increasing attention over the last few decades from the statistical community \citep{MR4507275, MR4279953, MR4695390}. Much of the existing research in FTS analysis is predicated on stationary models since the stationarity requirement is often adopted in time series modeling, leading to various models and methodologies.  Theoretical foundations of function spaces concerning linear processes were explored in \cite{Bosq2000}, which are crucial for modeling and forecasting functional time series. Statistical analysis of autoregressive processes within Banach spaces was emphasized in \cite{Bosq2002}, a specialized class of functional spaces. Similarly, \cite{DehlingSharipov2005} examined asymptotic properties of the sample mean and variance for autoregressive processes of order one in a separable Banach space with independent innovations. Extending this, \cite{Antoniadis2006} incorporated wavelet analysis alongside nonparametric kernel regression for FTS forecasting. 

However, many FTS exhibit nonstationarity. Even with detrending and deseasonalizing, the stationarity assumption is not always beneficial for modeling functional data. Many time series models, commonly observed in various physical phenomena and economic data, are non-stationary \citep{Daisuke2022, soukarieh2024weak, Bouzebda20234}. Conventional approaches become inappropriate when the (weak) stationarity assumption is violated. Specifically, global climate changes in meteorology affect the distribution of a region's daily temperature, precipitation, and cloud cover when viewed interconnectedly. In finance, an option’s implied volatility evolves over time depending on its moneyness \citep{vanDelft2018}. 
To address such nonstationarity, the locally stationary time series (LSTS) framework provides a practical modeling approach. 
The parameters of LSTS exhibit temporal dependence. These nonstationary processes can be approximated by locally stationary counterparts that remain stationary within smaller time windows. As a result, asymptotic theories can be developed to estimate time-dependent characteristics \citep{Dahlhaus1996, DahlhausSubbaRao2006, DAHLHAUS2012351}. Many works have addressed the intuitive concept of local stationarity. The seminal works in \cite{Dahlhaus1996, Dahlhaus2000, DahlhausSubbaRao2006} provide a strong foundation for the inference of LSTS. Parametric \citep{DahlhausSubbaRao2006, Fryzlewiczetal2008, DAHLHAUS2012351} and nonparametric \citep{vogt2012, Daisuke2022, SoukariehSISP} frameworks were developed to analyze LSTS. For instance, \cite{Dahlhaus1996} introduced a parametric framework by leveraging local periodograms to minimize the generalized Whittle function tailored to locally stationary models. Expanding the parametric framework beyond time-invariant assumptions, \cite{KooLinton2012} proposed a semiparametric model. This approach combined the strengths of parametric and nonparametric methods, allowing for greater flexibility in modeling complex time-varying structures without imposing stringent parametric constraints. Parallel to the parametric advancements, nonparametric approaches have garnered significant attention from researchers aiming to model conditional mean and variance functions without relying on predefined parametric forms. In particular, \cite{Daisuke2022} developed an estimation theory for nonparametric regression problems involving LSFTS. Central to nonparametric approaches in the literature is Nadaraya-Watson (NW) \citep{MR166874, MR185765} estimation procedure, which has been extensively used to estimate conditional mean functions.

Despite the proven efficacy of NW procedure in estimating conditional mean functions, its usefulness extends to conditional distribution estimation (CDE). Estimating accurately conditional distributions is essential in prediction and forecasting, as it comprehensively describes the conditional law for any given random variable. 
Due to its critical role in predictive modeling and inference, numerous studies have focused on developing robust estimation theories for CDE within FTS framework \citep{ferraty2005conditional, Hormmanetal2022, MR4783436}.  Two distinct approaches for estimating the conditional distribution of a target variable within a prediction set, given a functional covariate, were introduced in \cite{Hormmanetal2022}. The first method relies on the empirical distribution of the estimated model residuals, while the second involves fitting functional parametric models to the residuals. These approaches provide flexible frameworks for handling complex dependencies between functional covariates and target variables, particularly in high-dimensional settings where traditional parametric models may be inadequate.
Additionally, \cite{MR4783436} developed a local polynomial estimator for conditional CDF of a scalar target variable given a functional covariate in the context of stationary strongly mixing processes. They established asymptotic normality property. 

In dealing with CDE, optimal transport (OT) theory has emerged as a powerful mathematical framework for quantifying difference between probability distributions. OT measures the minimal cost required to transport one distribution to another, providing a meaningful metric for distributional comparison \citep{PeyreCuturi2020}. This approach addresses shortest-path problems by enabling the concurrent transportation of multiple items along geodesic curves or straight paths. Among various metrics derived from OT, Wasserstein distance stands out for its robustness and versatility in comparing probability distributions \citep{Villani2009OToldnew, Manoleetal2022, Tinioetal1.2024}. 

The landscape of LSTS analysis is enriched by parametric and nonparametric approaches, each contributing unique strengths to the modeling of time-evolving data. Parametric approaches offer structured and interpretable models for capturing dynamic behaviors through frameworks like those proposed in \cite{Dahlhaus1996, KooLinton2012}. Concurrently, nonparametric techniques, particularly those leveraging NW estimator and OT theory, provide flexible and robust tools for estimating conditional means, variances, and distributions. The integration of these procedures, underpinned by rigorous theoretical advancements, continues to enhance the capacity of statisticians and data scientists to model, predict, and infer from complex, non-stationary time series data. 

\paragraph{Contributions.} 
This paper investigates CDE for LSFTS characterized by weakly dependent sequences using NW estimation procedure. We establish convergence rates of NW estimator for conditional distribution of \( Y_{t, T} \mid X_{t, T} \) with respect to Wasserstein distance. Here, the response variable \( Y_{t, T} \) is scalar, while locally stationary covariate \( X_{t, T} \) resides in a semi-metric space \( \mathscr{H} \) endowed with a semi-metric \( \mathsf{D}(\cdot, \cdot) \). We perform numerical experiments on synthetic and real-world data to illustrate our theoretical findings.

\paragraph{Layout.}
The paper is structured as follows: Section \ref{sec:preliminaries} introduces the regression estimation problem and provides an overview of key concepts, including local stationarity, Wasserstein distance, small ball probability, and mixing conditions. In Section \ref{sec:theoretical_guarantees}, we present main theoretical results, which include \textit{(i)} the formulation of kernelized NW estimator, \textit{(ii)} the explicit convergence rates of  Wasserstein distance between NW estimator and true conditional distributions, and \textit{(iii)} the proposed method for bandwidth selection. 
The finite-sample performance of the model and estimation method is presented in  Section \ref{sec:numerical_experiments}  through the analysis of both simulated and real data. To preserve the flow of the presentation, all proofs are deferred to the appendices.

\paragraph{Notation.}
Throughout this paper, we adopt the following notations. Dirac measure at a point \( y \) is denoted by \( \delta_y \). For any real-valued random variable \( X \) and any \( q \geq 1 \), we denote \( L_q \)-norm of \( X \) by \( \|X\|_{L_q} \) and is defined as
\(
\|X\|_{L_q} = \left( \mathbb{E}\left[ |X|^q \right] \right)^{1/q}.
\)
We use the notation \( a_T \lesssim b_T \) to indicate that there exists a constant \( C \), independent of \( T \), such that
\(
a_T \leq C b_T.
\)
The constant \( C \) may vary unless specified otherwise. Similarly, \( a_T \sim b_T \) signifies that both \( a_T \lesssim b_T \) and \( b_T \lesssim a_T \) hold. For positive sequences \( \{a_T\} \) and \( \{b_T\} \), we write \( a_T = \mathcal{O}(b_T) \) provided that
\(
\lim_{T \to \infty} \frac{a_T}{b_T} \leq C
\)
for some constant \( C > 0 \). Additionally, \( a_T = \mathcal{O}(1) \) indicates that \( a_T \) is bounded. We denote \( a_T = o(b_T) \) if
\(
\lim_{T \to \infty} \frac{a_T}{b_T} = 0,
\)
and \( a_T = o(1) \) when \( a_T \) approaches zero. For a sequence of random variables \( \{X_T\} \) and a given sequence \( \{a_T\} \), we write \( X_T = \mathcal{O}_\mathbb{P}(a_T) \) if, for every \( \epsilon > 0 \), there exist constants \( C_\epsilon > 0 \) and \( T_0(\epsilon) \in \mathbb{N} \) such that for all \( T \geq T_0(\epsilon) \),
\(
\mathbb{P}[ \frac{|X_T|}{a_T} > C_\epsilon ] < \epsilon.
\)
Similarly, \( X_T = o_\mathbb{P}(a_T) \) indicates that
\(
\lim_{T \to \infty} \mathbb{P}[ \frac{|X_T|}{a_T} > \epsilon ] = 0
\)
for all \( \epsilon > 0 \), and \( X_T = o_\mathbb{P}(1) \) when \( X_T \) converges in probability to zero. For any \( a, b \in \mathbb{R} \), we define
\(
a \vee b = \max\{a, b\} \quad \text{and} \quad a \wedge b = \min\{a, b\}.
\)


\section{Background on LSTS and Wasserstein distance} \label{sec:preliminaries}

This section presents some preliminaries of LSTS and Wasserstein distance.  We then introduce small ball probability and delineate the mixing condition employed to assess weak dependency.

\subsection{Locally stationary time series} \label{ssub:background_of_locally_stationary_processes_and_wasserstein_distance}

Let $T \in \mathbb{N}$ and suppose that there exists a process of $T$ random variables $\{Y_{t,T}, X_{t,T}\}_{t=1, \ldots, T}$, where $Y_{t,T}$ is real-valued,  and $X_{t,T}$ belongs to a semi-metric space $\mathscr{H}$ with a semi-metric $\mathsf{D} (\cdot,\cdot)$. The semi-metric space $\mathscr{H}$ can be Banach or Hilbert spaces with norm $\norm{\cdot}$.    
We consider the following regression estimation problem:
\begin{align}
\label{eq:major_estimation_problem}
Y_{t,T} = m^\star\big(\frac tT, X_{t,T}\big) + \varepsilon_{t,T}, \text{ for all }t=1, \ldots, T,
\end{align}
where  $\{\varepsilon_{t,T}\}_{t\in \mathbb{Z}}$ is a sequence of independent and identically distributed (i.i.d.) random variables independent of $\{X_{t,T}\}_{t=1, \ldots, T}$, i.e., $\E[\varepsilon_t|X_{t,T}] = 0,$ for all $t=1, \ldots, T$. The response $Y_{t, T}$ is assumed to be integrable. The functional covariate $X_{t, T}$ is assumed to be locally stationary, that dynamically changes slowly over time and hence can be considered approximately stationary at local time.  Note that $m^\star\big(\frac tT, X_{t,T}\big) = \E[Y_{t,T}|X_{t,T}]$ is the \textit{oracle} conditional mean function in model~\eqref{eq:major_estimation_problem}, and does not depend on real-time $t$ but rather on rescaled time $u=\frac tT$. As the sample size $T$ goes to infinity, these $u$-points form a dense subset of the unit interval $[0,1]$. Hence, at all rescaled $u$-points, $m^\star$ is identified almost surely (a.s.) if it is continuous in the time direction. In LSTS, this rescaled time refers to transforming the original time scale. 

\paragraph{Example.} Consider the process $X_{t,T} = a \big(\frac{t}{T}\big) + \varepsilon_t$, $t\in \{1,\ldots,T\}$, $T\in \mathbb{N}$, where  $a(\cdot)$ is a continuous function $a: [0,1] \rightarrow \R$ and a sequence of i.i.d. random variables $\{\varepsilon_t\}_{t \in \mathbb{N}}$. The process $X_{t,T}$ behaves ``almost'' stationary for $t$ close to $t'$, for some  $t' \in \{1,\ldots,T\}$, that is, $a \big(\frac{t'}{T}\big) \approx a \big(\frac{t}{T}\big)$. However, this process is not weakly stationary. Local stationarity gives a more realistic concept that allows this kind of change \cite{SoukariehSISP}.

Let us now formally define the notion of LSTS. We adopt the definition given in \cite{Daisuke2022}.

\begin{definition}
\label{definition:locallystatpr}
An $\mathscr{H}$-valued process $\{X_{t,T}\}_{t=1,\ldots,T}$ is locally stationary if for each rescaled time point $u\in[0,1]$, there exists an associated $\mathscr{H}$-valued process $\{X_t(u)\}_{t=1,\ldots,T}$ verifying
\begin{align}\label{eqn: dis (X, Xu)}
\mathsf{D} \big({X_{t,T}, X_t(u)}\big) &\leq \Big(\Big|\frac tT - u\Big| + \frac 1T\Big) U_{t,T}(u) \text{  a.s.,}
\end{align}
where $\{U_{t,T}(u)\}_{t=1,\ldots,T}$ is a positive process such that $\E\big[(U_{t,T}(u))^\rho\big]< C_U$ for
some $\rho > 0$ and $C_U < \infty$ independent of $u, t,$ and $T$.
\end{definition}

From this definition, if an $\mathscr{H}$-valued process $\{X_{t,T}\}_{t=1,\ldots,T}$ is locally stationary, a strictly stationary process $\{X_t(u)\}_{t=1,\ldots,T}$ can always be found around each rescaled time $u$, which will be used to approximate $\{X_{t,T}\}_{t=1,\ldots,T}$. This approximation will result in a negligible difference between random variables $X_{t, T}$ and $X_t(u)$. Since the $\rho$-th moments of the positive random variables $U_{t,T}(u)$ are uniformly bounded, $U_{t,T}(u)=O_{\mathds{P}}(1)$ \citep{vogt2012}. Thus, we have
\begin{align*}
    \mathsf{D} \big({X_{t,T}, X_t(u)}\big) &= O_{\mathds{P}}\Big(\Big|\frac{t}{T}-u\Big|+\frac{1}{T}\Big).
\end{align*}
If $u=\frac{t}{T}$, then  $\mathsf{D} \big({X_{t,T}, X_t(\frac{t}{T})}\big)\leq \frac{C_U}{T}$.  The constant $\rho$ can be considered as an indicator of how well this approximation is being made: larger $\rho$ indicates a better approximation of $X_{t, T}$ by $X_t(u)$ and moderate bounds for their absolute difference.

Definition \ref{definition:locallystatpr} is consistent with the one given in \cite{vanDelft2018, vanDelftDette2023} when $\mathscr{H}$ is the Hilbert space $L_{\R}^2([0,1])$, with inner product $L_2$-norm:
\begin{equation*}
    \norm{f}_2 = \sqrt{\langle f,f \rangle},\quad \langle f,g \rangle = \int_0^1 f(t)g(t)\diff t,
\end{equation*}
where $f,g \in L_{\R}^2([0,1])$. Sufficient conditions were also provided so that an $L_{\R}^2([0,1])$-valued stochastic process $X_{t,T}$ satisfies (\ref{eqn: dis (X, Xu)}) with $\mathsf{D} (f,g) = \norm{f - g}_2$ and $\rho=2$. In \cite{vanDelft2018}, $L_{E}^p(I, \mu)$ is defined as the Banach space of all strongly measurable functions $f: I \rightarrow E$ with norms
\begin{equation*}
    \norm{f}_p = \norm{f}_{L_{E}^p(I, \mu)} = \Big( \int \norm{f(s)}_E^p \diff \mu(s) \Big)^\frac{1}{p},
\end{equation*}
for $1 \leq p < \infty$, and 
\begin{equation*}
    \norm{f}_\infty = \norm{f}_{L_{E}^\infty(I, \mu)} = \inf_{\mu(N)=0} \sup_{s\in I \setminus N} \norm{f(s)}_E,
\end{equation*}
for $p=\infty$. 

\subsection{Wasserstein distance}

Let $\mathcal{P}_r(\R)$ be the set of Borel probability measures in $\R$ having finite $r$-th moment $(r\geq 1)$, i.e., $$\mathcal{P}_r(\R)=\{\mu\in\mathcal{P}(\R): \int_\R |x|^r \mu(\mathrm{d}x) < \infty \}.$$ Given probability measures $\mu, \nu \in \mathcal{P}_r(\R)$, we calculate the distance between them using the $r$th-Wasserstein distance, $W_r(\mu,\nu)$, as follows
\begin{align}\label{def:W1_general}
    W_r(\mu,\nu) &= \Big(\inf_{\pi\in \Pi(\mu,\nu)} \iint_{\R \times \R} {|u- v|}^r\pi(\diff u, \diff v)\Big)^{1/r},
\end{align}
where $\Pi(\mu,\nu)$ denotes the set of probability measures on $\R \times \R$ with marginals $\mu$ and $\nu$. Optimal couplings always exist since $\R$ is a complete and separable metric space, where the infimum is, in fact, a minimum \citep{Villani2009OToldnew}. Equation (\ref{def:W1_general}) signifies that $W_r(\mu, \nu)$ is the infimum of expectation of the distance between two random variables over all possible couplings, i.e., $W_r(\mu,\nu) = \big( \inf_{U\sim \mu, \text{ } V\sim \nu}  \E [|U-V|^r]\big)^{1/r}$.

We can represent a simple optimal coupling by a probability inverse transform: given $\mu, \nu \in \mathcal{P}_r(\R)$, let $F_\mu(\cdot)$ and $F_\nu(\cdot)$ be the cumulative distribution functions (CDFs) and $F_\mu^{-1}(\cdot)$ and $F_\nu^{-1}(\cdot)$ be the respective generalized inverse or quantile functions defined as $F_\mu^{-1}(z):= \inf\{v\in\R: \mu ((-\infty,v])\geq z\}$ for all $z\in[0,1]$ (similarly for $F_\nu^{-1}(z)$). Then, assuming a random variable $Z$ uniformly distributed on $(0,1)$, an optimal coupling $(U,V)=(F_\mu^{-1}(Z),F_\nu^{-1}(Z))$ can be established \citep{DedeckerMerleved2017, Dombry2023}. Hence, the minimization problem (\ref{def:W1_general}) can be represented by
\begin{align*}
    W_r(\mu,\nu) &= \Big(\int_0^1 {\left|F_\mu^{-1}(z) - F_\nu^{-1}(z)\right|}^r \diff z \Big)^{1/r},
\end{align*}
in a one-dimensional context. Specifically, for $r=1$ and using a change of variable, the $1$-Wasserstein distance is represented as
\begin{align}\label{def:W1_cdf}
    W_1(\mu,\nu) &= \int_{\R} |F_\mu(v) - F_\nu(v)| \diff v.
\end{align}
Consequently, $W_1(\mu,\nu)$ can be considered as the $L_1$-distance between the CDFs $F_\mu(\cdot)$ and $F_\nu(\cdot)$.

\subsection{Small ball probability}\label{subsection: small ball probability}
The absence of a density function for functional random variables is a technical difficulty in infinite-dimensional spaces since we lack a universal reference such as the Lebesgue measure. We overcome this using the \textit{small ball probability} property. We control the concentration of probability measure of the functional variable on a  small ball using a function $\phi(r)$ defined as, for all $r>0$ and a fixed $x\in \mathscr{H}$, 
\begin{equation*}    \P[X\in B(x,r)] =: \phi_x(r) > 0.
\end{equation*}
where $B(x,r) = \{v\in \mathscr{H}: \mathsf{D} (x,v)\leq r. \}$
Assume that $r$ is a function of $T$ such that $r = r(T) \rightarrow 0$ as $T \rightarrow \infty$. If we take $T$ very large, $B(x, r)$ is then considered as a small ball; hence, $\P[X\in B(x,r)]$ is a small ball probability \citep{FerratyVieu2004}.  Unfortunately, obtaining $\P[X\in B(x,r)]$ is complicated \citep{Ferratyetal2007}. For a survey on the main results on small ball probability, refer to \cite{LiShao2001}. In most cases, it is fitting to suppose that, as $r\rightarrow 0$,
\begin{equation}\label{eqn: assumed factorization on small ball}
    \P[X\in B(x,r)] \sim \psi(x) \phi(r),
\end{equation}
where $\E[\psi(X)] = 1$, a necessary normalizing restriction to ensure the identifiability of the decomposition. There are two main reasons for conveniently assuming (\ref{eqn: assumed factorization on small ball}). First, the function $\psi(x)$ can be thought as a surrogate density of the functional $X$ and can be utilized in different frameworks where the surrogate density is estimated differently and is used for classification purposes. Second, the function $\phi(r)$ signifies the volumetric term that can be used to evaluate the complexity of the probability law of the process \citep{Bongiornoetal2018}. In the $d$-dimensional case $X\in \R^d$, we suppose $\phi(r) \sim  r^d$, which is commonly referred as the curse of dimensionality \citep{FerratyVieu2004, Ferratyetal2006}. The intrinsic nature of the probability effects, involving small balls, is apparent in infinite-dimensional framework. We give some examples of several forms of $\phi(r)$ that can also be found in \cite{Ferratyetal2006, SoukariehSISP}.

\paragraph{Fractional Brownian Motion.}
Considering the space $\mathcal{C}([0,1], \R)$ with the supremum norm and its Cameron-Martin associated space $\mathcal{F} = \mathcal{C}([0,1], \R)^{CM}$. 
    Using Theorems 3.1 and 4.6 in \cite{LiShao2001}, for $0 < \eta < 2$, we have
        \begin{equation*}
            \forall x \in \mathcal{F}, \quad C'_{x} \mathrm{e}^{r^{-2/\eta}} \leq \P [ \zeta^{FBM} \in B(x, r)] \leq C_{x} \mathrm{e}^{r^{-2/\eta}},
        \end{equation*}
    where $\zeta^{FBM}$ is the usual Fractional Brownian Motion with parameter $\eta$ and $B(x, r) = \{ \zeta^{FBM} \in \mathcal{F}: \norm{\zeta^{FBM} - x}_\infty \leq r \}$. In this example, for the Fractional Brownian process, we choose $\phi(r)$ of the form
        \begin{equation*}
            \phi^{FBM}(r) \sim \mathrm{e}^{r^{-2/\eta}}.
        \end{equation*}

\paragraph{Gaussian process.}
Next, let us consider the centered Gaussian process $\zeta^{GP} = \{\zeta_t^{GP}, 0\leq t \leq 1\}$, which can be expanded by Karhunen-Loève decomposition as
        \begin{equation*}
            \zeta_t^{GP} = \sum_{i=1}^\infty \sqrt{\lambda_i} Z_i f_i(t),
        \end{equation*}
    where $\lambda_i$'s are the eigenvalues of the covariance operator of $\zeta^{GP}$, $f_i$'s are the associated orthonormal eigenfunctions, and $Z_i$'s are independent standard normal real random variables. The orthogonal projection onto the subspace spanned by the eigenfunctions $\{f_1, \ldots, f_k\}$ is denoted by $\Xi_k$, for $k\in\mathbb{N}^*$. Define a semi-metric by
        \begin{equation*}
            \mathsf{D}(x,y) = \int_0^1 (\Xi_k(x-y)(t))^2 \diff t.
        \end{equation*}
    Using the Karhunen-Loève expansion, we get
        \begin{equation*}
            \mathsf{D} (\zeta^{GP},x) = \sum_{i=1}^k \Big(\sqrt{\lambda_i} Z_i - x_i \Big)^2 := \sum_{i=1}^k \chi_i^2,
        \end{equation*}
    where
        \begin{equation*}
            x_i = \int_0^1 x(t) f_i(t) \diff t, \quad i=1, \ldots, k.
        \end{equation*}
    $\mathsf{D} (\zeta^{GP},x)$ can be written in terms of the usual Euclidian norm on $\R^k$ of a vector $\chi = (\chi_1, \ldots, \chi_k)$. Since $\chi_i$'s are independent real random variables with density with respect to the Lebesgue measure, we have, for $B(x, r) = \{ \zeta^{GP} \in \mathcal{F}: \mathsf{D}(\zeta^{GP},x)<r \}$
        \begin{equation*}
            \P[\zeta^{GP} \in B(x, r)] \sim r^k.
        \end{equation*}

\paragraph{Ornstein-Unhlenbeck process.}
Lastly, considering the same space in (i) and the metric $\mathsf{D} (\cdot , \cdot)$ associated with the supremum norm
        \begin{equation*}
            \forall x \in \mathcal{C}([0,1], \R), \quad \norm{x}_\infty = \sup_{t\in [0,1]} |x(t)|.
        \end{equation*}
    We denote the Wiener measure on $\mathcal{C}([0,1], \R)$ by $\P^W$ and the associated functional Cameron-Martin space of $\mathcal{C}([0,1], \R)$ is given by $\mathcal{F} = \mathcal{C}([0,1], \R)^{CM}$. Moreover, we denote the standard Wiener process by $w$ and let us consider the Ornstein-Unhlenbeck process $\zeta^{OU}$ defined by $\zeta_0^{OU} = 0$ and by
        \begin{equation*}
            \diff \zeta_t^{OU} = \diff w_t - \frac{1}{2}\zeta_t^{OU}, \quad \forall t, 0<t\leq 1.
        \end{equation*}
    By \cite{Bogachev1998}, the small centered ball Wiener measures are known to be of the form
        \begin{equation*}
            \P^W [\norm{x}_\infty \leq r] \sim \frac{4}{\pi} \mathrm{e}^{-\pi^2/8r^2}.
        \end{equation*}
    The process $\zeta^{OU}$ has a probability measure that is absolutely continuous with respect to $\P^W$, we write
        \begin{equation*}
            \forall x\in\mathcal{F}, \quad \P^W[\norm{x - \zeta^{OU}}_{\infty} \leq r] \sim C_{x} \mathrm{e}^{-\pi^2/8r^2}.
        \end{equation*}
    For Ornstein-Uhlenbeck process, we choose
        \begin{equation*}
            \phi^{OU}(r) \sim \mathrm{e}^{-\pi^2/8r^2}.
        \end{equation*}

Since we deal with sequences exhibiting weak dependency, let us formally define the mixing coefficient considered in this paper.

\subsection{Mixing condition}

The degree of dependence between observations of a stochastic process as they become distant apart in time is measured using mixing coefficients. Mixing processes were introduced to generalize the law of large numbers for non-i.i.d. stochastic processes. For effective modeling and inference, selecting the appropriate mixing condition for a stochastic process is crucial \citep{Peligrad2002, DedeckerPrieur2005, Rio2017}. One of the mixing criteria usually considered is $\beta$-mixing. It has been applied to demonstrate moment inequalities and central limit theorems \citep{ Dedeckeretal2007, Bosq2012, Poinas2020}.

\begin{definition} \label{def: mixing}
    Let $(\Omega, \mathcal{A}, \mathds{P})$ be a probability space, $\mathcal{B}$ and $\mathcal{C}$ be subfields of $\mathcal{A}$, and set $\beta(\mathcal{B},\mathcal{C}) = \E [\sup_{C\in\mathcal{C}}|\P(C) - \P(C|\mathcal{B})|]$.  
    For any array $\{Z_{t,T}:1\leq t \leq T\}$, define the coefficient
    \begin{align*}  
        \beta(k) &= \sup_{1\leq t\leq T-k} \beta\big( \sigma(Z_{s,T}, 1\leq s\leq t), \sigma(Z_{s,T}, t+k\leq s \leq T)\big ),
    \end{align*}
    where $\sigma(Z)$ denotes the $\sigma$-algebra generated by $Z$. The array $\{Z_{t,T}\}$ is said to be $\beta$-mixing or absolutely regular mixing if $\beta(k)\rightarrow 0$ as $k \rightarrow \infty.$
\end{definition}

This definition implies that if a process is $\beta$-mixing, asymptotic independence can be attained when $k \rightarrow \infty$. It is a ``just right'' assumption in analyzing weakly dependent sequences \citep{Vidyasagar1997}. There are different forms of $\beta$-mixing, such as exponentially $\beta$-mixing $\beta(k) = \bigO \big(e^{-\gamma k}\big)$, for $\gamma>0$, and arithmetically $\beta$-mixing $\beta(k) = \bigO \big(k^{-\gamma}\big)$ \citep{Ferraty&Vieu2006}. Numerous common time series models, such as autoregressive moving average (ARMA) models \citep{Mokkadem1998}, generalized autoregressive conditional heteroscedastic (GARCH) models \citep{CarrascoChen2002}, and some Markov processes \citep{Doukhan1994}, are known to be $\beta$-mixing. 

\section{Nadaraya-Watson estimation with Wasserstein distance} \label{sec:theoretical_guarantees}

We denote the conditional probability distribution of $Y_{t,T}|X_{t,T}=x$ by $\pi_t^\star(\cdot|x)$ and its conditional CDF by $F_t^\star(\cdot|x)$, for a fixed $t\in \{1, \ldots, T\}$ and $x\in\mathscr{H}$. The mean conditional function reads as
\begin{equation*}
m^\star(\frac{t}{T}, x) = \E_{\pi_t^\star(\cdot|x)}[Y_{t,T}|X_{t,T}=x] = \int_{-\infty}^\infty y \,\diff\pi_t^\star(y|x).
\end{equation*}
Setting $K_1, K_2$ two $1$-dimensional basic kernel functions and $h$ a bandwidth that depends on the sample size $T$, i.e., $h=h(T)$ with $h(T)\rightarrow 0$ as $T\rightarrow\infty$. 
For ease of notation, we set the scaled kernels $K_{h,i}(\cdot) = K_i(\frac{\cdot}{h}),$ for $ i=1,2$. Next, we define the considered NW estimator.

\begin{definition}
\label{definition: pi_hat}
The NW estimator of $\pi_t^\star(\cdot|x)$ is given by 
\begin{equation*}
    \hat{\pi}_t(\cdot|x) = \sum_{a=1}^T\omega_{a}(\frac t T,x)\delta_{Y_{a,T}},
\end{equation*}
where
\begin{align}\label{def: weights}
  \omega_{a}(\frac{t}{T}, x)=\frac{\displaystyle  K_{h,1}(\frac{t}{T} - \frac{a}{T})K_{h,2}(\mathsf{D} (x,X_{a,T}))}{\displaystyle \sum_{a=1}^TK_{h,1}(\frac{t}{T} - \frac{a}{T})K_{h,2}(\mathsf{D} (x,X_{a,T}))}.
 \end{align}
NW estimator of the conditional CDF $F^\star_t(y|x)$ can be written as, for all $y \in \R$, 
\begin{equation}\label{eqn:CDF of pi-hat}
    \hat{F}_t(y|x)=\sum_{a=1}^T\omega_{a}(\frac{t}{T},x)\mathds{1}_{Y_{a,T}\leq y}. \end{equation}
\end{definition}
This definition extends the estimator considered in \cite{Tinioetal1.2024} to a functional covariate $X_{t,T}$. The weights $\{\omega_a(u,x)\}_{a= 1, \ldots, T}$ are assumed to be measurable functions of $x$, $X_{a,T}$, and $u$ but do not depend on $Y_{a,T}$. Note that in \cite{Daisuke2022}, the NW estimator of $m^\star(u,x)$ is given by 
\begin{equation}\label{eqn: mhat definition}
    \hat m(u,x) = \sum_{a=1}^T   \omega_{a}(u, x) Y_{a,T}.
\end{equation} 
\paragraph{Remark.} We are using two kernel functions: one with respect to the rescaled time $u=\frac{t}{T}$ and the other in the direction of the functional $X_{t,T}$. To appropriately assign weights $\omega_a(\frac{t}{T},x)$, we smooth with respect to the rescaled time and the space-direction of the covariates $X_{t,T}$ to analyze the local behavior of the data \citep{vogt2012}. We consider a single bandwidth $h$ for the kernels $K_{h,i}(\cdot)$; however, $h$ could also be different for $K_{h,1}(\cdot)$ and $K_{h,2}(\cdot)$ \citep{Silverman1998}. \\

Next, let us establish the underlying assumptions used for our main results.

\subsection{Assumptions}
 
The following assumptions are conventional in the literature of LSTS~\citep{vogt2012, Daisuke2022, SoukariehSISP} and CDE~\citep{Halletal1999, metmous2021nonparametric, MR4783436}. 

\begin{assumption}[Local stationarity]
\label{Assumption: X is lsp}
Assume that the $\mathscr{H}$-valued process $\{X_{t,T}\}_{t=1, \ldots, T}$ 
is locally stationary approximated by $\{X_t(u)\}_{t=1, \ldots, T}$ for each time point $u\in [0,1]$. 
\end{assumption}

\begin{assumption}[Kernel functions]
\label{Assumption: kernel functions}
    The kernel $K_{1}(\cdot)$ is symmetric about zero, bounded and has compact support, that is, $K_1(v)=0$ for all $|v|>C_1$ for some $C_1<\infty$. On the other hand, the kernel $K_{2}(\cdot)$ is bounded, and has a compact support in $[0,1]$ such that $0<K_{2}(0)$ and $K_{2}(1)=0$.     In addition, $K_{2}'(v) = \diff K_{2}(v)/\diff v$ exists on $[0,1]$, satisfying $C_1' \leq K_{2}'(v) \leq C_2'$, for real constants $-\infty < C_1' < C_2' <0$. Moreover, for $i=1,2$, $K_i(\cdot)$, is Lipschitz continuous, that is, $|K_i(v)-K_i(v')|\leq L_i|v-v'|$ for some $L_i<\infty$ and all $v,v'\in \R$. We further assume the following:
    \begin{equation}\label{eqn: some properties of K}
        \int K_{i}(z)\diff z = 1,  
        \quad \text{and}
        \quad \int zK_{1}(z)\diff z = 0. 
            \end{equation}
\end{assumption}

Assumption \ref{Assumption: X is lsp} formalizes the property of the functional covariate $X_{t, T}$ as locally stationary. Assumption \ref{Assumption: kernel functions} is standard in literature. The conditions that $K_{i}$ is compactly supported and Lipschitz implies that the kernel function has a bounded rate of change and is essential in obtaining upper bounds. First condition in (\ref{eqn: some properties of K}) is a normalization, ensuring that the kernel can be interpreted as a probability density function. We assume that $K_2(\cdot)$ is compactly supported in $[0,1]$; that is, it is a kernel of type II \citep{Ferraty&Vieu2006}. Second condition implies that $K_1(\cdot)$ is symmetric around the origin, and it ensures that it does not introduce first-order linear bias when applied to the data.  

\begin{assumption}[Small ball probability]
\label{Assumption: small ball}
 Let $B(x, h) = \{v\in \mathscr{H}: \mathsf{D} (x,v)\leq h\}$ denote a ball centered at $x\in \mathscr{H}$ with radius $h$. We assume that for all $u \in [0,1]$, $x\in \mathscr{H}$, and $h>0$, there exists positive constants $C' < C$, such that
    \begin{equation}\label{eqn: CDF bounded by small ball}
        0 < C' \phi(h)\psi(x) \leq \P[X_t(u)\in B(x,h)] =: F_u(h;x) \leq C \phi(h)\psi(x),
    \end{equation}
    where $\phi(0)\rightarrow 0$ and $\phi(u)$ is absolutely continuous in a neighborhood of the origin, and $\psi(x)$ is a nonnegative functional in $x\in \mathscr{H}$. 
\end{assumption}

\begin{assumption}[Regularity condition on $h$ and $\phi(h)$] \label{Assumption: bandwidth}
    Assume that as $T\rightarrow \infty$, the bandwidth $h$ satisfies $T^\frac{1}{2} h\phi(h) \rightarrow \infty$.
\end{assumption}

Assumption \ref{Assumption: small ball} gives condition on the distributional behavior of the variables. Equation (\ref{eqn: CDF bounded by small ball}) controls the behavior of the small ball probability around zero. The small ball probability can be approximately expressed as the product of two independent functions $\phi(\cdot)$ and $\psi(\cdot)$. This condition corresponds to the assumption used in \cite{Gasseretal1998, Daisuke2022, Bouzebda20234}. On the other hand, Assumption \ref{Assumption: bandwidth} indicates that $h$ should converge slower to zero, for instance, at a polynomial rate, i.e., $h= \bigO(T^{-\xi})$, for small $\xi>0$. As $h$ approaches zero, $\phi(h)$ also goes to zero. Assumption \ref{Assumption: bandwidth} is a strengthening of the condition in \cite{Daisuke2022} that $Th\phi(h) \rightarrow\infty$ and is needed to guarantee our resulting convergence rates. 
With this, for Fractal-type processes, Assumption \ref{Assumption: bandwidth} holds true when we choose $h \sim T^{-\xi}$ for $0 < \xi < \frac{1}{2(1+\tau_0)}$ and $\phi(h) \sim h^{\tau_0}$ for some $ \tau_0 > 1$ \citep{Bouzebda20234, AguaBouzebda2024}. To see different expressions of the function $\phi(h)$, one may refer to \cite{Ferraty&Vieu2006, Bogachev1998} for some discussions on fractal-type processes, \cite{MayerWolf1993} for diffusion processes, and \cite{LiShao2001} for general Gaussian processes. We have given examples of the forms of $\phi(\cdot)$ in Subsection \ref{subsection: small ball probability}.

\begin{assumption}[Conditional CDF] \label{assumption: CDF}
    The conditional CDF $F_{\cdot}^\star(\cdot|\cdot)$ is Lipschitzian, i.e., $\big| F_a^\star(\cdot|x) - F_t^\star (\cdot| x') \big| \leq L_{F^\star} \big( \mathsf{D} (x, x') + \big|\frac{a}{T} - \frac{t}{T}\big|\big)$, for some $L_{F^\star}<\infty$, and for all $a, t \in \{1, \ldots, T\},$ $x, x'\in \mathscr{H}$. 
\end{assumption}

The conditional CDF $F_{\cdot}^\star(\cdot|\cdot)$ should behave smoothly and not change much as the observation does, as assumed in \cite{Bouananietal2018, metmous2021nonparametric, Tinioetal1.2024}. We do not assume that the conditional CDF is twice differentiable, in contrast to~\cite{Halletal1999, ferraty2005conditional, OtneimTjostheim2016}.

\begin{assumption}[Mixing condition] \label{Assumption: mixing}
    The process $\{(X_{t,T}, \varepsilon_{t,T})\}$ is arithmetically $\beta$-mixing satisfying $\beta(k) \leq Ak^{-\gamma}$ for some $A>0$ and $\gamma>2$. We also assume that for some $p>2$ and $\zeta>1-\frac{2}{p}$,
    \begin{align}\label{eqn: infinite sum of betas is finite}
        \sum_{k=1}^{\infty} k^\zeta \beta(k)^{1-\frac{2}{p}} <\infty. 
    \end{align}
\end{assumption}

\begin{assumption}[Blocking condition]\label{assumption: blocking}
    There exists a sequence of positive integers $\{q_T\}$ satisfying $q_T \rightarrow \infty$ and $q_T = o\big(\sqrt{Th\phi(h)}\big)$, as $T \rightarrow \infty$.
\end{assumption}

For dependent sequences estimation technique, Assumptions \ref{Assumption: mixing} and \ref{assumption: blocking} are helpful. A more robust type of independence between far-off observations in a process is the $\beta$-mixing \citep{Bradley2005, Rio2017, Poinas2020}. The decay of the regular mixing coefficient $\beta(k)$ is highlighted by condition (\ref{eqn: infinite sum of betas is finite}). 
Bernstein's blocking approach was utilized to create independent blocks in the proof of Theorem \ref{Theorem: convergence of EW1} \citep{Bernstein1927}. Assumption \ref{assumption: blocking} defines the big block size as proportional to $q_T$.

\subsection{Convergence in Wasserstein distance}

We establish the convergence rate of NW estimator $\hat{\pi}_t(\cdot|x)$ wrt the Wasserstein distance. Theorem \ref{Theorem: convergence of EW1} below generalizes the convergence results in \cite{Tinioetal1.2024} to the functional or infinite-dimensional setting.

\begin{theorem}\label{Theorem: convergence of EW1}
Suppose Assumptions \ref{Assumption: X is lsp} - \ref{assumption: blocking} are satisfied and define $I_h = [C_1h, 1 - C_1h]$. Then
\begin{align*}
    \sup_{x\in \mathscr{H}, \frac{t}{T}\in I_h} \E[W_1(\hat{\pi}_t(\cdot|x),\pi_t^\star(\cdot|x))]
    &= \bigO_\P\Big(\frac{1}{T^\frac{1}{2} h\phi(h)} + h\Big).
\end{align*}
\end{theorem}

This convergence rate is comparable with Theorem 1 in \cite{Tinioetal1.2024} for the $d$-dimensional covariate case. However, we do not have the bias term involving $\rho$ that comes from approximating $X_{t, T}$ by a locally stationary $X_t(\frac{t}{T})$. This rate depends on the bandwidth $h$ and the small ball probability $\phi(h)$, where $\phi(h)\rightarrow 0$ as $h\rightarrow 0$, as highlighted in Assumption \ref{Assumption: bandwidth}. We defer the proof to Appendix \ref{appendix: proof of convergence of EW1}, which follows similar steps of the proof of Theorem 1 in \cite{Tinioetal1.2024}.

For the i.i.d. case $(Y_t, X_t)_{t=1,\ldots, T}$, where $Y_t$ is scalar, and $X_t$ is functional, \cite{metmous2021nonparametric} provided a convergence result for their proposed conditional CDF estimator with surrogate data. The proposed estimator involves three rescaled kernel functions: a one-dimensional kernel $K(\frac{\cdot}{h_K})$ to account the functional $X_t$, an integrated Kernel $H(\frac{\cdot}{h_H})$ that acts as a CDF of $Y_t$, and a two-dimensional kernel $W(\frac{\cdot}{a_T}, \frac{\cdot}{a_T})$ to account for the surrogate variable, with bandwidths $h_K$, $h_H$, and $a_T$, respectively. This estimator converges to the true conditional distribution of order $\bigO(h_K^{c_1} + h_H^{c_2} + a_T^{c_1}) + \bigO \bigg(\sqrt{\frac{\log d_T}{T\phi(a_T)}} \bigg) + \bigO \bigg(\sqrt{\frac{\log T}{T\phi(h_K)}} \bigg)$, where $c_1, c_2 > 0$ and $d_T$ satisfies $\frac{\log^2 T}{T\phi(a_T)} \leq d_T \leq \frac{T\phi(a_T)}{\log T}$. If there is no surrogate data and the integrated kernel $H$ is replaced by an indicator function, this convergence rate becomes $\bigO(h_K^{c_1}) + \bigO \bigg(\sqrt{\frac{\log T}{T\phi(h_K)}} \bigg)$, which is comparable to the result above.

\begin{corollary}\label{Remark: bound EW_s-s}
    Suppose Assumptions \ref{Assumption: X is lsp} - \ref{assumption: blocking} are satisfied and $Y_{t, T}$ is uniformly bounded by $M>0$. Then, for $r\geq 1$,
    \begin{align*}    \sup_{x\in \mathscr{H}, \frac{t}{T}\in I_h} \E[W_r^r(\hat{\pi}_t(\cdot|x),\pi_t^\star(\cdot|x))]
    &= \bigO_\P\Big(\frac{1}{T^\frac{1}{2} h\phi(h)} + h\Big).
    \end{align*}
\end{corollary}

Proof of Corollary \ref{Remark: bound EW_s-s} is shown in Appendix \ref{appendix: proof of EW_s-s}. For the i.i.d case, \cite{bobkovledoux2019} (Theorem 5.3) showed that $\E[W_r^r(\mu_T, \mu)] = \bigO((T+2)^{-\frac{r}{2}})$, for $r\geq 1$, where $\mu_T$ is the empirical measure of an i.i.d sample $(X_t)_{t\geq 1}$ with common law $\mu$. 

\begin{corollary}\label{corollary: convergence of the 2nd moment}
Let Assumptions \ref{Assumption: X is lsp} - \ref{assumption: blocking} be satisfied. Then
\begin{align*}
    \sup_{x\in \mathscr{H}, \frac{t}{T}\in I_h} \norm{W_1\big(\hat{\pi}_t(\cdot|x), \pi_t^\star(\cdot|x)\big)}_{L_2}
    &= \bigO_\P\Big(\frac{1}{T^\frac{1}{2} h\phi(h)} + h\Big).
\end{align*}

\end{corollary}

Proof of Corollary \ref{corollary: convergence of the 2nd moment} is based on Minkowski's integral inequality: for any $r\geq 1$,
\begin{align}\label{eqn: minkowski's inequality}
    {\Big\|\int \big|\hat{F}_t(y|x)-F^\star_{t}(y|x)\big|\diff y\Big\|}_{L_r} \leq \int {\big\|\hat{F}_t(y|x)-F^\star_{t}(y|x)\big\|}_{L_r} \diff y.
\end{align}
The remainder of the proof adheres to the same lines used for Theorem \ref{Theorem: convergence of EW1}, refer to Appendix \ref{appendix: proof of convergence of the 2nd moment}. The following proposition shows that the NW conditional mean function estimator $\hat{m}$ warrants

\begin{proposition}\label{prop: |mhat-m| leq W1}
    Let Assumptions \ref{Assumption: X is lsp} - \ref{assumption: blocking} be satisfied and $\hat{m}(\frac{t}{T},x)$ be defined by (\ref{eqn: mhat definition}). Then
    \begin{align*}
        \sup_{x\in \mathscr{H}, \frac{t}{T}\in I_h} \E \big[|\hat{m}(\frac{t}{T},x) - m^\star(\frac{t}{T}, x)| \big]
        &= \bigO_\P\Big(\frac{1}{T^\frac{1}{2} h\phi(h)} + h\Big).
    \end{align*}
\end{proposition}

Proof of Proposition \ref{prop: |mhat-m| leq W1} is detailed in Appendix \ref{appendix: proof of mhat leq W1}. Similar to Proposition 1 in \cite{Tinioetal1.2024}, this result indicates that Wasserstein distance can be used to obtain the convergence rate of $\hat{m}(u, x)$. The bound of the Wasserstein distance is slower than $\hat{m}(u,x)$ since we are examining differences between distributions, not just differences between conditional means \citep{Tinioetal1.2024}. This rate is comparable to Theorem 3.1 in \cite{Daisuke2022} with convergence rate of order $\bigO_\P \big(\sqrt{\frac{\log T}{Th\phi(h)}} + h^{2\wedge \beta} \big)$. We remark that a similar component for the bias term can be obtained if we assume that $F_\cdot^\star(\cdot)$ is twice differentiable and satisfies the Hölder condition.

\begin{proposition}\label{corollary: convergence of EW1 chosen h}
Suppose $X_t(u)$ is a fractal-type process and Assumptions \ref{Assumption: X is lsp} - \ref{assumption: blocking} are satisfied. Let the bandwidth be chosen to be $h = \bigO(T^{-\xi})$, and the small ball probability take the form $\phi(h)=h^{\tau_0}$, where $0 < \xi < \frac{1}{2(1+\tau_0)}$ and $\tau_0 >1$. Then
\begin{align*}
    \sup_{x\in \mathscr{H}, \frac{t}{T}\in I_h} \E[W_1(\hat{\pi}_t(\cdot|x),\pi_t^\star(\cdot|x))]
    &= \bigO_\P\Big(\frac{1}{T^{\frac{1}{2}-\xi(1+\tau_0)}} + \frac{1}{T^{\xi}}\Big).
\end{align*}
\end{proposition}

As demonstrated in Appendix \ref{appendix: proof of convergence of EW1 chosen h}, by setting $h = \bigO(T^{-\xi})$ and $\phi(h)=h^{\tau_0}$, proof of Proposition \ref{corollary: convergence of EW1 chosen h} follows immediately from the proof of Theorem \ref{Theorem: convergence of EW1}.

\subsection{Bandwidth selection criterion}
\label{sec:bandwidth selection}

In nonparametric kernel estimation, especially NW, the bandwidth must be suitably selected for the estimator to perform well. Bandwidth selection methods have already been established and developed in \cite{RachdiVieu2007,  MR4507275}. This paper considers the leave-one-out cross-validation procedure used in \cite{Benhennietal2007, RachdiVieu2007}. For any fixed $i \in \{1, \ldots, T\}$, we define
\begin{equation}\label{eqn: leave-one-out conditional CDF}
    \hat{m}_{i}(\frac{t}{T}, x)=\sum_{a=1; a \neq i}^T\omega_{a}(\frac{t}{T},x)Y_{a, T},
\end{equation}
where $\omega_{a}(\frac{t}{T},x)$ is given by (\ref{def: weights}). Equation (\ref{eqn: leave-one-out conditional CDF}) is regarded as the leave-out-$(X_{i,T}, Y_{i,T})$ estimator of $m_i^\star(\frac{t}{T},x)$. To minimize the quadratic loss function, we introduce the following criterion
\begin{equation}\label{eqn: CV w tilde}
    CV(y, x, h) := \frac{1}{T} \sum_{i = 1}^T \big( Y_{i,T} -  \hat{m}_{i}(\frac{t}{T}, x) \big)^2 \widetilde{g}(X_{i,T}),
\end{equation}
for some non-negative weight function $\widetilde{g}(\cdot)$. As highlighted in \cite{RachdiVieu2007}, we choose a bandwidth $\hat{h}$ among $h\in [a_T, b_T]$ that minimizes (\ref{eqn: CV w tilde}). For bandwidths that are locally chosen by data-driven method, according to \cite{Benhennietal2007}, we replace (\ref{eqn: CV w tilde}) by 
\begin{equation*}    CV(y, x, h) := \frac{1}{T} \sum_{i = 1}^T \big( Y_{i,T} -  \hat{m}_{i}(\frac{t}{T}, x)  \big)^2 \hat{g}(X_{i,T}).
\end{equation*}
In practice, for $i\in \{1, \ldots, T\}$, we take the uniform global weights $\widetilde{g}(X_{i,T}) = 1$, or the local weights 
\begin{equation*}
    \hat{g}(X_{i,T}, x) = 
    \begin{cases}
        1 & \text{if $\mathsf{D} (X_{i,T}, x)\leq h$},\\
        0 & \text{otherwise}.
    \end{cases}
\end{equation*}

\section{Numerical experiments}
\label{sec:numerical_experiments}

To illustrate the convergence of NW estimator wrt Wasserstein distance, we conduct numerical experiments using synthetic and real-world datasets. 

\subsection{Synthetic data}

We generate samples $(X_{t,T}, Y_{t,T})_{t=1,\ldots,T}$ using examples provided in \cite{vanDelft2018}. We consider two locally stationary processes.

\paragraph{Generation of functional covariates.}
We generate the functional covariate from a Hilbert space $\mathscr{H} = L_{\R}^2([0,1])$, using the following examples:\\

\begin{figure}[htbp]
\centering
\begin{subfigure}[b]{0.48\textwidth} 
    \centering
    \includegraphics[width=\textwidth]{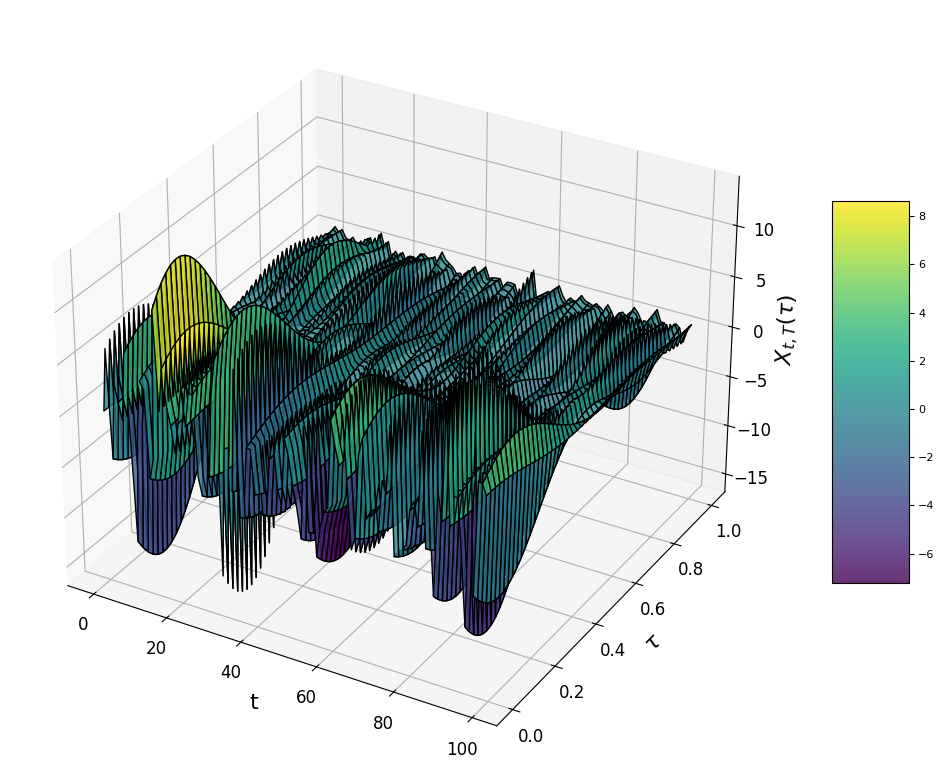}
    \caption{$X_{t,T}(\tau)$ for all $t$ and some $\tau$}
    \label{fig:covariates_gaussian_tvFAR1_3d}
\end{subfigure}
\hfill
\begin{subfigure}[b]{0.48\textwidth}
    \centering
    \begin{subfigure}[b]{\textwidth} 
        \includegraphics[width=.9\textwidth]{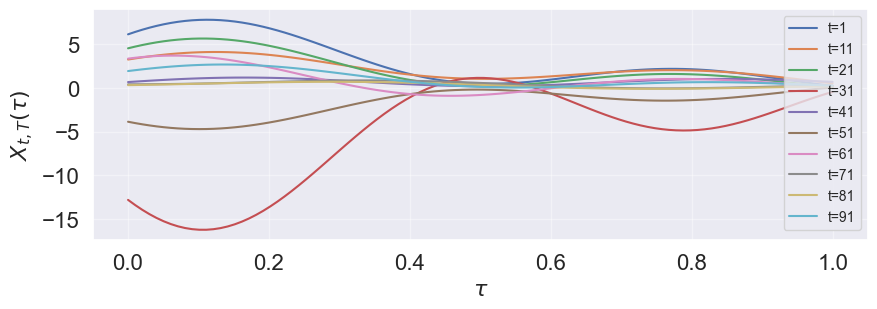}
        \caption{$X_{t,T}(\tau)$ at $\tau$ given some $t$}
        \label{fig:covariates_gaussian_tvFAR1_different_t}
    \end{subfigure}
    \vspace{2mm} 
    \begin{subfigure}[b]{\textwidth} 
        \includegraphics[width=.9\textwidth]{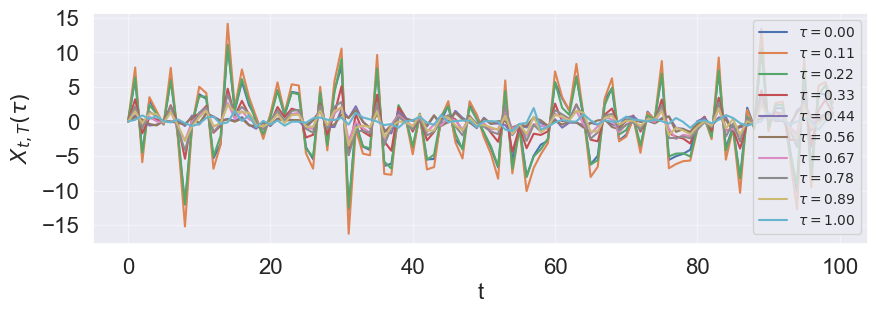}
        \caption{$X_{t,T}(\tau)$ at $t$ given some $\tau$}
        \label{fig:covariates_gaussian_tvFAR1_different_tau}
    \end{subfigure}
\end{subfigure}
\caption{Realizations of Gaussian tvFAR(1) $X_{t,T}(\tau)$ for all $t$ and some $\tau$ for $T=100$ with $J=7$ and $N=100$ discretization points of $\tau\in[0,1]$.}
\label{fig: covariates Gaussian tvFAR(1)}
\end{figure}

\noindent{\scshape Example 1. Gaussian tvFAR(1).} 
We consider the time-varying functional autoregressive process of order 1, tvFAR(1), with Gaussian noise represented by
\begin{equation}\label{eqn: tvFAR1 process}
    X_{t,T}(\tau) = B_{t/T}(X_{t-1, T})(\tau) + \eta_t(\tau), \quad \tau \in [0,1], \quad t=1,\ldots,T,
\end{equation}
with a linear operator $B_{t/T}$ indexed by rescaled time $u=\frac{t}{T}$ and innovation function $\eta_t$. The innovation $\eta_t$ is a linear combination of the Fourier basis function $(\psi_j)_{j\in\mathds{N}}$ with coefficients $\langle \eta_t, \psi_j \rangle$ that are generated from independent zero-mean Gaussian distribution with $j$th coefficient having variance $(\pi(j-1.5))^{-2}$, that is,
\begin{equation*}
    \eta_t = \sum_{j=1}^\infty \langle \eta_t, \psi_j \rangle \psi_j \quad \text{with } \langle \eta_t, \psi_j \rangle \sim \mathcal{N}\big(0, (\pi(j-1.5))^{-2}\big),
\end{equation*}
where 
\begin{equation*}
    \psi_j(\tau) = 
    \begin{cases}
        \sqrt{2} \sin(\pi j\tau), & \text{ for odd } j,\\
        \sqrt{2} \cos(\pi j\tau), & \text{ for even } j.
    \end{cases}
\end{equation*} 
In application, we truncate an infinite-dimensional series at some $J$ basis functions. Now, instead of decomposing $X_{t,T}$ on the basis $(\psi_j)_{j\in\mathds{N}}$ by $\sum_{j=1}^\infty \langle X_{t,T}, \psi_j \rangle \psi_j $, it can be represented by an approximate finite-dimensional $X_{t,T}$:
$$X_{t, T} = \sum_{j=1}^J \langle X_{t,T}, \psi_j \rangle \psi_j.$$
Hence, $\boldsymbol{X}_{t,T} \approx \boldsymbol{B}_{t/T} \boldsymbol{X}_{t-1,T} + \boldsymbol{\eta}_t$, $t=1,\ldots,T$, where $\boldsymbol{X}_{t,T} = (\langle X_{t,T}, \psi_1 \rangle, \ldots, \langle X_{t,T}, \psi_J \rangle)'$, $\boldsymbol{\eta}_t = (\langle \eta_t, \psi_1 \rangle, \ldots, \langle \eta_t, \psi_J \rangle)'$, and $\boldsymbol{B}_{t/T} = (\langle B_{t/T}(\psi_i), \psi_j \rangle)_{1\leq i,j \leq J}$. In this example, the matrix $\boldsymbol{B}_{t/T}$ is defined as $\boldsymbol{B}_{t/T} = \frac{0.4 \boldsymbol{A}_{t/T}}{\|\boldsymbol{A}_{t/T}\|_\infty}$, where $\boldsymbol{A}_{t/T}$ is a $J\times J$ matrix with entries $A_{t/T}(i,j)$ that are mutually independent zero-mean Gaussian random variables with variance $\frac{t}{T i^6} + (1-\frac{t}{T})\e^{-j-i}$ and $\| A\|_\infty = \sup_{\|x\|\leq 1}\|A x\|$ is a Schatten $\infty$-norm. Figure \ref{fig: covariates Gaussian tvFAR(1)} shows the plot of $X_{t,T}(\tau)$ for $T=100$. This example was also used in \cite{AguaBouzebda2024}. \\

\begin{figure}[htbp]
\centering
\begin{subfigure}[b]{0.48\textwidth} 
    \centering
    \includegraphics[width=\textwidth]{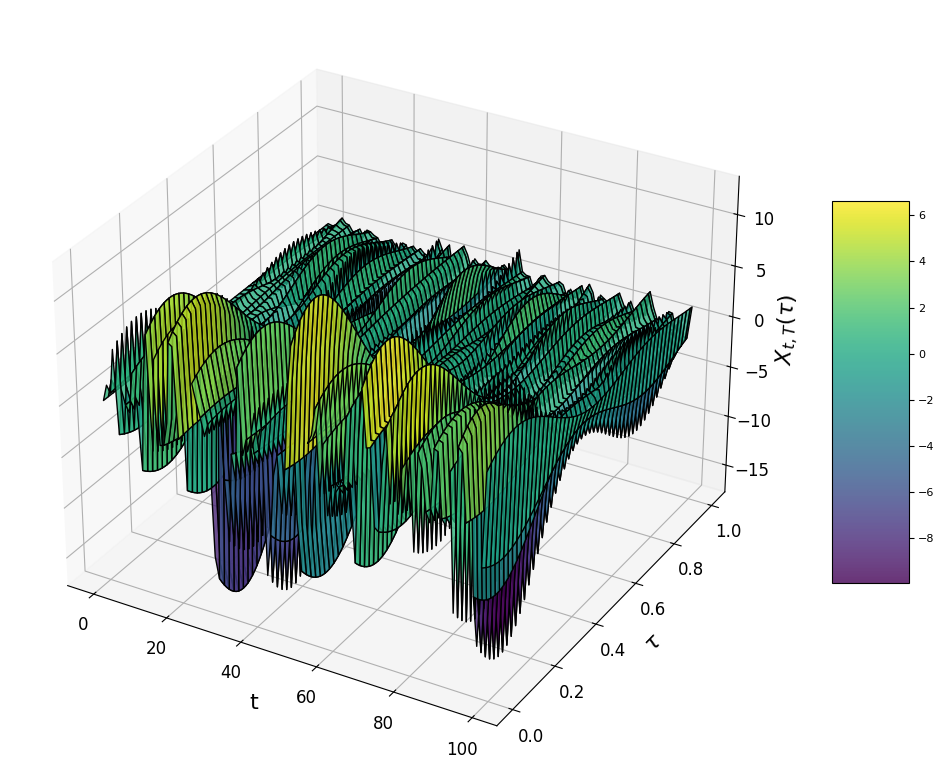}
    \caption{$X_{t,T}(\tau)$ for all $t$ and some $\tau$}
    \label{fig:covariates_gaussian_tvFAR2_3d}
\end{subfigure}
\hfill
\begin{subfigure}[b]{0.48\textwidth}
    \centering
    \begin{subfigure}[b]{\textwidth} 
        \includegraphics[width=.9\textwidth]{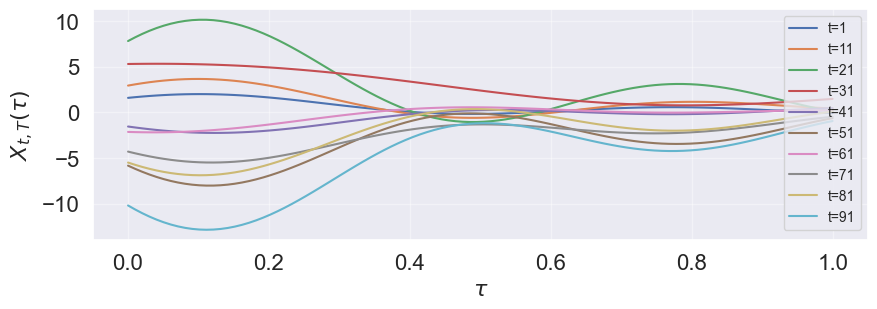}
        \caption{$X_{t,T}(\tau)$ at $\tau$ given some $t$}
        \label{fig:covariates_gaussian_tvFAR2_different_t}
    \end{subfigure}
    \vspace{2mm} 
    \begin{subfigure}[b]{\textwidth} 
        \includegraphics[width=.9\textwidth]{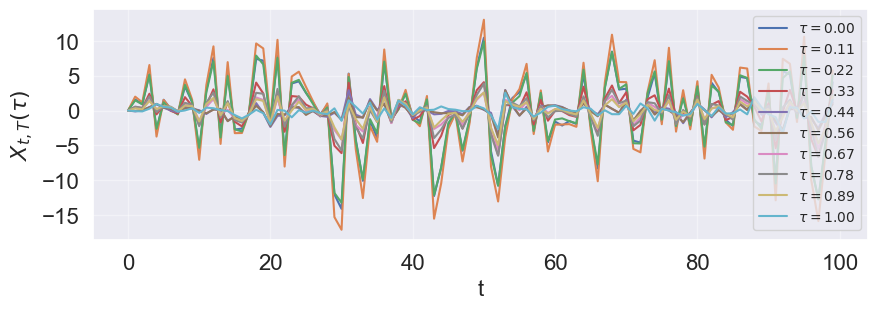}
        \caption{$X_{t,T}(\tau)$ at $t$ given some $\tau$}
        \label{fig:covariates_gaussian_tvFAR2_different_tau}
    \end{subfigure}
\end{subfigure}
\caption{Realizations of Gaussian tvFAR(2) $X_{t,T}(\tau)$ for all $t$ and some $\tau$ for $T=100$ with $J=7$ and $N=100$ discretization points of $\tau\in[0,1]$.}
\label{fig: covariates Gaussian tvFAR(2)}
\end{figure}

\noindent{\scshape Example 2. Gaussian tvFAR(2).}
We next consider the time-varying functional autoregressive process of order 2, tvFAR(2), with Gaussian noise defined by
\begin{equation}\label{eqn: tvFAR2 process}
    X_{t,T}(\tau) = B_{t/T, 1}(X_{t-1, T})(\tau) + B_{t/T, 2}(X_{t-2, T})(\tau) + \eta_t(\tau), \quad \tau \in [0,1], \quad t=1,\ldots,T,
\end{equation}
where $B_{t/T, 1}$ and $B_{t/T, 2}$ are linear operators indexed by the rescaled time $u=\frac{t}{T}$ and innovation function $\eta_t$ is a linear combination of the Fourier basis function $(\psi_j)_{j\in\mathds{N}}$. The parameters are set similarly to (\ref{eqn: tvFAR1 process}) with $\boldsymbol{B}_{t/T,1} = \frac{0.4\cos(1.5 - \cos(\pi \frac{t}{T})) \boldsymbol{A}_{t/T,1}}{\|\boldsymbol{A}_{t/T,1}\|_\infty}$ and $\boldsymbol{B}_{t/T,2} = \frac{-0.5 \boldsymbol{A}_{t/T,2}}{\|\boldsymbol{A}_{t/T,2}\|_\infty}$, where $A_{t/T,1}(i,j)$ and $A_{t/T,2}(i,j)$ are mutually independent-centered Gaussian random variables with variances $\e^{-(i-3)-(j-3)}$ and $1/(i^4 + j)$, respectively. Realizations of $X_{t,T}(\tau)$ in this example are depicted in Figure \ref{fig: covariates Gaussian tvFAR(2)}.\\

\begin{figure}[ht]
\centering
\begin{subfigure}[b]{.8\textwidth}
\centering
\includegraphics[height=.8in, width=.8\textwidth]{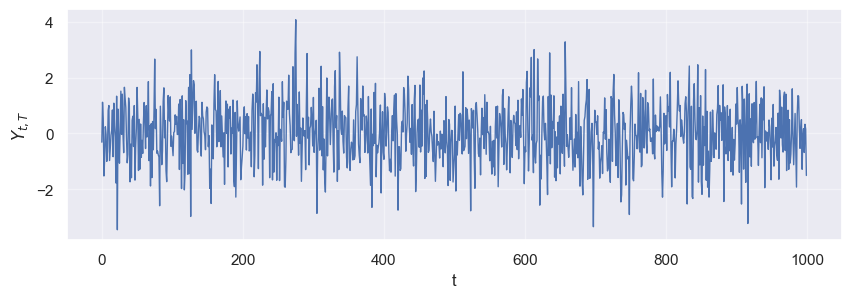} 
\caption{Using Gaussian tvFAR(1) $X_{t, T}$}
\label{fig:Y_gaussian_tvFAR1}
\end{subfigure}
\begin{subfigure}[b]{.8\textwidth}
\centering
\includegraphics[height=.8in, width=.8\textwidth]{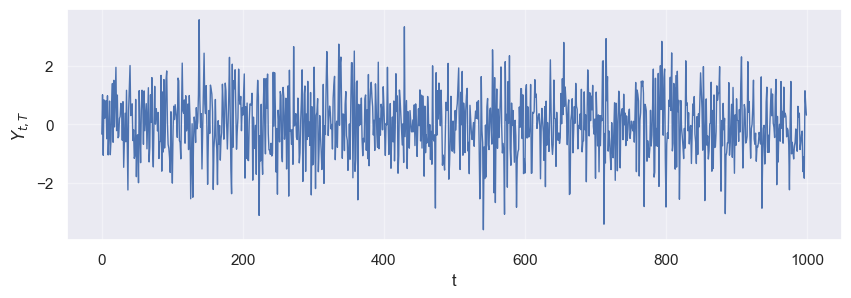} 
\caption{Using Gaussian tvFAR(2) $X_{t, T}$}
\label{fig:Y_gaussian_tvFAR2}
\end{subfigure}
\caption{Time plots of response variables for $T=1000$}
\label{fig: response}
\end{figure}

\paragraph{Generation of response variables.}
Using locally stationary covariates $X_{t, T}$ in Examples 1 and 2, the response variables $Y_{t, T}$ are generated by (\ref{eq:major_estimation_problem}) with $\varepsilon_{t, T} \sim \mathcal{N}(0,1)$ and 
$$m^\star(\frac{t}{T},x) = 2.5\sin(2\pi \frac{t}{T})\int_0^1 \cos(\pi x(\tau))\diff \tau.$$
Figure \ref{fig: response} shows the time plots of the responses for each process using $T=1000$, whose values remain tight with constant mean.

\paragraph{Monte Carlo simulations.}
Using an identical Monte Carlo simulation process in \cite{Tinioetal1.2024}, we calculate the NW estimator and true conditional probability distribution for a fixed time $t \in \{1, \ldots, T\}$. Each process is replicated using $L = 500$, and as described in Algorithm \ref{alg: simulated data} of \cite{Tinioetal1.2024}, for each $l\in\{1, \ldots, L\}$, we compute the NW conditional CDF at a given time $t$. We calculate the average NW and the empirical conditional CDFs using these $L$ replications. We then quantify the corresponding Wasserstein distance.

\begin{figure}[ht]
\centering
\begin{subfigure}[b]{0.49\textwidth}
\centering
\includegraphics[width=\textwidth]{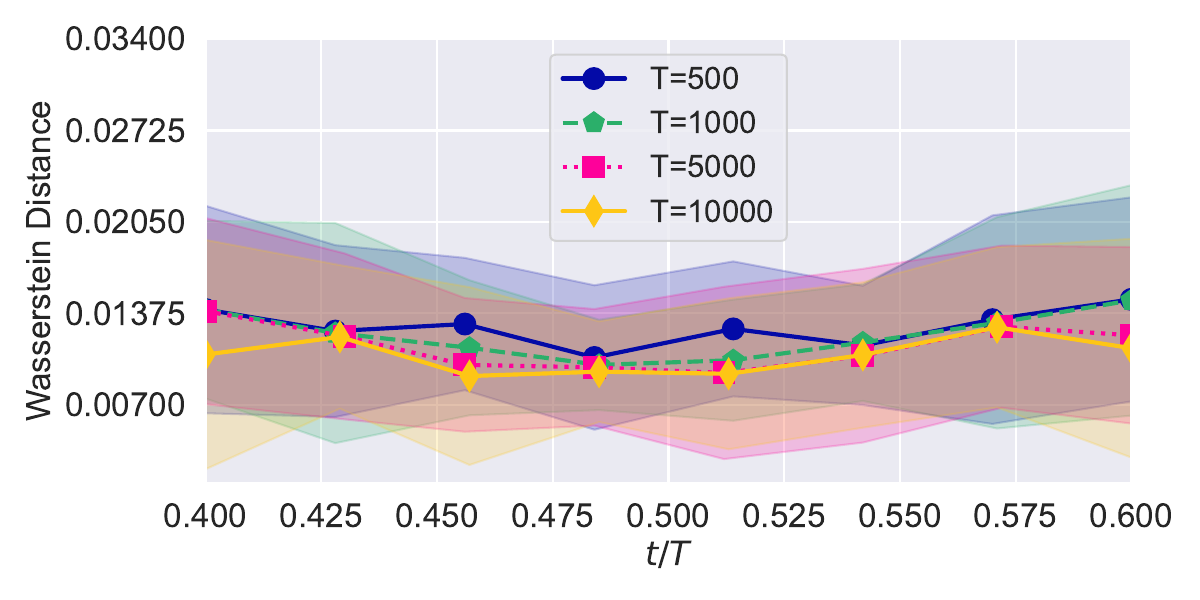} 
\caption{\centering Using Gaussian tvFAR(1) $X_{t, T}$; $K_1= \texttt{Uniform}, K_2 = \texttt{Silverman}$ }
\label{fig:convergence_gaussian_tvAR1}
\end{subfigure}
\begin{subfigure}[b]{0.49\textwidth}
\centering
\includegraphics[width=\textwidth]{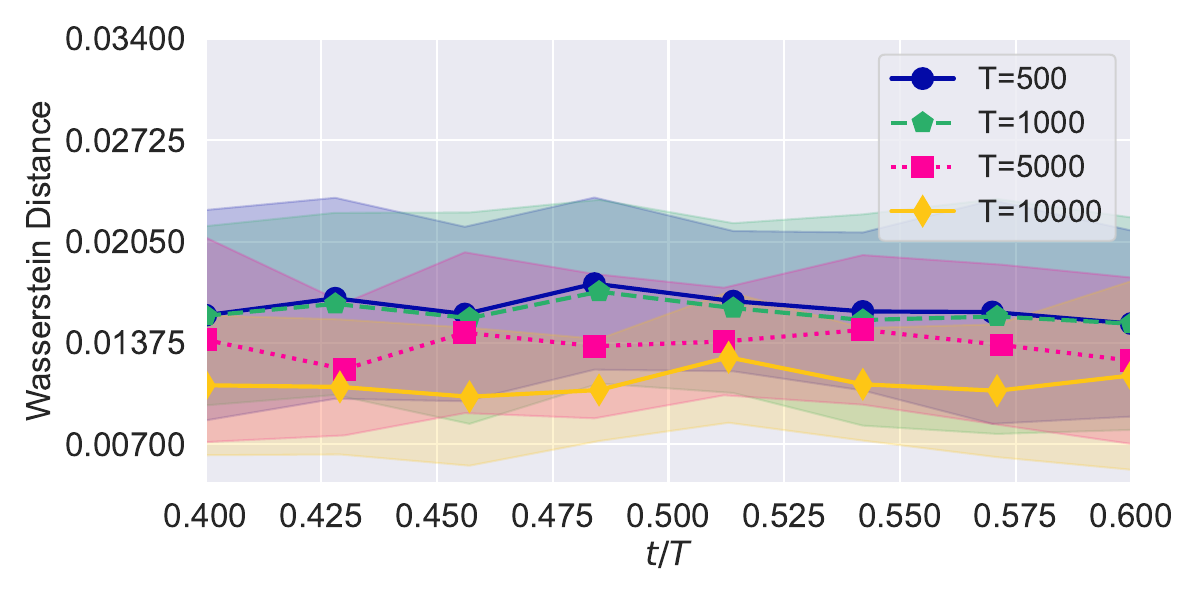} 
\caption{\centering Using Gaussian tvFAR(2) $X_{t, T}$; $K_1= \texttt{Tricube}, K_2 = \texttt{Gaussian}$ }
\label{fig:convergence_gaussian_tvAR2}
\end{subfigure}
\vspace{5pt}
\caption{Wasserstein distances $\pm$ standard deviation at different $u = \frac{t}{T}$ for $T = 500, 1000, 5000, 10000$ using $L=500$ replications and 50 Monte Carlo runs. }
\label{fig: convergence at various t increasing T}
\end{figure}

We obtain the expected Wasserstein distance between the underlying conditional distributions by conducting 50 Monte Carlo runs of Algorithm \ref{alg: simulated data}. To produce functional covariates, we select $N=100$ discretization points of $\tau \in [0,1]$ and set $J = 7$ since results do not vary much wrt $J$ \citep{AguaBouzebda2024}.  As specified in Figure \ref{fig: convergence at various t increasing T}, we use different kernels $K_1$ and $K_2$ for the chosen processes. Increasing sample sizes $T = 500, 1000, 5000, 10000$ are set. The bandwidths are chosen using the cross-validation method introduced in Section \ref{sec:bandwidth selection}. Our results are valid when $\frac{t}{T} \in I_h$, hence, we fix $t$ such that $\frac{t}{T} \in I_h = [C_1h, 1-C_1h]$ with constant $C_1 = 1$ for time kernels $K_1$ belonging to \texttt{Uniform} and \texttt{Tricube}.

Figure \ref{fig: convergence at various t increasing T} depicts the expected Wasserstein distances for each identified process. Wasserstein distance decreases as the sample size $T$ increases. This emphasizes that NW estimator captures the true distribution better as $T$ grows larger; it provides more representative distributions with reduced deviation from the true distribution. Remarkably, the largest sample size, $T=10000$, consistently achieves the minimum expected Wasserstein distance. This behavior is consistent across both processes under investigation. 

\subsection{Real-world data}    

We use two real-world datasets: sea surface temperature (SST) and Nikkie225. To handle these datasets, we employ the same method in \cite{ferraty2005conditional}, described below. \\

\begin{figure}[htbp]
\centering
\begin{subfigure}[b]{0.55\textwidth}
    \centering
    \includegraphics[width=\textwidth, height=1.2in]{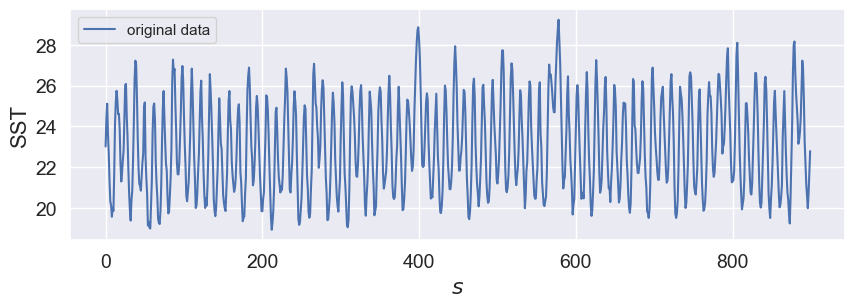} 
    \caption{Original time series $Z(s)$ ($n=900$)}
    \label{fig:SST_original}
\end{subfigure}
\begin{subfigure}[b]{0.47\textwidth}
    \centering
    \includegraphics[width=\textwidth, height=1.2in]{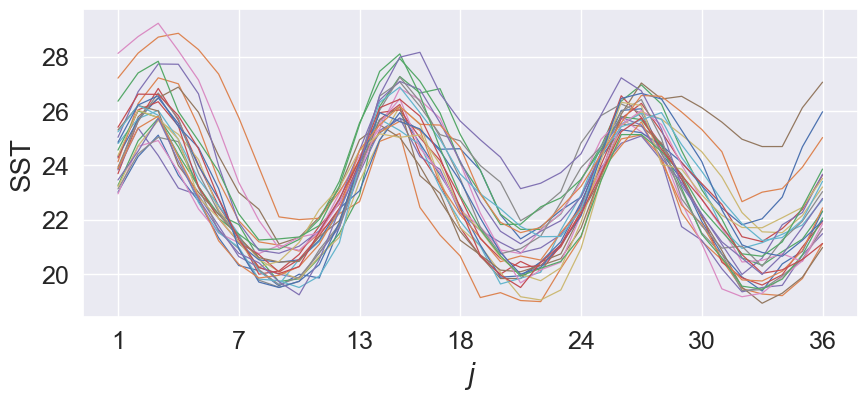} 
    \captionsetup{justification=centering} 
    \caption{Continuous sample curves $X_{t, T}(j)$;\\
    $t=1,\ldots, 25$}
    \label{fig:SST_covariate}
\end{subfigure}
\begin{subfigure}[b]{0.47\textwidth}
    \centering
    \includegraphics[width=\textwidth, height=1.2in]{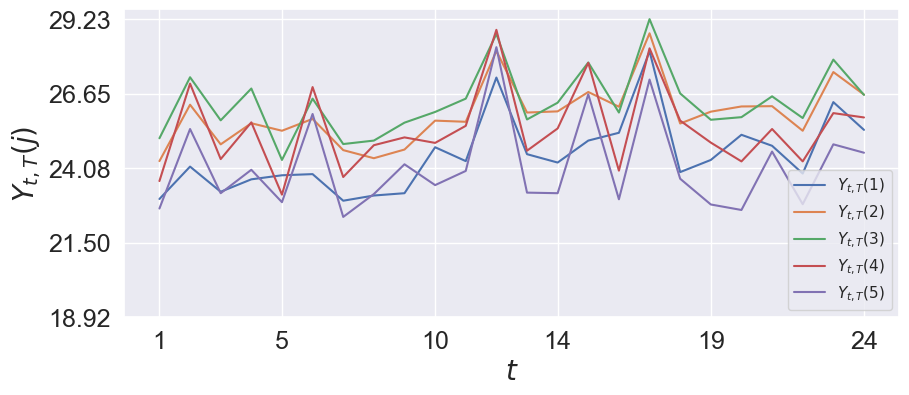} 
    \captionsetup{justification=centering} 
    \caption{Response $Y_{t, 24}(j)$;\\
    for some $j=1,\ldots,5$}
    \label{fig:SST_response}
\end{subfigure}
\caption{SST monthly time series from Jan. 1950 - Dec. 2024}
\label{fig:SST}
\end{figure}

\noindent{\scshape Example 3. SST data\footnote{Obtained from
\url{https://www.cpc.ncep.noaa.gov/data/indices/}}.} This dataset is used for climate monitoring and research, which is continuously updated by the National Centers for Environmental Information (NCEI) \citep{ kilpatrick2015decade}. We take the index from Niño 1+2 region with coordinates 0°- 10°South latitude and 90°West - 80°West longitude. This region covers the eastern equatorial Pacific near the coast of South America and is important for monitoring El Niño and La Niña events. SST contains 900 monthly data points from January 1950 to December 2024, depicted in Figure \ref{fig:SST_original}.

\paragraph{Constructing covariates.} To construct $X_{t, T}$, we treat the original series as 25 continuous sample curves, each containing 36 monthly observations as plotted in Figure \ref{fig:SST_covariate}. Particularly, we let the SST observed for $n=900$ months be $\{Z(s)\}_{s=1,\ldots,n}$, and build, $\forall j\in\{1, \ldots, 36\}$,
$$z_{t, T}(j) = Z(36(t-1)+j).$$
The covariates are then constructed as $X_{t, T} = (z_{t, T}(1),\ldots,z_{t, T}(36))$ that corresponds to the variations of SST for $t = 1, \ldots, 25$. We consider the original time series as 25 dependent functional covariates, $X_{1, 25}, \ldots, X_{25, 25}$, which are individually observed at 36 discretized points.

\paragraph{Constructing response.} We construct the response variables by 
$$Y_{t, T}(j) = z_{t+1, T}(j) = Z(36t + j),$$ 
for a fixed $j \in \{ 1,\ldots, 36\}$ and $t= 1,\ldots, 24$. This enables us to generate new 24 functional pairs $\{(X_{t, 24}, Y_{t, 24}(j))\}_{t=1,\ldots,24}$ \citep{ferraty2005conditional}. Figure \ref{fig:SST_response} presents sample plots of $Y_{t,24}(j)$ for some $j$.\\

\begin{figure}[htbp]
\centering
\begin{subfigure}[b]{0.55\textwidth}
    \centering
    \includegraphics[width=\textwidth, height=1.2in]{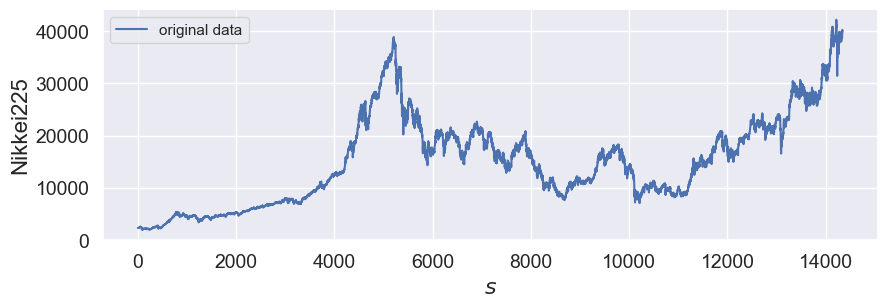} 
    \caption{Original time series $Z(s)$ ($n=14340$)}
    \label{fig:Nikkei225_original}
\end{subfigure}
\begin{subfigure}[b]{0.47\textwidth}
    \centering
    \includegraphics[width=\textwidth, height=1.2in]{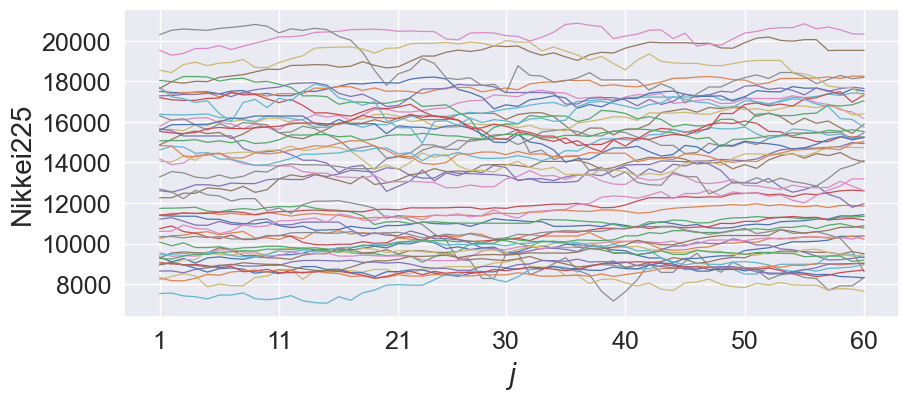} 
    \captionsetup{justification=centering} 
    \caption{Continuous sample curves $X_{t, T}(j)$;\\
    $t=151, \ldots, 200$}
    \label{fig:Nikkei225_covariate}
\end{subfigure}
\begin{subfigure}[b]{0.47\textwidth}
    \centering
    \includegraphics[width=\textwidth, height=1.2in]{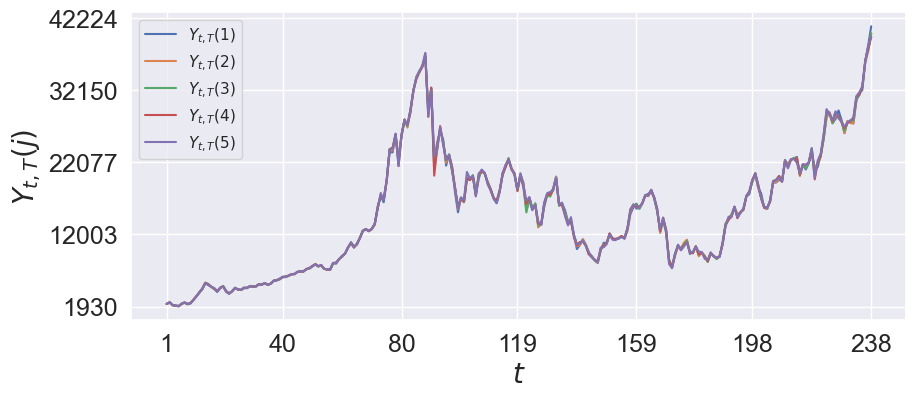} 
    \captionsetup{justification=centering} 
    \caption{Response $Y_{t, 238}(j)$;\\
    for some $j=1, \ldots, 5$}
    \label{fig:Nikkei225_response}
\end{subfigure}
\caption{Nikkie225 time series from Jan. 14, 1970 - Dec. 31, 2024}
\label{fig:Nikkei225}
\end{figure}

\noindent{\scshape Example 4. Nikkei225 data\footnote{Obtained from \url{https://fred.stlouisfed.org/series/NIKKEI225}}.} We next use the Nikkei stock market index dataset or Nikkei225, a key indicator of the Japanese stock market's overall health. The index tracks the performance of 225 large and active companies listed on the Tokyo Stock Exchange (TSE) \citep{bjorgvinsson2024introducing}. We consider 14340 Nikkei225 data points covering January 14, 1971 to December 31, 2024, plotted in Figure \ref{fig:Nikkei225_original}. We construct 239 continuous sample curves by segmenting the original time series $\{Z(s)\}_{s=1,\ldots,14340}$ by 60 observations. Figure \ref{fig:Nikkei225_covariate} reflects 50 examples of the generated continuous sample curves.  We use the same method in Example 5 to generate the functional pairs $\{(X_{t, 238}, Y_{t, 238}(j))\}_{t = 1,\ldots, 238}$ where $j\in\{1, \ldots,60\}$. The behavior of the response variable is plotted in Figure \ref{fig:Nikkei225_response}.\\

We create copies of these datasets using the same method, Algorithm \ref{alg: real data} used in \cite{Tinioetal1.2024}, that relies on Gaussian smoothed procedure \cite{Nietertetal2021}. For a chosen $j$th continuous sample curve, we add $Z_{t, T} \sim \mathcal{N}(0,\sigma^2)$ to each data observation $Y_{t,T}$ with $\sigma > 0$, for all $t \in \{1, \ldots, T\}$. The Gaussian-smoothed datasets are replicated $L=500$ times. We then calculate NW conditional CDF for each replicate at a specific time point $t$. We measure the Wasserstein distance between the average NW and the empirical conditional CDFs.

\begin{figure}[ht]
    \centering
    \begin{subfigure}[t]{\linewidth}
    \begin{subfigure}{0.32\textwidth}
        \centering
        \subcaption*{$T = 24$}
        \includegraphics[width=\linewidth]{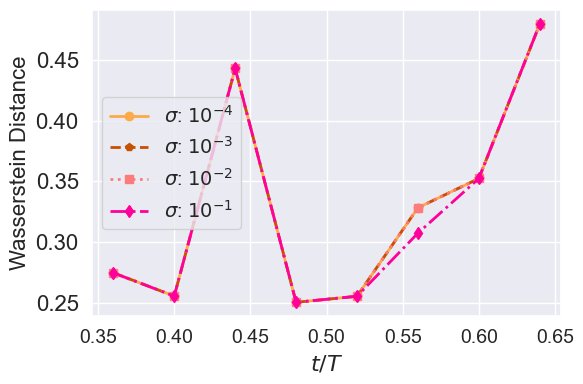}
    \end{subfigure}
    \begin{subfigure}{0.32\textwidth}
        \centering
        \subcaption*{$T = 74$}
        \includegraphics[width=\linewidth]{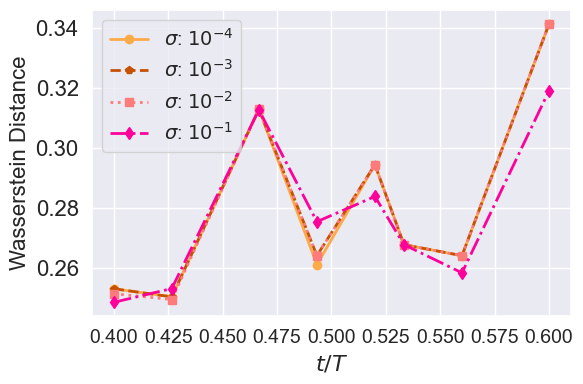}
    \end{subfigure}
    \begin{subfigure}{0.32\textwidth}
        \centering
        \subcaption*{$T = 149$}
        \includegraphics[width=\linewidth]{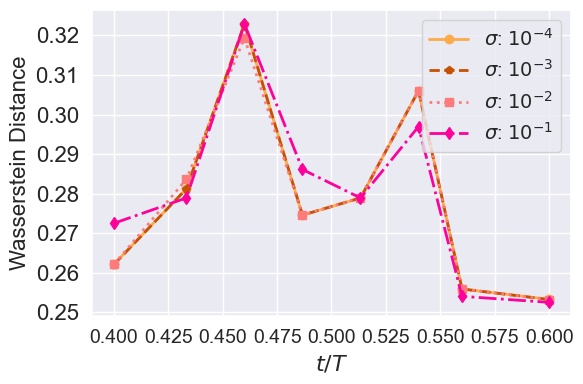}
    \end{subfigure}
    \subcaption{SST}
    \end{subfigure}
                \begin{subfigure}[t]{\linewidth}
    \begin{subfigure}{0.32\textwidth}
        \centering
        \subcaption*{$T = 238$}
        \includegraphics[width=\linewidth]{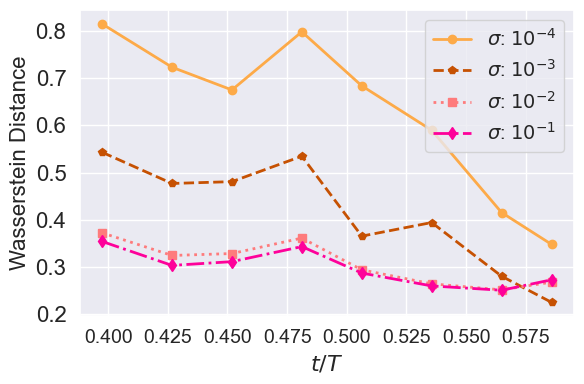}
    \end{subfigure}
    \begin{subfigure}{0.32\textwidth}
        \centering
        \subcaption*{$T = 477$}
        \includegraphics[width=\linewidth]{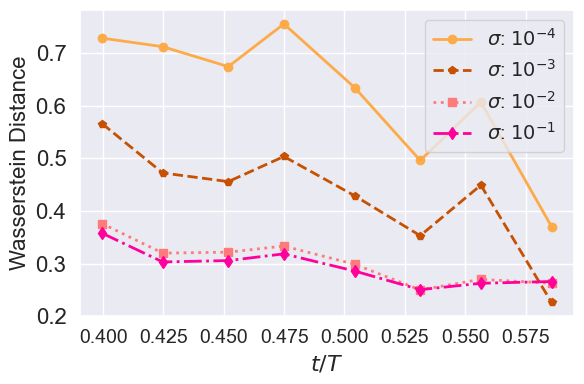}
    \end{subfigure}
    \begin{subfigure}{0.32\textwidth}
        \centering
        \subcaption*{$T = 955$}
        \includegraphics[width=\linewidth]{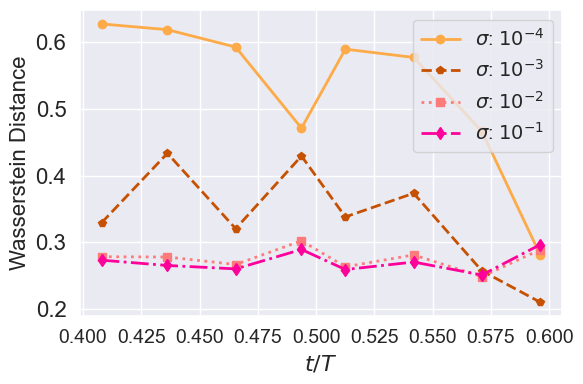}
    \end{subfigure}
    \subcaption{Nikkie225}
    \end{subfigure}
            \caption{Wasserstein distance at various $u = \frac{t}{T}$ with different smoothness level $\sigma$, $K_1= \texttt{Uniform} \text{ and } K_2 = \texttt{Silverman}$ at increasing $T$ using $L=500$ replications. }
    \label{fig: real datasets W1 various sigmas}
\end{figure}

We refine the segmentation of each dataset to increase the sample size, $T$.  By dividing the 900 monthly SST observations into segments of 12 and 6 months, we generate sample sizes of $T=74$ and $T=149$, respectively.  Similarly, we split the 14340 observations of Nikkei225 into 30 and 15 intervals, yielding sample sizes of $T=477$ and $T=955$, respectively. Hence, in this experiment, we set $T=24, 74, 149$ for SST and $T=238, 477, 955$ for Nikkei225.  We set different values of the smoothing parameter $\sigma \in  \{10^{-1}, 10^{-2}, 10^{-3}, 10^{-4}\}$. We use Uniform and Silverman kernels for $K_1$ and $K_2$, respectively, to quantify NW conditional CDF. Like synthetic data experiments, the bandwidths are selected using a cross-validation method. We select $t$ such that $\frac{t}{T} \in [h, 1-h]$ since we use a uniform kernel for $K_1$. For SST, we selected $j = 21$ from $\{1, \ldots, 36\}$, $j= 5$ from $\{1, \ldots, 12\}$, and $j=3$ from $\{1, \ldots, 6\}$. Then, for Nikkei225, we fixed $j = 10$ from $\{1, \ldots, 60\}$, $j= 10$ from $\{1, \ldots, 30\}$, and $j=5$ from $\{1, \ldots, 15\}$. The resulting Wasserstein distances are shown in Figure \ref{fig: real datasets W1 various sigmas} that depicts similar behavior with the results in \cite{Tinioetal1.2024}. For each dataset, Wasserstein distances for larger sample sizes are slightly lower and are higher for $\sigma\rightarrow 0$.

\section{Conclusion} \label{sec:conclusion}
We proposed a NW conditional distribution estimator for LSFTS and established its convergence rates with respect to Wasserstein distance. These rates depend on the bandwidth \( h \) and the small ball probability \( \phi(h) \). We provided the convergence rates for a fractal-type process with \( h = \bigO (T^{-\xi}) \) and \( \phi(h) = h^{\tau_0} \), for \( 0 < \xi < \frac{1}{2(1+\tau_0)} \) and \( \tau_0 > 1 \). Numerical synthetic and real-world data experiments were conducted, supported by a data-generating algorithm designed to calculate the NW estimator.

This work also outlines promising directions for future research. One avenue involves modifying the basic indicator function to an integrated kernel \( H_g(y - Y_{t,T}) \), where \( H \) is a smooth cumulative distribution function (CDF) and \( H_g(y - Y_{t,T}) \) serves as a local weighting function with bandwidth \( g \), analogous to \( h \). Another possible extension is amending the NW estimator to handle missing data. While expanding our results to encompass functional ergodic data would be highly valuable, it requires substantial mathematical advancements and lies beyond the current scope of this paper.

\paragraph{Acknowledgements.} Mr. Tinio acknowledges the support provided by the Department of Science and Technology - Science Education Institute (DOST-SEI) in partnership with Campus France through the PhilFrance-DOST Scholarship grant.

\appendix 

\clearpage

\section{Numerical experiment algorithms}
We use the following algorithms, which are based on the approach presented in \cite{Tinioetal1.2024}, to generate data and calculate NW.\\

\LinesNotNumbered
\begin{algorithm}[H]
\small
\DontPrintSemicolon
\SetNlSty{textbf}{}{.}
\SetKwInOut{Input}{input}
\SetKwInOut{Return}{return}
\caption{\mbox{Data generating and NW estimation for synthetic data \cite{Tinioetal1.2024}}} 
\label{alg: simulated data}
\nl \Input{ sample size $T$, time point $t \in \{1, \ldots, T\} $, $N$ spatial discretization points of $\tau \in [0,1]$, $J$ basis functions, number of replications $L$, based kernels $K_1(\cdot), K_2(\cdot)$, bandwidth $h;$}
\nl \For{$l = 1, \ldots, L$}{
    \# \texttt{Generate $l$-th replication process} $\{Y_{a, T}^{(l)}\}_{a=1, \ldots, T}$ \texttt{with functional covariates} $\{X_{a, T}^{(l)}\}_{a=1, \ldots, T}$ \texttt{constructed using} (\ref{eqn: tvFAR1 process}) or (\ref{eqn: tvFAR2 process})\\
    \For{$a=1,\ldots,T$}{
   $Y_{a,T}^{(l)}  \gets  m^\star\big(\frac aT, X^{(l)}_{a,T}\big) + \varepsilon_{a,T}^{(l)};$\\
    }
     \# \texttt{Calculate $l$-th NW conditional CDF estimator}\\
     $\displaystyle \hat{F}_{t}^{(l)}(y|x) \gets \sum_{a=1}^T \omega_{a}(\frac tT,x)  \mathds{1}_{Y_{a,T}^{(l)}\leq y};$
}
\# \texttt{Calculate average NW estimator }\\
\nl $\displaystyle \hat{F}_{t}^L(y|x) \gets  \frac{1}{L} \sum_{l=1}^L \hat{F}_{t}^{(l)}(y|x);$\\
 \# \texttt{Calculate empirical conditional CDF}\\
 \nl $ \displaystyle F_t^L(y|x) \gets \frac{1}{L} \sum_{l=1}^L \mathds{1}_{Y_{t,T}^{(l)}\leq y};$\\
\nl \Return{$W_1(\hat{F}_t^L(y|x),  F_t^L(y|x) );$}
\end{algorithm}

\vspace{6mm}

\LinesNotNumbered
\begin{algorithm}[H]
\small
\DontPrintSemicolon
\SetNlSty{textbf}{}{.}
\SetKwInOut{Input}{input}
\SetKwInOut{Return}{return}
\caption{\mbox{Gaussian smoothed procedure and NW estimation for real datasets \cite{Tinioetal1.2024}}}
\label{alg: real data}
\nl \Input{ real dataset $\{(X_{a, T}, Y_{a, T}(j))\}_{a=1, \ldots, T}$ for fixed $j$,  $\sigma >0$, time point $t \in \{1, \ldots, T\} $, number of replications $L$, based kernels $K_1(\cdot), K_2(\cdot)$, bandwidth $h;$}
\nl \For{$l = 1, \ldots, L$}{
    \# \texttt{Generate $l$-th replication} $\{Y_{a, T}^{(l)}\}_{a=1, \ldots, T}$\\
    \For{$a=1,\ldots,T$}{
    $Y_{a,T}^{(l)} \gets  Y_{a,T}(j) + Z_{a,T}^{(l)} $, where $Z_{a,T}^{(l)} \sim \mathcal{N}(0, \sigma^2);$\\
    }
     \# \texttt{Calculate $l$-th NW conditional CDF estimator}\\
     $\displaystyle \hat{F}_{t}^{(l)}(y|x) \gets \sum_{a=1}^T \omega_{a}(\frac tT,x)  \mathds{1}_{Y_{a,T}^{(l)}\leq y};$
}
                      \# \texttt{Calculate average NW estimator }\\
\nl $\displaystyle \hat{F}_{t}^L(y|x) \gets  \frac{1}{L} \sum_{l=1}^L \hat{F}_{t}^{(l)}(y|x);$\\
 \# \texttt{Calculate empirical conditional CDF}\\
 \nl $ \displaystyle F_t^L(y|x) \gets \frac{1}{L} \sum_{l=1}^L \mathds{1}_{Y_{t,T}^{(l)}\leq y};$\\
\nl \Return{$W_1(\hat{F}_t^L(y|x),  F_t^L(y|x) );$}
\end{algorithm}

\clearpage
\section{Proofs of the main results}\label{Appendix: Proofs of the main results}

We begin with the following propositions that will be useful in the succeeding proofs. For the sake of completeness and consistency, the lines of proofs are adapted from \cite{Tinioetal1.2024} where we introduce the semi-metric $\mathsf{D}(\cdot, \cdot)$.

\begin{proposition}\label{Lemma: E of K2}
    Let Assumptions \ref{Assumption: X is lsp} to \ref{assumption: CDF} hold. 
    Then, for $a,t \in \{1, \ldots, T\}$, the following inequalities hold:
    \begin{enumerate}[label=(\roman*)]
    \item $\E \big[ K_{h,2}(\mathsf{D} (x, X_{a,T})) - K_{h,2}\big(\mathsf{D} \big(x, X_a\big(\frac{a}{T}\big)\big) \big)\big] \leq \frac{L_2 C_U}{Th}.$
    \item $\E \left[ K_{h,2}(\mathsf{D} (x, X_{a,T}))\right] \leq \frac{L_2 C_U}{Th}  + C_d\phi(h)\psi(x).$
    \item $\lefteqn{K_{h,1}\big( \frac{t}{T} - \frac{a}{T} \big) \E \Big[ K_{h,2}(\mathsf{D} (x, X_{t,T})) [\mathds{1}_{Y_{a,T}\leq y} - F_t^\star(\cdot|x)] \Big]}\\
    \hspace{3ex}\qquad \leq (C_1 + C_2) L_{F^\star} K_{h,1}\big( \frac{t}{T} - \frac{a}{T} \big) \Big\{ \frac{L_2 C_U}{T}  + C_dh\phi(h)\psi(x) \Big\},$
    \end{enumerate}
    where $\psi(x)$ is a nonnegative functional in $x\in \mathscr{H}$.
\end{proposition}

\begin{proposition}\label{lemma: J1 is Op(1)}
    Let Assumptions \ref{Assumption: X is lsp} - \ref{Assumption: bandwidth} hold, then
    \begin{align*}
        J_{t,T}^{-1}(\frac{t}{T}, x) = \Big(\frac{1}{Th\phi(h)}\sum_{a=1}^T K_{h,1}\big(\frac{t}{T} - \frac{a}{T}\big) K_{h,2}(\mathsf{D} (x, X_{t,T}))\Big)^{-1} &= \bigO_\P(1).
    \end{align*}
\end{proposition}

\begin{proposition}\label{prop: control of square of sums}
    Let Assumptions \ref{Assumption: X is lsp} - \ref{assumption: blocking} be satisfied. For $y\in \R$ and $x\in \mathscr{H}$, define
    \begin{align*}
        Z_{t,T}(y, x) &= \frac{1}{Th\phi(h)}\sum_{a=1}^T K_{h,1}\big(\frac{t}{T} - \frac{a}{T}\big) K_{h,2}(\mathsf{D} (x, X_{t,T})) \big[\mathds{1}_{Y_{a,T}\leq y} - F^\star_{t}(y|x)\big],
    \end{align*}
    then
    \begin{align*}
    \E\big[Z_{t,T}^2 (y, x)\big]
        =  \bigO\Big(\frac{1}{Th^2\phi^2(h)} + h^2\Big).
    \end{align*}
\end{proposition}

The proofs of Propositions \ref{Lemma: E of K2} to \ref{prop: control of square of sums} are shown in Appendix \ref{Appendix: Proofs of propositions}.

\subsection{Proof of Theorem \ref{Theorem: convergence of EW1}}\label{appendix: proof of convergence of EW1}

Recall that $\pi_t^\star(\cdot|x)$ is the probability measure of the random variable $Y_{t,T}|X_{t,T}=x$ with conditional CDF $F^\star_{t}(y|x) = \mathds{P}[Y_{t,T}\leq y|X_{t,T}=x]$. Observe that, by the definition of $W_1$ given in (\ref{def:W1_cdf}),
\begin{align*}
    \E[W_1(\hat{\pi}_t(\cdot|x),\pi_t^\star(\cdot|x))] 
        &= \E\big[\int \big|\hat{F}_t(y|x)-F^\star_{t}(y|x)\big|\diff y \big] \\
        &= \int \E\big[\big|\hat{F}_t(y|x)-F^\star_{t}(y|x)\big|\big]\diff y,
\end{align*}
using Fubini's theorem. Observe that, using (\ref{def: weights}) and (\ref{eqn:CDF of pi-hat}),
\newpage
\begin{align}\label{eqn: Fhat - Fstar}
    \hat{F}_t(y|x)-F^\star_{t}(y|x)
    &= \frac{\sum_{a=1}^T K_{h,1}\big(\frac{t}{T} - \frac{a}{T}\big) K_{h,2}(\mathsf{D} (x,X_{a,T})) \mathds{1}_{Y_{a,T}\leq y} }{\sum_{a=1}^T K_{h,1}\big(\frac{t}{T} - \frac{a}{T}\big) K_{h,2}(\mathsf{D} (x,X_{a,T}))} - F^\star_{t}(y|x)\nonumber\\
        &= \frac{\frac{1}{Th\phi(h)}\sum_{a=1}^T K_{h,1}\big(\frac{t}{T} - \frac{a}{T}\big) K_{h,2}(\mathsf{D} (x,X_{a,T})) \big[\mathds{1}_{Y_{a,T}\leq y} - F^\star_{t}(y|x)\big]}{\frac{1}{Th\phi(h)}\sum_{a=1}^T K_{h,1}\big(\frac{t}{T} - \frac{a}{T}\big) K_{h,2}(\mathsf{D} (x,X_{a,T}))}.\nonumber\\
\end{align}
Further, by applying Cauchy-Schwarz inequality, we obtain
\begin{align}\label{eqn: EW1 cauchy_schwarz sums}
    \lefteqn{\E[W_1(\hat{\pi}_t(\cdot|x),\pi_t^\star(\cdot|x))]}\nonumber\\ 
        &= \int \E\Big[\Big| \frac{\frac{1}{Th\phi(h)}\sum_{a=1}^T K_{h,1}\big(\frac{t}{T} - \frac{a}{T}\big) K_{h,2}(\mathsf{D} (x,X_{a,T})) \big[\mathds{1}_{Y_{a,T}\leq y} - F^\star_{t}(y|x)\big]}{\frac{1}{Th\phi(h)}\sum_{a=1}^T K_{h,1}\big(\frac{t}{T} - \frac{a}{T}\big) K_{h,2}(\mathsf{D} (x,X_{a,T}))}\Big| \Big] \diff y\nonumber\\
        &\leq \int \Big( \E\Big[\Big(\frac{1}{\frac{1}{Th\phi(h)}\sum_{a=1}^T K_{h,1}\big(\frac{t}{T} - \frac{a}{T}\big) K_{h,2}(\mathsf{D} (x,X_{a,T}))}\Big)^2\Big]\Big)^{\frac{1}{2}}\nonumber\\
        &\quad \times \Big(\E\Big[\Big(\frac{1}{Th\phi(h)}\sum_{a=1}^T K_{h,1}\big(\frac{t}{T} - \frac{a}{T}\big) K_{h,2}(\mathsf{D} (x,X_{a,T})) \big[\mathds{1}_{Y_{a,T}\leq y} - F^\star_{t}(y|x)\big]  \Big)^2\Big]\Big)^{\frac{1}{2}} \diff y.\nonumber\\
\end{align}
Let $J_{t,T}(\frac{t}{T}, x) = \frac{1}{Th\phi(h)}\sum_{a=1}^T K_{h,1}\big(\frac{t}{T} - \frac{a}{T}\big) K_{h,2}(\mathsf{D} (x,X_{a,T}))$. Using Proposition \ref{lemma: J1 is Op(1)}, $J_{t,T}^{-1}(\frac{t}{T}, x)  = \bigO_\P(1)$. Hence, the first term in  (\ref{eqn: EW1 cauchy_schwarz sums}) becomes
\begin{align}\label{eqn: J inv Op1}
    \Big(\E\Big[\Big(\frac{1}{\frac{1}{Th\phi(h)}\sum_{a=1}^T K_{h,1}\big(\frac{t}{T} - \frac{a}{T}\big) K_{h,2}(\mathsf{D} (x,X_{a,T}))}\Big)^2\Big]\Big)^{\frac{1}{2}}
    &= \bigO_{\mathds{P}}(1).
\end{align} 
Additionally, from Proposition \ref{prop: control of square of sums}, the second term in (\ref{eqn: EW1 cauchy_schwarz sums}) is shown to be $\bigO\big( \frac{1}{T^\frac{1}{2} h\phi(h)} + h \big)$. Therefore, from (\ref{eqn: EW1 cauchy_schwarz sums}) and combining (\ref{eqn: bound E of S square}) and (\ref{eqn: J inv Op1}), we have
\begin{align*}
    \E\big[W_1\big(\hat{\pi}_t(\cdot|x), \pi_t^\star(\cdot|x)\big)\big]
    &= \bigO_\P\Big(\frac{1}{T^\frac{1}{2} h\phi(h)} + h\Big).
\end{align*}

\subsection{Proof of Corollary \ref{Remark: bound EW_s-s}}
\label{appendix: proof of EW_s-s}
Using the definition of $W_1$ and noting that $y\in [-M,M]$, we have
\begin{align*}
    W_r^r(\hat{\pi}_t(\cdot|x),\pi_t^\star(\cdot|x))
                            &\leq (2M)^{r-1} \int_{-M}^M |\hat{F}_t(y|x) - F_t^\star (y|x)| \diff y.
\end{align*}
So,
\begin{align*}
    \E[W_r^r(\hat{\pi}_t(\cdot|x),\pi_t^\star(\cdot|x))]
    &\leq (2M)^{r-1} \E\Big[ \int_{-M}^M |\hat{F}_t(y|x) - F_t^\star (y|x)| \diff y \Big]\\
    &\leq (2M)^{r-1} \E[W_1(\hat{\pi}_t(\cdot|x),\pi_t^\star(\cdot|x))].
\end{align*}
By Theorem \ref{Theorem: convergence of EW1}, we get the desired result.

\subsection{Proof of Corollary \ref{corollary: convergence of the 2nd moment}}\label{appendix: proof of convergence of the 2nd moment}

Again, we used the definition of $W_1$ given by (\ref{def:W1_cdf}). Additionally, by using Minkowski's integral inequality given by (\ref{eqn: minkowski's inequality}), for $r=2$ we have
\begin{align*}    \norm{W_1\big(\hat{\pi}_t(\cdot|x), \pi_t^\star(\cdot|x)\big)}_{L_2}
    &= {\Big\|\int_\R \big|\hat{F}_t(y|x)-F^\star_{t}(y|x)\big|\diff y\Big\|}_{L_2} \nonumber\\
        &\leq \int_\R {\big\|\hat{F}_t(y|x)-F^\star_{t}(y|x)\big\|}_{L_2} \diff y \nonumber\\
        &= \int_\R \big( \E\big[\big(\hat{F}_t(y|x)-F^\star_{t}(y|x) \big)^2 \big]\big)^\frac{1}{2} \diff y \nonumber\\
                    &= \int_\R \Big( \E \Big[ \Big( \frac{Z_{t,T}(y, x)}{J_{t,T}(\frac{t}{T}, x)} \Big)^2 \Big] \Big)^\frac{1}{2} \diff y,
\end{align*}
using (\ref{eqn: Fhat - Fstar}) and (\ref{eqn: Z tT}). However, using Proposition \ref{lemma: J1 is Op(1)}, $J_{t,T}^{-1}(\frac{t}{T}, x) = \bigO_\P(1)$. So
\begin{align*}
    \norm{W_1\big(\hat{\pi}_t(\cdot|x), \pi_t^\star(\cdot|x)\big)}_{L_2}
    &\lesssim \int_\R \big( \E\big[ Z_{t,T}^2(y, x) \big] \big)^\frac{1}{2} \diff y \\
        &\lesssim \int_\R \Big( \frac{1}{Th^2\phi^2(h)} + h^2 \Big)^\frac{1}{2} \diff y,
\end{align*}
by Proposition \ref{prop: control of square of sums}. Therefore,
\begin{align*}
    \norm{W_1\big(\hat{\pi}_t(\cdot|x), \pi_t^\star(\cdot|x)\big)}_{L_2}
    &= \bigO_\P\Big(\frac{1}{T^\frac{1}{2} h\phi(h)} + h\Big).
\end{align*}

\subsection{Proof of Proposition \ref{prop: |mhat-m| leq W1}}
\label{appendix: proof of mhat leq W1}

    Observe that
    \begin{align*}
        |\hat{m}(\frac{t}{T},x) - m^\star(\frac{t}{T}, x)|
        &= |\E[\hat{Y}_{t,T}|X_{t,T} = x] - \E[Y_{t,T}|X_{t,T}=x]|\\
                &= \Big|\int_\R \hat{y}\diff \hat{\pi}_t(\cdot|x) - \int_\R y \diff \pi^\star_t(\cdot|x) \Big|\\
                &\leq \sup_{f\in\mathcal{F}} \Big|\int_\R f \diff \hat{\pi}_t(\cdot|x) - \int_\R f \diff \pi^\star_t(\cdot|x) \Big|\\
                &= W_1(\hat{\pi}_t(\cdot|x), \pi^\star_t(\cdot|x)).
    \end{align*}
    The duality formula of the Kantorovich-Rubinstein distance is used in the last equality (see Remark 6.5 in \cite{Villani2009OToldnew}), where $\mathcal{F}$ is the set of all continuous functions satisfying the Lipschitz condition $\norm{f}_{Lip}\leq 1$, i.e., $\sup_{y\neq y'} \frac{|f(y) - f(y')|}{|y - y'|} \leq 1$. Hence,
    \begin{align*}
        \E \big[|\hat{m}(\frac{t}{T},x) - m^\star(\frac{t}{T}, x)|\big] 
        \leq \E \big[W_1(\hat{\pi}_t(\cdot|x), \pi^\star_t(\cdot|x))\big].
    \end{align*}
This finishes the proof.

\subsection{Proof of Proposition \ref{corollary: convergence of EW1 chosen h}}\label{appendix: proof of convergence of EW1 chosen h}

If $h = \bigO(T^{-\xi})$ and $\phi(h) = h^{\tau_0}$, for $\tau_0>1$, then directly from Theorem \ref{Theorem: convergence of EW1},
\begin{align*}
    \E\big[W_1\big(\hat{\pi}_t(\cdot|x), \pi_t^\star(\cdot|x)\big)\big] &\lesssim \frac{1}{T^\frac{1}{2} h\phi(h)} + h\\
        &\lesssim \frac{1}{T^\frac{1}{2} T^{-\xi} T^{-\xi \tau_0}} + \frac{1}{T^{\xi}}\\
        &\lesssim \frac{1}{T^{\frac{1}{2} - \xi(1+\tau_0)}} + \frac{1}{T^{\xi}},
\end{align*}
which goes to zero if $0< \xi < \frac{1}{2(1+\tau_0)}$.

\section{Proofs of Propositions~\ref{Lemma: E of K2},  ~\ref{lemma: J1 is Op(1)}, and ~\ref{prop: control of square of sums}}
\label{Appendix: Proofs of propositions}

\subsection{Proof of Proposition \ref{Lemma: E of K2}}\label{Appendix: Proof of E of K2}

\paragraph{\textit{(i)}} Using Assumption \ref{Assumption: kernel functions}, we note that $K_2$ is Lipshitz. In addition, by Assumption \ref{Assumption: X is lsp} and when $u=\frac{t}{T}$, $\mathsf{D} \big(X_{a,T}, X_a\big(\frac{a}{T}\big)\big) \leq \frac{1}{T}U_{t,T}\big(\frac{a}{T}\big)$, where $\E\big[\big(U_{t,T}\big(\frac{a}{T}\big)\big)^\rho\big]<C_U$. So, 
\begin{align*}
    \lefteqn{\E \big[ K_{h,2}(\mathsf{D} (x, X_{a,T})) - K_{h,2}\big(\mathsf{D} \big(x, X_a\big(\frac{a}{T}\big)\big) \big)\big]}\\
        &= \E \big[ K_{2}\big(\frac{\mathsf{D} (x, X_{a,T})}{h}\big) - K_{2}\big(\frac{\mathsf{D} \big(x, X_a\big(\frac{a}{T}\big)\big)}{h} \big)\big]\\
        &\leq \frac{L_2}{h} \E \big[ \big| \mathsf{D} (x,  X_{a,T}) - \mathsf{D} \big(x,  X_a\big(\frac{a}{T}\big)\big)\big| \big]\\
        &\leq \frac{L_2}{h} \E \big[ \big|  \mathsf{D} \big( X_{a,T}, X_a\big(\frac{a}{T}\big) \big) \big| \big]\\
        &\leq \frac{L_2}{h} \E \big[\big| 
    \frac{1}{T}U_{t,T}\big(\frac{a}{T}\big) \big|\big]\\
    &\leq \frac{L_2 C_U}{Th}.
\end{align*}
\paragraph{\textit{(ii)}} We have
\begin{align*}
    \lefteqn{\E \big[ K_{h,2}(\mathsf{D} (x, X_{a,T}))\big]}\\
        &= \E \big[K_{h,2}(\mathsf{D} (x, X_{a,T})) - K_{h,2}\big(\mathsf{D} \big(x, X_a\big(\frac{a}{T}\big)\big) \big) + K_{h,2}\big(\mathsf{D} \big(x, X_a\big(\frac{a}{T}\big)\big) \big) \big]\\
        &= \E \big[K_{h,2}(\mathsf{D} (x, X_{a,T})) - K_{h,2}\big(\mathsf{D} \big(x, X_a\big(\frac{a}{T}\big)\big) \big)\big] + \E\big[K_{h,2}\big(\mathsf{D} \big(x, X_a\big(\frac{a}{T}\big)\big) \big) \big]\\
        &\leq \frac{L_2 C_U}{Th}  + \E\big[K_{h,2}\big(\mathsf{D} \big(x, X_a\big(\frac{a}{T}\big)\big) \big)\big],
\end{align*}
using \textit{(i)}. Now using Assumption \ref{Assumption: small ball},
\begin{align*}    \E\big[K_{h,2}\big(\mathsf{D} \big(x, X_a\big(\frac{a}{T}\big)\big) \big)\big]
    &\leq \E\Big[ \mathds{1}_{\mathsf{D} \big(x, X_a\big(\frac{a}{T}\big)\big) \leq h} \Big]
        = \P\big[ X_a\big(\frac{a}{T}\big) \in B(x,h)\big]\nonumber\\
        &= F_{t/T}(h;x)
    \leq C_d\phi(h)\psi(x).
\end{align*}
Hence,
\begin{align*}
    \E \big[ K_{h,2}(\mathsf{D} (x, X_{t,T}))\big]
    &\leq \frac{L_2 C_U}{Th}  + C_d\phi(h)\psi(x).
\end{align*}
\paragraph{\textit{(iii)}} Note that using Assumption \ref{assumption: CDF}, $\big| F_a^\star(y|X_{a,T}) - F_t^\star(y|x) \big| \leq L_{F^\star} \big( \mathsf{D} (x, X_{a,T}) + \big| \frac{a}{T} - \frac{t}{T} \big|\big)$. Now see that
\begin{align*}
    \lefteqn{K_{h,1}\big( \frac{t}{T} - \frac{a}{T} \big) \E \Big[ K_{h,2}(\mathsf{D} (x, X_{a,T})) [\mathds{1}_{Y_{a,T}\leq y} - F_t^\star(y|x)] \Big]}\\
        &\leq K_{h,1}\big( \frac{t}{T} - \frac{a}{T} \big) \E \Big[ K_{h,2}(\mathsf{D} (x, X_{a,T})) \E\Big[\big(\mathds{1}_{Y_{a,T}\leq y} - F_t^\star(y|x) \big) \Big| X_{a,T} \Big] \Big]\\
        &\leq K_{h,1}\big( \frac{t}{T} - \frac{a}{T} \big) \E \Big[ K_{h,2}(\mathsf{D} (x, X_{a,T})) \big| F_a^\star(y|X_{a,T}) - F_t^\star(y|x) \big| \Big]\\
        &\leq L_{F^\star} K_{h,1}\big( \frac{t}{T} - \frac{a}{T} \big)  \E \Big[ K_{h,2}(\mathsf{D} (x, X_{a,T})) \big( \mathsf{D} (x, X_{a,T}) + \big| \frac{a}{T} - \frac{t}{T} \big|\big) \Big].
\end{align*}
However, using Assumption \ref{Assumption: kernel functions}, $\mathsf{D} (x, X_{a,T}) \leq C_2 h$ otherwise, $K_{h,2}(\mathsf{D} (x, X_{a,T})) = 0$. Additionally, $\big| \frac{a}{T} - \frac{t}{T} \big| \leq C_1 h$ otherwise, $K_{h,1}\big(\big| \frac{a}{T} - \frac{t}{T} \big|\big) = 0$. So,
\begin{align*}
    \lefteqn{K_{h,1}\big( \frac{t}{T} - \frac{a}{T} \big) \E \Big[ K_{h,2}(\mathsf{D} (x, X_{a,T})) [\mathds{1}_{Y_{a,T}\leq y} - F_t^\star(y|x)] \Big]}\\
        &\leq L_{F^\star} K_{h,1}\big( \frac{t}{T} - \frac{a}{T} \big) \Big\{ C_2 h \E \Big[ K_{h,2}(\mathsf{D} (x, X_{a,T})) \Big] + C_1 h \E \Big[ K_{h,2}(\mathsf{D} (x, X_{a,T})) \Big] \Big\}\\
        &\leq (C_1 + C_2) L_{F^\star} h K_{h,1}\big( \frac{t}{T} - \frac{a}{T} \big) \E \Big[K_{h,2}(\mathsf{D} (x, X_{a,T}))\Big]\\
        &\leq (C_1 + C_2) L_{F^\star} K_{h,1}\big( \frac{t}{T} - \frac{a}{T} \big) \Big\{ \frac{L_2 C_U}{T}  + C_dh\phi(h)\psi(x) \Big\},
\end{align*}
using \textit{(ii)}.

\subsection{Proof of Proposition \ref{lemma: J1 is Op(1)}}\label{Appendix: Proof of J1 is Op(1)}

By applying Theorem 3.1 in \cite{Daisuke2022}, $\Big|J_{t,T} (\frac{t}{T},x) -\E\left[J_{t,T}(\frac{t}{T},x)\right]\Big| = \bigO_\P\Big(\sqrt{\frac{\log T}{Th\phi(h)}}\Big)$. Additionally, using Assumption \ref{Assumption: X is lsp}, $J_{t,T}(\frac{t}{T},x)$ can be decomposed as $J_{t,T}(\frac{t}{T},x) = \widetilde{J}_{t,T}(\frac{t}{T},x) + \Bar{J}_{t,T}(\frac{t}{T},x)$. So,

\begin{align*}
    \Big|J_{t,T}(\frac{t}{T},x)\Big| &= \Big| J_{t,T}(\frac{t}{T},x) - \E[J_{t,T}(\frac{t}{T},x)] + \E[J_{t,T}(\frac{t}{T},x)] \Big| \\
        &\leq \Big| J_{t,T}(\frac{t}{T},x) - \E[J_{t,T}(\frac{t}{T},x)]\Big| + \Big| \E[J_{t,T}(\frac{t}{T},x)] \Big| \\
        &\leq \bigO_\P\Big(\sqrt{\frac{\log T}{Th\phi(h)}}\Big) + \Big|\E[J_{t,T}(\frac{t}{T},x)]\Big|\\
        &\leq \bigO_\P\Big(\sqrt{\frac{\log T}{Th\phi(h)}}\Big) + \Big|\E[\widetilde{J}_{t,T}(\frac{t}{T},x) + \Bar{J}_{t,T}(\frac{t}{T},x)]\Big|\\
        &\leq \bigO_\P\Big(\sqrt{\frac{\log T}{Th\phi(h)}}\Big) + \Big|\E[\widetilde{J}_{t,T}(\frac{t}{T},x)]\Big| + \Big|\E[\Bar{J}_{t,T}(\frac{t}{T},x)]\Big|,
\end{align*}
where
\begin{align*}
    \widetilde{J}_{t,T}(\frac{t}{T},x) 
        &= \frac{1}{Th\phi(h)}\sum_{a=1}^T K_{h,1}\big(\frac{t}{T} - \frac{a}{T}\big) K_{h,2}\big(\mathsf{D} \big(x, X_a\big(\frac{a}{T}\big)\big)\big),
\end{align*}
and
\begin{align*}
    \Bar{J}_{t,T}(\frac{t}{T}, x) &= \frac{1}{Th\phi(h)}\sum_{a=1}^T K_{h,1}\big(\frac{t}{T} - \frac{a}{T}\big) \big\{ K_{h,2}(\mathsf{D} (x, X_{a,T})) - K_{h,2}\big(\mathsf{D} \big(x, X_a\big(\frac{a}{T}\big)\big)\big)\big\}.
\end{align*}
Now, let us first observe $\E[\Bar{J}_{t,T}(\frac{t}{T}, x)]$. Using Assumptions \ref{Assumption: X is lsp} and \ref{Assumption: kernel functions} together with Proposition \ref{Lemma: E of K2}.\textit{i}, we have
\begin{align*}
    \E[\Bar{J}_{t,T}(\frac{t}{T}, x)]
    &= \E\big[\frac{1}{Th\phi(h)}\sum_{a=1}^T K_{h,1}\big(\frac{t}{T} - \frac{a}{T}\big) \big\{ K_{h,2}(\mathsf{D} (x, X_{a,T})) - K_{h,2}\big(\mathsf{D} \big(x, X_a\big(\frac{a}{T}\big)\big)\big)\big\}\big]\\
        &= \frac{1}{Th\phi(h)}\sum_{a=1}^T K_{h,1}\big(\frac{t}{T} - \frac{a}{T}\big) \E\big[\big\{ K_{h,2}(\mathsf{D} (x, X_{a,T})) - K_{h,2}\big(\mathsf{D} \big(x, X_a\big(\frac{a}{T}\big)\big)\big)\big\}\big]\\
        &\leq \frac{1}{Th\phi(h)}\sum_{a=1}^T K_{h,1}\big(\frac{t}{T} - \frac{a}{T}\big)\big(\frac{L_2 C_U}{Th}\big)\\
        &\leq \frac{L_2 C_U}{Th\phi(h)}
    \underbrace{\frac{1}{Th}\sum_{a=1}^T K_{h,1}\big(\frac{t}{T} - \frac{a}{T}\big)}_{\bigO(1)}\\ 
        &\leq \frac{L_2 C_U}{Th\phi(h)}\\
    &\lesssim \frac{1}{Th\phi(h)},
\end{align*}
which converges to zero using Assumption \ref{Assumption: bandwidth}. 
In the lines above, $\frac{1}{Th}\sum_{a=1}^T K_{h,1}\big(\frac{t}{T} - \frac{a}{T}\big) = \bigO(1)$ since, using Lemma B.2 in \cite{vogt2012}, for $I_h = [C_1h, 1 - C_1h]$,
\begin{align}\label{eqn: O1 sum}
    \frac{1}{Th}\sum_{a=1}^T K_{h,1}\big(\frac{t}{T} - \frac{a}{T}\big) 
    &\leq \sup_{u\in I_h} \Big|\frac{1}{Th}\sum_{a=1}^T K_{h,1}\big(u - \frac{a}{T}\big) \Big|\nonumber\\
        &\leq \sup_{u\in I_h} \Big|\frac{1}{Th}\sum_{a=1}^T K_{h,1}\big(u - \frac{a}{T}\big) - 1 \Big| + 1 \nonumber\\
        &= \bigO\Big(\frac{1}{Th^2}\Big) + o(h) + 1
        = \bigO(1).
\end{align}
On the other hand,
\begin{align*}
    \E[\widetilde{J}_{t,T}(\frac{t}{T},x)]
        &= \E\big[\frac{1}{Th\phi(h)}\sum_{a=1}^T K_{h,1}\big(\frac{t}{T} - \frac{a}{T}\big) K_{h,2}\big(\mathsf{D} \big(x, X_a\big(\frac{a}{T}\big)\big)\big)\big]\\
        &= \frac{1}{Th\phi(h)}\sum_{a=1}^T K_{h,1}\big(\frac{t}{T} - \frac{a}{T}\big)\E\big[K_{h,2}\big(\mathsf{D} \big(x, X_a\big(\frac{a}{T}\big)\big)\big)\big].
\end{align*}
Using equation (4.3) in \cite{Ferraty&Vieu2006}, we have
\begin{align*}
    \E[\widetilde{J}_{t,T}(\frac{t}{T},x)]
        &= \frac{1}{Th\phi(h)}\sum_{a=1}^T K_{h,1}\big(\frac{t}{T} - \frac{a}{T}\big)\E\big[\mathds{1}_{(\mathsf{D} (x, X_a(\frac{a}{T})))\leq h}\big]\\
        &= \frac{1}{Th\phi(h)}\sum_{a=1}^T K_{h,1}\big(\frac{t}{T} - \frac{a}{T}\big) F_{t/T}(h;x)\\
    &\geq \frac{1}{\phi(h)}\underbrace{\frac{1}{Th}\sum_{a=1}^T K_{h,1}\big(\frac{t}{T} - \frac{a}{T}\big)}_{\bigO(1)} c_d\phi(h) \psi(x) \quad \text{(using Assumption \ref{Assumption: small ball})}\\
        &\sim \psi(x) > 0,
\end{align*}
which implies that $\E[\widetilde{J}_{t,T}(\frac{t}{T},x)] > 0$. Therefore,

\begin{align*}
    \frac{1}{J_{t,T}(\frac{t}{T},x)}
    &= \frac{1}{ o_{\mathds{P}}(1) +o(1)+\E[\widetilde{J}_{t,T}(\frac{t}{T},x)]}\\
    &= \bigO_{\mathds{P}}(1).
\end{align*}

\subsection{Proof of Proposition \ref{prop: control of square of sums}}\label{Appendix: proof of control of square of sums}

Let
\begin{align}\label{eqn: Z tT}
    Z_{t,T}(y, x) := \frac{1}{Th^{d+1}}\sum_{a=1}^T K_{h,1}\big(\frac{t}{T} - \frac{a}{T}\big) Z_{a,t,T}(y, x),
\end{align}
where $$Z_{a,t,T}(y, x) = K_{h,2}(\mathsf{D} (x, X_{a,T})) \big[\mathds{1}_{Y_{a,T}\leq y} - F^\star_{t}(y|x)\big].$$ 
Applying Bernstein's big-block and small-block procedure on $Z_{t,T}(y, x)$, we partition the set $\{1,\ldots,T\}$ into $2v_{T}+1$ independent subsets: $v_T$ big blocks of size $r_T$, $v_T$ small blocks of size $s_T$, and a remainder block of size $T-v_T(r_T+s_T)$, where $v_T = \lfloor\frac{T}{r_T + s_T}\rfloor$. To establish independence between the blocks, we need to place the asymptotically negligible small blocks in between two consecutive big blocks. This procedure was also used in \citep{Masry2005, MR4783436}. So, we decompose $Z_{t,T} (y, x)$ as
\begin{align}\label{eqn: bernstein blocking}
    Z_{t,T} (y, x) &= \Lambda_{t,T}(y, x) + \Pi_{t,T}(y, x) + \Xi_{t,T}(y, x)\nonumber\\
        &:= \sum_{l=0}^{v_T-1} \Lambda_{l, t,T}(y, x) + \sum_{l=0}^{v_T-1} \Pi_{l, t,T}(y, x) + \Xi_{t,T}(y, x),
\end{align}
where
\begin{align*}
    \Lambda_{l, t,T}(y, x) &= \frac{1}{Th^{d+1}} \sum_{a = l(r_T+s_T)+1}^{l(r_T+s_T)+r_T} K_{h,1}\big(\frac{t}{T} - \frac{a}{T}\big) Z_{a,t,T}(y, x),\\
        \Pi_{l, t,T}(y, x) &= \frac{1}{Th^{d+1}} \sum_{a = l(r_T+s_T)+r_T+1}^{(l+1)(r_T+s_T)} K_{h,1}\big(\frac{t}{T} - \frac{a}{T}\big) Z_{a,t,T}(y, x),
\end{align*}
and
\begin{align*}
    \Xi_{t,T}(y, x) &= \frac{1}{Th^{d+1}} \sum_{a = v_T(r_T+s_T)+1}^T K_{h,1}\big(\frac{t}{T} - \frac{a}{T}\big) Z_{a,t,T}(y, x).
\end{align*}
Let us define the size of the big blocks as $r_T = \lfloor\sqrt{Th\phi(h)}/q_T\rfloor$, where $q_T$ satisfies Assumption \ref{assumption: blocking}, i.e., $q_T = o(\sqrt{Th\phi(h)})$. This further implies that there exists a sequence of positive integers $\{q_T\}$, $q_T \rightarrow \infty$, such that $q_Ts_T = o\big(\sqrt{Th\phi(h)}\big)$. Additionally, as $T\rightarrow \infty$,
\begin{align}\label{eqn: blocking asymptotics}
    \frac{s_T}{r_T}\rightarrow 0, \qquad \text{and} \qquad \frac{r_T}{T}\rightarrow 0.
\end{align}
Note that defining $r_T = \lfloor\sqrt{Th\phi(h)}/q_T\rfloor$ immediately implies that $r_T = o\big(\sqrt{Th\phi(h)}\big)$. Additionally, note that $s_T = o(r_T)$ and $v_T = o(q_T\sqrt{Th\phi(h)})$.
Now,  
\begin{align*}
    \E\big[ Z_{t,T}^2 &(y, x)\big] 
    = \E\big[\Lambda_{t,T}^2(y, x)\big] + \E\big[\Pi_{t,T}^2(y, x)\big] + \E\big[\Xi_{t,T}^2(y, x)\big]\\
    & + 2 \Big\{ \E\big[\Lambda_{t,T}(y, x) \Pi_{t,T}(y, x) \big] + \E\big[\Lambda_{t,T}(y, x) \Xi_{t,T}(y, x) \big] + \E\big[\Pi_{t,T}(y, x) \Xi_{t,T}(y, x) \big] \Big\}.
\end{align*}
However, the defined size of big blocks and the relation (\ref{eqn: blocking asymptotics}) ensure that the blocks are asymptotically independent and the sums of small blocks and the remainder block are asymptotically negligible. Consequently, we can neglect the last terms in the previous equation. Hence, we have
\begin{align*}
    \E\big[ Z_{t,T}^2 (y, x)\big] 
    \approx \E\big[\Lambda_{t,T}^2(y, x)\big] + \E\big[\Pi_{t,T}^2(y, x)\big] + \E\big[\Xi_{t,T}^2(y, x)\big].
\end{align*}
For convenience of notation, in the succeeding steps, we let the dependency on $y$ and $x$ be implicit.\\

\noindent \textcolor{magenta}{\underline{\it Step 1. Control of the big blocks.}} First, let us start by dealing with $\E\big[\Lambda_{t,T}^2\big]$. One has

\begin{align*}
    \E\big[\Lambda_{t,T}^2\big]
        &= \sum_{l=0}^{v_T-1}\E\big[\Lambda_{l,t,T}^2\big] + \sum_{\substack{l=0\\ \quad \mathrlap{l\neq l'}}}^{v_T-1}\sum_{l'=0}^{v_T-1} \E[\Lambda_{l,t,T}]\E[\Lambda_{l',t,T}] \\
        &= \frac{1}{(Th\phi(h))^2} \sum_{l=0}^{v_T-1} \E\Big[\Big(\sum_{a=l(r_T+s_T)+1}^{l(r_T+s_T)+r_T} K_{h,1}\big(\frac{t}{T} - \frac{a}{T}\big)Z_{a,t,T}\Big)^2\Big]\\
        &\qquad + \frac{1}{(Th\phi(h))^2} \sum_{\substack{l=0\\ \quad \mathrlap{l\neq l'}}}^{v_T-1}\sum_{l'=0}^{v_T-1} \sum_{a=l(r_T+s_T)+1}^{l(r_T+s_T)+r_T} \sum_{b=l'(r_T+s_T)+1}^{l'(r_T+s_T)+r_T} K_{h,1}\big(\frac{t}{T} - \frac{a}{T}\big)\\
        &\hspace{8cm} \times K_{h,1}\big(\frac{t}{T} - \frac{b}{T}\big) \E\big[Z_{a,t,T}Z_{b,t,T}\big].
\end{align*}\\
Observe that
\begin{align*}
    \E\big[\Lambda_{t,T}^2\big]
        &= \frac{1}{(Th\phi(h))^2}\sum_{l=0}^{v_T-1} \sum_{a=l(r_T+s_T)+1}^{l(r_T+s_T)+r_T} K_{h,1}^2\big(\frac{t}{T} - \frac{a}{T}\big)\E\big[Z_{a,t,T}^2\big]\\
        &\qquad + \frac{1}{(Th\phi(h))^2} \sum_{l=0}^{v_T-1} \sum_{\substack{a=l(r_T+s_T)+1\\ \qquad \quad \mathrlap{|a-b|>0}}}^{l(r_T+s_T)+r_T} \sum_{b=l(r_T+s_T)+1}^{l(r_T+s_T)+r_T} K_{h,1}\big(\frac{t}{T} - \frac{a}{T}\big)\\
        &\hspace{8cm} \times K_{h,1}\big(\frac{t}{T} - \frac{b}{T}\big) \E\big[Z_{a,t,T}Z_{b,t,T}\big]\\
        &\qquad \quad + \frac{1}{(Th\phi(h))^2} \sum_{\substack{l=0\\ \quad \mathrlap{l\neq l'}}}^{v_T-1}\sum_{l'=0}^{v_T-1} \sum_{a=l(r_T+s_T)+1}^{l(r_T+s_T)+r_T} \sum_{b=l'(r_T+s_T)+1}^{l'(r_T+s_T)+r_T} K_{h,1}\big(\frac{t}{T} - \frac{a}{T}\big)\\
        &\hspace{8cm} \times K_{h,1}\big(\frac{t}{T} - \frac{b}{T}\big) \E\big[Z_{a,t,T}Z_{b,t,T}\big]\\
        &=: {\mathsf \Sigma}_1^{\Lambda} + {\mathsf \Sigma}_2^{\Lambda} + {\mathsf \Sigma}_3^{\Lambda}.
\end{align*}\\
\noindent \textcolor{magenta}{\underline{\it Step 1.1. Control of $\mathlarger{\mathsf \Sigma}_1^{\Lambda}$.}} Considering $\mathlarger{\mathsf \Sigma}_1^{\Lambda}$, we have
\begin{align*}
    {\mathsf \Sigma}_1^{\Lambda} &= \frac{1}{(Th\phi(h))^2}\sum_{l=0}^{v_T-1} \sum_{a=l(r_T+s_T)+1}^{l(r_T+s_T)+r_T} K_{h,1}^2\big(\frac{t}{T} - \frac{a}{T}\big)\E\big[Z_{a,t,T}^2\big]\nonumber\\
        &= \frac{1}{(Th\phi(h))^2}\sum_{l=0}^{v_T-1} \sum_{a=l(r_T+s_T)+1}^{l(r_T+s_T)+r_T} K_{h,1}^2\big(\frac{t}{T} - \frac{a}{T}\big) \E\Big[  K_{h,2}^2(\mathsf{D} (x, X_{a,T}))(\mathds{1}_{Y_{a,T}\leq y} - F_t^\star (y|x))^2 \Big].
\end{align*}
By Proposition  \ref{Lemma: E of K2}.\textit{iii}, we have
\begin{align*}    \lefteqn{K_{h,1}\big(\frac{t}{T} - \frac{a}{T}\big)\E\Big[K_{h,2}^2(\mathsf{D} (x, X_{a,T}))(\mathds{1}_{Y_{a,T}\leq y} - F_t^\star (y|x))^2 \Big]}\nonumber\\
        &\leq 2 C_2 K_{h,1}\big(\frac{t}{T} - \frac{a}{T}\big)\E\Big[K_{h,2}(\mathsf{D} (x, X_{a,T})) \big| \mathds{1}_{Y_{a,T}\leq y} - F_t^\star (y|x)\big| \Big]\\
        &\leq 2 C_2 (C_1 + C_2) L_{F^\star} K_{h,1}\big( \frac{t}{T} - \frac{a}{T} \big) \Big\{ \frac{L_2 C_U}{T}  + C_dh\phi(h)\psi(x) \Big\}\nonumber\\
        &\lesssim K_{h,1}\big( \frac{t}{T} - \frac{a}{T} \big) \Big( \frac{1}{T}  + h\phi(h) \Big).
\end{align*}
So
\begin{align*}
    {\mathsf \Sigma}_1^{\Lambda} 
    &\lesssim \frac{1}{T^2 h^2 \phi^2(h)} \Big(\frac{1}{T}  + h\phi(h)\Big) \sum_{l=0}^{v_T-1} \sum_{a=l(r_T+s_T)+1}^{l(r_T+s_T)+r_T} K_{h,1}^2\big(\frac{t}{T} - \frac{a}{T}\big) \nonumber\\
        &\leq \frac{C_1}{Th\phi^2(h)} \Big(\frac{1}{T}  + h\phi(h)\Big) \underbrace{\frac{1}{Th}\sum_{a=1}^T K_{h,1}\big(\frac{t}{T} - \frac{a}{T}\big)}_{\bigO(1)},
\end{align*}
using (\ref{eqn: O1 sum}). So we have
\begin{align}\label{eqn: big block 1}
    {\mathsf \Sigma}_1^{\Lambda} 
    &\lesssim \frac{1}{Th\phi^2(h)} \Big(\frac{1}{T}  + h\phi(h)\Big)\nonumber\\
        &\lesssim \frac{1}{T^2 h \phi^2(h)} + \frac{1}{T \phi(h)}\nonumber\\
        &\lesssim \frac{1}{Th\phi^2(h)}.
\end{align}
\noindent\textcolor{magenta}{ \underline{\it Step 1.2. Control of $\mathlarger{\mathsf \Sigma}_2^{\Lambda}$.}}
On the other hand,
\begin{align*}
    {\mathsf \Sigma}_2^{\Lambda} &= \frac{1}{(Th\phi(h))^2} \sum_{l=0}^{v_T-1} \sum_{\substack{a=l(r_T+s_T)+1\\ \qquad \quad \mathrlap{|a-b|>0}}}^{l(r_T+s_T)+r_T} \sum_{b=l(r_T+s_T)+1}^{l(r_T+s_T)+r_T} K_{h,1}\big(\frac{t}{T} - \frac{a}{T}\big) K_{h,1}\big(\frac{t}{T} - \frac{b}{T}\big) \E\big[Z_{a,t,T} Z_{b,t,T}\big]\\
        &= \frac{1}{(Th\phi(h))^2} \sum_{l=0}^{v_T-1} \sum_{\substack{a=l(r_T+s_T)+1\\ \qquad \quad \mathrlap{|a-b|>0}}}^{l(r_T+s_T)+r_T} \sum_{b=l(r_T+s_T)+1}^{l(r_T+s_T)+r_T} K_{h,1}\big(\frac{t}{T} - \frac{a}{T}\big) K_{h,1}\big(\frac{t}{T} - \frac{b}{T}\big) \text{$\C$ov}\big(Z_{a,t,T}, Z_{b,t,T}\big)\\
        &\qquad \quad + \frac{1}{(Th\phi(h))^2} \sum_{l=0}^{v_T-1} \sum_{\substack{a=l(r_T+s_T)+1\\ \qquad \quad \mathrlap{|a-b|>0}}}^{l(r_T+s_T)+r_T} \sum_{b=l(r_T+s_T)+1}^{l(r_T+s_T)+r_T} K_{h,1}\big(\frac{t}{T} - \frac{a}{T}\big)\\
        &\hspace{8cm} \times K_{h,1}\big(\frac{t}{T} - \frac{b}{T}\big) \E\big[Z_{a,t,T}\big] \E\big[Z_{b,t,T}\big]\\
        &:= {\mathsf \Sigma}_{21}^{\Lambda} + {\mathsf \Sigma}_{22}^{\Lambda}.
\end{align*}
\\
\noindent \textcolor{magenta}{\underline{\it Step 1.2.1. Control of $\mathlarger{\mathsf \Sigma}_{21}^{\Lambda}$.}} Looking at $\mathlarger{\mathsf \Sigma}_{21}^{\Lambda}$, we have
\begin{align*}
    {\mathsf \Sigma}_{21}^{\Lambda} &= \frac{1}{(Th\phi(h))^2} \sum_{l=0}^{v_T-1} \sum_{\substack{a=l(r_T+s_T)+1\\ \qquad \quad \mathrlap{|a-b|>0}}}^{l(r_T+s_T)+r_T} \sum_{b=l(r_T+s_T)+1}^{l(r_T+s_T)+r_T} K_{h,1}\big(\frac{t}{T} - \frac{a}{T}\big) K_{h,1}\big(\frac{t}{T} - \frac{b}{T}\big)\\
     &\hspace{9cm} \times \text{$\C$ov}\big(Z_{a,t,T}, Z_{b,t,T}\big)\\
        &= \frac{1}{(Th\phi(h))^2} \sum_{l=0}^{v_T-1} \sum_{\substack{n_1=1 \\  \mathrlap{|n_1-n_2|>0}}}^{r_T} \sum_{n_2=1}^{r_T} K_{h,1}\big(\frac{t}{T} - \frac{\lambda + n_1}{T}\big)K_{h,1}\big(\frac{t}{T} - \frac{\lambda + n_2}{T}\big) \\
     &\hspace{9cm} \times \text{$\C$ov}\big(Z_{\lambda + n_1,t,T}, Z_{\lambda + n_2,t,T}\big) \\
        &\leq \frac{1}{(Th\phi(h))^2} \sum_{l=0}^{v_T-1} \sum_{\substack{n_1=1 \\  \mathrlap{|n_1-n_2|>0}}}^{r_T} \sum_{n_2=1}^{r_T} K_{h,1}\big(\frac{t}{T} - \frac{\lambda + n_1}{T}\big)K_{h,1}\big(\frac{t}{T} - \frac{\lambda + n_2}{T}\big) \\
     &\hspace{9cm} \times \big|\text{$\C$ov}\big(Z_{\lambda + n_1,t,T}, Z_{\lambda + n_2,t,T}\big)\big|, 
\end{align*}
where $\lambda = l(r_T+s_T)$. Note that, by Assumption \ref{Assumption: mixing}, $\{X_{t,T},\varepsilon_{t,T}\}$ is regularly mixing. So using Davydov's inequality \cite{Davydov1973} and Lemma 1 in \cite{Tinioetal1.2024}, $\beta(\sigma(X_{\lambda+n_1,t,T}), \sigma(X_{\lambda+n_2,t,T})) \leq \beta(|n_1 - n_2|)$. Then, for $p>2$, we have
\begin{align}\label{eqn: cov within blocks}
    \lefteqn{K_{h,1}\big(\frac{t}{T} - \frac{\lambda+n_1}{T}\big) K_{h,1}\big(\frac{t}{T} - \frac{\lambda+n_2}{T}\big) \Big|\text{$\C$ov}\big(Z_{\lambda+n_1,t,T}, Z_{\lambda+n_2,t,T}\big)\Big|}\nonumber\\
        &\leq 8 K_{h,1}\big(\frac{t}{T} - \frac{\lambda+n_1}{T}\big) K_{h,1}\big(\frac{t}{T} - \frac{\lambda+n_2}{T}\big) \big\| Z_{\lambda+n_1,t,T}\big\|_{L_p} \big\| Z_{\lambda+n_2,t,T}\big\|_{L_p}  \beta(\sigma(X_{\lambda+n_1,t,T}), \sigma(X_{\lambda+n_2,t,T}))^{1-\frac{2}{p}}\nonumber\\
        &\leq 8 K_{h,1}\big(\frac{t}{T} - \frac{\lambda+n_1}{T}\big) K_{h,1}\big(\frac{t}{T} - \frac{\lambda+n_2}{T}\big)\nonumber\\
        &\quad \times \Big(\E\Big[\Big| K_{h,2}(\mathsf{D} (x, X_{\lambda+n_1,T}))(\mathds{1}_{Y_{\lambda+n_1,T}\leq y} - F_t^\star (y|x))\Big|^p\Big]\Big)^\frac{1}{p}\nonumber\\
        &\quad \times \Big(\E\Big[\Big|K_{h,2}(\mathsf{D} (x, X_{\lambda+n_2,T}))(\mathds{1}_{Y_{\lambda+n_2,T}\leq y} - F_t^\star (y|x))\Big|^p\Big]\Big)^\frac{1}{p} \beta(|n_1 - n_2|)^{1-\frac{2}{p}}\nonumber\\
        &\leq K_{h,1}\big(\frac{t}{T} - \frac{\lambda+n_1}{T}\big) K_{h,1}\big(\frac{t}{T} - \frac{\lambda+n_2}{T}\big)\nonumber\\
        &\quad \times \Big(C_2^{p-1} 2^{p-1}\E\Big[ K_{h,2}(\mathsf{D} (x, X_{\lambda+n_1,T}))\big|\mathds{1}_{Y_{\lambda+n_1,T}\leq y} - F_t^\star (y|x)\big|\Big]\Big)^\frac{1}{p}\nonumber\\
        &\quad \times \Big(C_2^{p-1} 2^{p-1} \E\Big[ K_{h,2}(\mathsf{D} (x, X_{\lambda+n_2,T}))\big|\mathds{1}_{Y_{\lambda+n_2,T}\leq y} - F_t^\star (y|x)\big|\Big]\Big)^\frac{1}{p} \beta(|n_1 - n_2|)^{1-\frac{2}{p}}\nonumber\\
        &\lesssim K_{h,1}\big(\frac{t}{T} - \frac{\lambda+n_1}{T}\big) \Big( \frac{1}{T}  + h\phi(h) \Big)^\frac{1}{p}  K_{h,1}\big(\frac{t}{T} - \frac{\lambda+n_2}{T}\big) \Big(\frac{1}{T}  + h\phi(h)\Big)^\frac{1}{p} \beta(|n_1 - n_2|)^{1-\frac{2}{p}}\nonumber\\
        &\lesssim K_{h,1}\big(\frac{t}{T} - \frac{\lambda+n_1}{T}\big) K_{h,1}\big(\frac{t}{T} - \frac{\lambda+n_2}{T}\big) \Big(\frac{1}{T}  + h\phi(h)\Big)^{\frac{2}{p}} \beta(|n_1 - n_2|)^{1-\frac{2}{p}},
\end{align}
using Proposition \ref{Lemma: E of K2}.\textit{iii}. In consequence,
\begin{align*}
    {\mathsf \Sigma}_{21}^{\Lambda} &\lesssim \frac{1}{(Th\phi(h))^2} \Big(\frac{1}{T}  + h\phi(h)\Big)^{\frac{2}{p}}  \sum_{l=0}^{v_T-1} \sum_{\substack{n_1=1 \\  \mathrlap{|n_1-n_2|>0}}}^{r_T} \sum_{n_2=1}^{r_T}  K_{h,1}\big(\frac{t}{T} - \frac{\lambda + n_1}{T}\big) K_{h,1}\big(\frac{t}{T} - \frac{\lambda + n_2}{T}\big)  \beta(|n_1-n_2|)^{1-\frac{2}{p}}\nonumber\\
                &\leq \frac{C_1^2}{T^2 h^2 \phi^2(h)} \Big(\frac{1}{T}  + h\phi(h)\Big)^{\frac{2}{p}}  \sum_{l=0}^{v_T-1} \sum_{\substack{n_1=1 \\  \mathrlap{|n_1-n_2|>0}}}^{r_T} \sum_{n_2=1}^{r_T} \beta(|n_1-n_2|)^{1-\frac{2}{p}}.
\end{align*}
Using Assumption \ref{Assumption: mixing}, $\sum_{k=1}^{\infty} k^{\zeta} \beta(k)^{1-\frac{2}{p}}<\infty$, which can be expressed as $\sum_{k=1}^{r_T} k^{\zeta} \beta(k)^{1-\frac{2}{p}} + \sum_{k=r_T+1}^{\infty} k^{\zeta} \beta(k)^{1-\frac{2}{p}}$. Now, observe that letting $k = |n_1 -n_2|$ yields
\begin{align*}
    \sum_{\substack{n_1=1 \\  \mathrlap{|n_1-n_2|>0}}}^{r_T} \sum_{n_2=1}^{r_T} \beta(|n_1 - n_2|)^{1-\frac{2}{p}}
        &= \sum_{n_1=1}^{r_T} \Big( \sum_{n_2>n_1}^{r_T} \beta(n_2 - n_1)^{1-\frac{2}{p}} + \sum_{n_2<n_1}^{r_T} \beta(n_1 - n_2)^{1-\frac{2}{p}}\Big)\\
        &= \sum_{n_1=1}^{r_T}\sum_{k>0}^{r_T-n_1} \beta(k)^{1 - \frac{2}{p}} + \sum_{n_2=1}^{r_T}\sum_{k>0}^{r_T-n_2} \beta(k)^{1 - \frac{2}{p}}\\
        &= 2 \sum_{n=1}^{r_T}\sum_{k>0}^{r_T-n} \beta(k)^{1 - \frac{2}{p}}
        \leq 2r_T \sum_{k=1}^{r_T} \beta(k)^{1 - \frac{2}{p}}\\
        &\lesssim  r_T \sum_{k=1}^{r_T} k^{\zeta} \beta(k)^{1 - \frac{2}{p}}
        \leq r_T \sum_{k=1}^{\infty} k^{\zeta} \beta(k)^{1 - \frac{2}{p}},
\end{align*}
since $k^\zeta \geq 1$ for $\zeta > 1 - \frac{2}{p}$, where $p>2$. Hence
\begin{align}\label{eqn: big block 2_1}
    {\mathsf \Sigma}_{21}^{\Lambda} 
    &\leq \frac{C_1^2 r_T}{T^2 h^2 \phi^2(h)} \Big(\frac{1}{T}  + h\phi(h)\Big)^{\frac{2}{p}}  \sum_{l=0}^{v_T-1}\sum_{k=1}^{\infty} k^\zeta \beta(k)^{1-\frac{2}{p}} \nonumber\\
        &\lesssim \frac{v_Tr_T}{T^2 h^2 \phi^2(h)} \Big(\frac{1}{T}  + h\phi(h)\Big)^{\frac{2}{p}} \sum_{k=1}^{\infty} k^{\zeta} \beta(k)^{1-\frac{2}{p}} \nonumber\\
        &\lesssim \frac{1}{T h^2 \phi^2(h)} \Big(\frac{1}{T}  + h\phi(h)\Big)^{\frac{2}{p}}, \quad \text{since } v_Tr_T \leq \frac{T}{r_T}r_T =T, \nonumber\\
        &= \Big(\frac{1}{T^{p} h^{2p} \phi^{2p}(h)} \Big(\frac{1}{T}  + h\phi(h)\Big)^2 \Big)^\frac{1}{p}
        \lesssim \Big(\frac{1}{T^{p} h^{2p} \phi^{2p}(h)}   \Big(\frac{1}{T^2} + h^2 \phi^2(h) \Big) \Big)^\frac{1}{p}\nonumber\\
        &\lesssim \Big(\frac{1}{T^{2+p} h^{2p} \phi^{2p}(h) } + \frac{1}{T^p h^{2p-2} \phi^{2p-2}(h)}\Big)^\frac{1}{p}
        \lesssim \Big(\frac{1}{T^{p} h^{2p} \phi^{2p}(h) }\Big)^\frac{1}{p}
        \lesssim \frac{1}{Th^2\phi^2(h)}.
\end{align}
\\
\noindent {\textcolor{magenta}{\underline{\it Step 1.2.2. Control of $\mathlarger{\mathsf \Sigma}_{22}^{\Lambda}$.}}}
Considering $\mathlarger{\mathsf \Sigma}_{22}^{\Lambda}$, see that
\begin{align*}
    {\mathsf \Sigma}_{22}^{\Lambda} &= \frac{1}{(Th\phi(h))^2} \sum_{l=0}^{v_T-1} \sum_{\substack{a=l(r_T+s_T)+1\\ \qquad \quad \mathrlap{|a-b|>0}}}^{l(r_T+s_T)+r_T} \sum_{b=l(r_T+s_T)+1}^{l(r_T+s_T)+r_T} K_{h,1}\big(\frac{t}{T} - \frac{a}{T}\big) K_{h,1}\big(\frac{t}{T} - \frac{b}{T}\big)  \E\big[Z_{a,t,T}\big] \E\big[Z_{b,t,T}\big]\nonumber\\
        &= \frac{1}{(Th\phi(h))^2} \sum_{l=0}^{k_T-1} \sum_{\substack{n_1=1 \\  \mathrlap{|n_1-n_2|>0}}}^{r_T} \sum_{n_2=1}^{r_T} K_{h,1}\big(\frac{t}{T} - \frac{\lambda + n_1}{T}\big)K_{h,1}\big(\frac{t}{T} - \frac{\lambda + n_2}{T}\big) \E\big[Z_{\lambda + n_1,t,T}\big]\E\big[Z_{\lambda + n_2,t,T}\big]\nonumber\\
        &= \frac{1}{(Th\phi(h))^2} \sum_{l=0}^{v_T-1} \sum_{\substack{n_1=1 \\  \mathrlap{|n_1-n_2|>0}}}^{r_T} \sum_{n_2=1}^{r_T} K_{h,1}\big(\frac{t}{T} - \frac{\lambda + n_1}{T}\big)K_{h,1}\big(\frac{t}{T} - \frac{\lambda + n_2}{T}\big)\nonumber\\
        &\hspace{6cm} \times \E\Big[K_{h,2}(\mathsf{D} (x, X_{\lambda+n_1,T}))(\mathds{1}_{Y_{\lambda+n_1,T}\leq y} - F_t^\star (y|x))\Big]\nonumber\\
        &\hspace{6cm} \times \E\Big[K_{h,2}(\mathsf{D} (x, X_{\lambda+n_2,T}))(\mathds{1}_{Y_{\lambda+n_2,T}\leq y} - F_t^\star (y|x))\Big].
\end{align*}
By Proposition \ref{Lemma: E of K2}.\textit{iii}, for $i=1,2$, $K_{h,1}\big(\frac{t}{T} - \frac{\lambda + n_i}{T}\big)\E\big[K_{h,2}(\mathsf{D} (x, X_{\lambda+n_i,T}))(\mathds{1}_{Y_{\lambda + n_i,T}\leq y} - F_t^\star (y|x))\big] \lesssim K_{h,1}\big(\frac{t}{T} - \frac{\lambda + n_i}{T}\big) \big( \frac{1}{T}  + h\phi(h) \big)$, then
\begin{align}\label{eqn: big block 2_2}
    {\mathsf \Sigma}_{22}^{\Lambda} 
        &\lesssim \frac{1}{(Th\phi(h))^2} \Big(\frac{1}{T}  + h\phi(h)\Big)^2  \sum_{l=0}^{v_T-1} \sum_{\substack{n_1=1 \\  \mathrlap{|n_1-n_2|>0}}}^{r_T} \sum_{n_2=1}^{r_T} K_{h,1}\big(\frac{t}{T} - \frac{\lambda + n_1}{T}\big) K_{h,1}\big(\frac{t}{T} - \frac{\lambda + n_2}{T}\big)\nonumber\\
        &\leq \frac{C_1}{T h \phi^2(h)}  \Big(\frac{1}{T}  + h\phi(h)\Big)^2 \underbrace{\frac{1}{Th} \sum_{a=1}^{T} K_{h,1}\big(\frac{t}{T} - \frac{a}{T}\big)}_{\bigO(1)} \nonumber
        \lesssim  \frac{1}{Th\phi^2(h)} \Big(\frac{1}{T}  + h\phi(h)\Big)^2\nonumber\\
        &\lesssim \frac{1}{Th\phi^2(h)} \Big( \frac{1}{T^2}  + h^2 \phi^2(h) \Big)
        \lesssim \frac{1}{T^3 h \phi^2(h)} + \frac{h}{T}
        \lesssim  \frac{1}{T h \phi^2(h)}.
\end{align}
\noindent {\textcolor{magenta}{\underline{\it Step 1.2.2. Control of $\mathlarger{\mathsf \Sigma}_{22}^{\Lambda}$.}}}
Considering $\mathlarger{\mathsf \Sigma}_{22}^{\Lambda}$, see that
\begin{align*}
    {\mathsf \Sigma}_{3}^{\Lambda} &= \frac{1}{(Th\phi(h))^2} \sum_{\substack{l=0\\ \quad \mathrlap{l\neq l'}}}^{v_T-1}\sum_{l'=0}^{v_T-1} \sum_{a=l(r_T+s_T)+1}^{l(r_T+s_T)+r_T} \sum_{b=l'(r_T+s_T)+1}^{l'(r_T+s_T)+r_T} K_{h,1}\big(\frac{t}{T} - \frac{a}{T}\big) K_{h,1}\big(\frac{t}{T} - \frac{b}{T}\big) \E[Z_{a,t,T}Z_{b,t,T}]\\
        &= \frac{1}{(Th\phi(h))^2} \sum_{\substack{l=0\\ \quad \mathrlap{l\neq l'}}}^{v_T-1}\sum_{l'=0}^{v_T-1} \sum_{a=l(r_T+s_T)+1}^{l(r_T+s_T)+r_T} \sum_{b=l'(r_T+s_T)+1}^{l'(r_T+s_T)+r_T} K_{h,1}\big(\frac{t}{T} - \frac{a}{T}\big) K_{h,1}\big(\frac{t}{T} - \frac{b}{T}\big) \text{$\C$ov}\big(Z_{a,t,T},Z_{b,t,T}\big)\\
    &\qquad \quad + \frac{1}{(Th\phi(h))^2} \sum_{\substack{l=0\\ \quad \mathrlap{l\neq l'}}}^{v_T-1}\sum_{l'=0}^{v_T-1} \sum_{a=l(r_T+s_T)+1}^{l(r_T+s_T)+r_T} \sum_{b=l'(r_T+s_T)+1}^{l'(r_T+s_T)+r_T} K_{h,1}\big(\frac{t}{T} - \frac{a}{T}\big) K_{h,1}\big(\frac{t}{T} - \frac{b}{T}\big)\\
        &\hspace{8cm} \times  \E[Z_{a,t,T}] \E[Z_{b,t,T}]\\
        &=: {\mathsf \Sigma}_{31}^{\Lambda} + {\mathsf \Sigma}_{32}^{\Lambda}.
\end{align*}
\noindent {\textcolor{magenta}{\underline{\it Step 1.3.1 Control of $\mathlarger{\mathsf \Sigma}_{31}^{\Lambda}$.}}} Looking at $\mathlarger{\mathsf \Sigma}_{31}^{\Lambda}$, we have
\begin{align*}
    {\mathsf \Sigma}_{31}^{\Lambda} &= \frac{1}{(Th\phi(h))^2} \sum_{\substack{l=0\\ \quad \mathrlap{l\neq l'}}}^{v_T-1}\sum_{l'=0}^{v_T-1} \sum_{a=l(r_T+s_T)+1}^{l(r_T+s_T)+r_T} \sum_{b=l'(r_T+s_T)+1}^{l'(r_T+s_T)+r_T} K_{h,1}\big(\frac{t}{T} - \frac{a}{T}\big) K_{h,1}\big(\frac{t}{T} - \frac{b}{T}\big) \text{$\C$ov}\big(Z_{a,t,T},Z_{b,t,T}\big)\\
        &= \frac{1}{(Th\phi(h))^2} \sum_{\substack{l=0 \\ \quad \mathrlap{l\neq l'}}}^{v_T-1}\sum_{l'=0}^{v_T-1}\sum_{n_1 =1}^{r_T}\sum_{n_2=1}^{r_T} K_{h,1}\big(\frac{t}{T} - \frac{\lambda + n_1}{T}\big)K_{h,1}\big(\frac{t}{T} - \frac{\lambda'+n_2}{T}\big) \text{$\C$ov}\big(Z_{\lambda + n_1,t,T},Z_{\lambda'+n_2,t,T}\big),
\end{align*}
where $\lambda = l(r_T+s_T)$ and $\lambda' = l'(r_T+s_T)$, however, for $l\neq l'$, see that
\begin{align*}
    |\lambda - \lambda' + n_1 - n_2| &\geq | l(r_T+s_T) - l'(r_T+s_T) + n_1 - n_2 |\\
        &\geq |(l-l')(r_T+s_T) + n_1 - n_2|\\
                &> s_T,
\end{align*}
since $n_1,n_2\in \{1, \ldots, r_T\}$. So if we let $m = \lambda+n_1$ and $m'=\lambda'+n_2$, we have
\begin{align*}
    {\mathsf \Sigma}_{31}^{\Lambda} &= \frac{1}{(Th\phi(h))^2} \sum_{\substack{m=1 \\ \qquad \mathrlap{|m-m'|> s_T}}}^{v_T(r_T+s_T)-s_T}\sum_{m'= 1}^{v_T(r_T+s_T)-s_T} K_{h,1}\big(\frac{t}{T} - \frac{m}{T}\big) K_{h,1}\big(\frac{t}{T} - \frac{m'}{T}\big)  \text{$\C$ov}\big(Z_{m,t,T},Z_{m',t,T}\big)\\
        &\leq \frac{1}{(Th\phi(h))^2} \sum_{\substack{m=1 \\  \mathrlap{|m-m'|> s_T}}}^{T}\sum_{m'=1}^{T} K_{h,1}\big(\frac{t}{T} - \frac{m}{T}\big) K_{h,1}\big(\frac{t}{T} - \frac{m'}{T}\big) \big|\text{$\C$ov}\big(Z_{m,t,T},Z_{m',t,T}\big)\big|.
\end{align*}
Now, using (\ref{eqn: cov within blocks}), we have
\begin{align*}
    \lefteqn{ K_{h,1}\big(\frac{t}{T} - \frac{m}{T}\big) K_{h,1}\big(\frac{t}{T} - \frac{m'}{T}\big) \big|\text{$\C$ov}\big(Z_{m,t,T},Z_{m',t,T}\big)\big|}\\
        &\lesssim K_{h,1}\big(\frac{t}{T} - \frac{m}{T}\big) K_{h,1}\big(\frac{t}{T} - \frac{m'}{T}\big) \big(\frac{1}{T} + h\phi(h) \Big)^{\frac{2}{p}} \beta(|m - m'|)^{1-\frac{2}{p}}.
\end{align*}
Thus
\begin{align*}
    {\mathsf \Sigma}_{31}^{\Lambda} &\lesssim \frac{1}{T^2 h^2 \phi^2(h)} \Big(\frac{1}{T} + h\phi(h)\Big)^{\frac{2}{p}} \sum_{\substack{m=1 \\ \mathrlap{|m-m'|> s_T}}}^{T}\sum_{m'=1}^{T} K_{h,1}\big(\frac{t}{T} - \frac{m}{T}\big) K_{h,1}\big(\frac{t}{T} - \frac{m'}{T}\big) \beta(|m -m'|)^{1-\frac{2}{p}}\nonumber\\
        &\leq\frac{C_1^2}{T^2 h^2 \phi^2(h)} \Big(\frac{1}{T} + h\phi(h) \Big)^{\frac{2}{p}} \sum_{\substack{m=1 \\ \mathrlap{|m-m'|> s_T}}}^{T}\sum_{m'=1}^{T} \beta(|m - m'|)^{1-\frac{2}{p}}.
\end{align*}
By Assumption \ref{Assumption: mixing}, $\sum_{k=1}^{\infty} k^{\zeta} \beta(k)^{1-\frac{2}{p}}<\infty$. Now, observe that letting $k = |m - m'|$ yields
\begin{align*}
    \sum_{\substack{m=1 \\  \mathrlap{|m-m'|> s_T}}}^{T} \sum_{m'=1}^{T} \beta(|m - m'|)^{1-\frac{2}{p}}
        &\leq C \sum_{k = s_T+1}^{T} \beta(k)^{1-\frac{2}{p}}
        \lesssim \frac{1}{k^\zeta} \sum_{k = s_T+1}^{T} k^\zeta \beta(k)^{1-\frac{2}{p}}\\
        &\leq  \frac{1}{s_T^\zeta} \sum_{k = s_T+1}^{T} k^\zeta \beta(k)^{1-\frac{2}{p}}, \quad \text{since } k > s_T,\\
        &\leq \frac{1}{s_T^\zeta} \sum_{k = s_T+1}^{\infty} k^\zeta \beta(k)^{1-\frac{2}{p}},
\end{align*}
since $\beta(k)\geq 0$ and $\big(\frac{k}{s_T}\big)^\zeta \geq 1$ for $\zeta > 1 - \frac{2}{p}$, where $p>2$. So
\begin{align}\label{eqn: big block 3_1}
    {\mathsf \Sigma}_{31}^{\Lambda} 
        &\leq \frac{C_1^2}{s_T^\zeta T^2h^2 \phi^2(h)} \Big(\frac{1}{T} + h\phi(h)\Big)^{\frac{2}{p}} \sum_{k=s_T+1}^{\infty} k^{\zeta} \beta(k)^{1-\frac{2}{p}}\nonumber\\
        &\lesssim \frac{1}{T^2 h^2 \phi^2(h)} \Big(\frac{1}{T} + h\phi(h)\Big)^{\frac{2}{p}}, \quad \text{since } \frac{1}{s_T^\zeta} \leq 1, \nonumber\\
        &\lesssim \Big(\frac{1}{T^{2p} h^{2p} \phi^{2p}(h)}   \Big(\frac{1}{T} + h\phi(h)\Big)^2 \Big)^\frac{1}{p}
        \lesssim \Big(\frac{1}{T^{2p} h^{2p} \phi^{2p}(h)}   \Big(\frac{1}{T^2}   + h^2 \phi^2(h) \Big) \Big)^\frac{1}{p}\nonumber\\
        &\lesssim \Big(\frac{1}{T^{2p+2} h^{2p} \phi^{2p}(h)} + \frac{1}{T^{2p} h^{2p-2} \phi^{2p-2}(h)}\Big)^\frac{1}{p}\nonumber\\
        &\lesssim \Big(\frac{1}{T^{2p} h^{2p} \phi^{2p}(h)}\Big)^\frac{1}{p} 
        \lesssim \frac{1}{T^2 h^2 \phi^2(h)}.
\end{align}

\noindent {\textcolor{magenta}{\underline{\it Step 1.3.2 Control of $\mathlarger{\mathsf \Sigma}_{32}^{\Lambda}$.}}} In view of $\mathlarger{\mathsf \Sigma}_{32}^{\Lambda}$, observe that
\begin{align*}
    {\mathsf \Sigma}_{32}^{\Lambda} &= \frac{1}{(Th\phi(h))^2} \sum_{\substack{l=0\\ \quad \mathrlap{l\neq l'}}}^{v_T-1}\sum_{l'=0}^{v_T-1} \sum_{a=l(r_T+s_T)+1}^{l(r_T+s_T)+r_T} \sum_{b=l'(r_T+s_T)+1}^{l'(r_T+s_T)+r_T} K_{h,1}\big(\frac{t}{T} - \frac{a}{T}\big) K_{h,1}\big(\frac{t}{T} - \frac{b}{T}\big)\\
        &\hspace{8cm} \times  \E[Z_{a,t,T}] \E[Z_{b,t,T}]\\
        &= \frac{1}{(Th\phi(h))^2} \sum_{\substack{l=0 \\ \quad \mathrlap{l\neq l'}}}^{v_T-1}\sum_{l'=0}^{v_T-1}\sum_{n_1 =1}^{r_T}\sum_{n_2=1}^{r_T} K_{h,1}\big(\frac{t}{T} - \frac{\lambda + n_1}{T}\big)K_{h,1}\big(\frac{t}{T} - \frac{\lambda'+n_2}{T}\big)\\
        &\hspace{8cm} \times \E[Z_{\lambda + n_1,t,T}] \E[Z_{\lambda'+n_2,t,T}].
\end{align*}
Similarly, for $l\neq l'$, $|\lambda - \lambda' + n_1 - n_2| > s_T$, then
\begin{align*}
    {\mathsf \Sigma}_{32}^{\Lambda} &\leq \frac{1}{(Th\phi(h))^2} \sum_{\substack{m=1\\ \mathrlap{|m-m'|> s_T}}}^{T}\sum_{m'=1}^{T} K_{h,1}\big(\frac{t}{T} - \frac{m}{T}\big)K_{h,1}\big(\frac{t}{T} - \frac{m'}{T}\big)\E[Z_{m,t,T}] \E[Z_{m',t,T}]\nonumber\\
        &= \frac{1}{(Th\phi(h))^2} \sum_{\substack{m=1\\ \mathrlap{|m-m'|> s_T}}}^{T}\sum_{m'=1}^{T} K_{h,1}\big(\frac{t}{T} - \frac{m}{T}\big)K_{h,1}\big(\frac{t}{T} - \frac{m'}{T}\big)  \E\Big[K_{h,2}(\mathsf{D} (x, X_{m,T}))(\mathds{1}_{Y_{m,T}\leq y} - F_t^\star (y|x))\Big]\nonumber\\
        &\quad \times  \E\Big[K_{h,2}(\mathsf{D} (x, X_{m',T}))(\mathds{1}_{Y_{m',T}\leq y} - F_t^\star (y|x))\Big].
\end{align*}
Using Proposition \ref{Lemma: E of K2}.\textit{iii}, $K_{h,1}\big(\frac{t}{T} - \frac{m}{T}\big)\E\big[K_{h,2}(\mathsf{D} (x, X_{m,T}))(\mathds{1}_{Y_{m,T}\leq y} - F_t^\star (y|x))\big] \lesssim K_{h,1}\big(\frac{t}{T} - \frac{m}{T}\big) \big( \frac{1}{T}  + h\phi(h) \big)$, then
\begin{align}\label{eqn: big block 3}
    {\mathsf \Sigma}_{32}^{\Lambda} 
        &\lesssim \frac{1}{(Th\phi(h))^2} \Big(\frac{1}{T}  + h\phi(h)\Big)^2 \sum_{\substack{m=1\\ \mathrlap{|m-m'|> s_T}}}^{T}\sum_{m'=1}^{T} K_{h,1}\big(\frac{t}{T} - \frac{m}{T}\big)K_{h,1}\big(\frac{t}{T} - \frac{m'}{T}\big)  \nonumber\\
        &\leq \frac{1}{\phi^2(h)} \Big(\frac{1}{T}  + h\phi(h)\Big)^2 \underbrace{\frac{1}{Th}\sum_{m=1}^{T} K_{h,1}\big(\frac{t}{T} - \frac{m}{T}\big)}_{\bigO(1)} \underbrace{\frac{1}{Th}\sum_{m'=1}^{T}  K_{h,1}\big(\frac{t}{T} - \frac{m'}{T}\big)}_{\bigO(1)}\nonumber\\
        &\lesssim \frac{1}{\phi^2(h)} \Big(\frac{1}{T} + h\phi(h)\Big)^2 \lesssim \frac{1}{\phi^2(h)} \Big(\frac{1}{T^2} + h^2 \phi^2(h)\Big)\nonumber\\
        &\lesssim \frac{1}{T^2 \phi^2(h)} + h^2, 
\end{align}
which goes to zero as $T\rightarrow \infty$ using Assumption \ref{Assumption: bandwidth}. Hence, comparing (\ref{eqn: big block 1}), (\ref{eqn: big block 2_1}), (\ref{eqn: big block 2_2}), (\ref{eqn: big block 3_1}), and (\ref{eqn: big block 3}), we have
\begin{align}\label{eqn: big blocks overall order}
    \E\big[\Lambda_{t,T}^2\big]
        &\lesssim  \frac{1}{Th^2\phi^2(h)} + h^2.
\end{align}
\noindent {\textcolor{magenta}{\underline{\it Step 2. Control of the small blocks.}}} Next, we deal with the small blocks. See that

\begin{align*}
    \E\big[\Pi_{t,T}^2\big]
    &= \E \Big[\sum_{l=0}^{v_T-1} \Pi_{l,t,T}^2 + \sum_{\substack{l=0 \\ \quad \mathrlap{l\neq l'}}}^{v_T-1}\sum_{l'=0}^{v_T-1} \Pi_{l,t,T}\Pi_{l',t,T} \Big]\\
        &= \E\Big[\frac{1}{(Th\phi(h))^2}\sum_{l=0}^{v_T-1} \Big(\sum_{a=l(r_T+s_T)+r_T+1}^{(l+1)(r_T+s_T)} K_{h,1}\big(\frac{t}{T} - \frac{a}{T}\big) Z_{a,t,T}\Big)^2\Big]\\
        &\qquad + \frac{1}{(Th\phi(h))^2} \sum_{\substack{l=0 \\ \quad \mathrlap{l\neq l'}}}^{v_T-1}\sum_{l'=0}^{v_T-1} \sum_{a=l(r_T+s_T)+r_T+1}^{(l+1)(r_T+s_T)} \sum_{b=l'(r_T+s_T)+r_T+1}^{(l'+1)(r_T+s_T)} K_{h,1}\big(\frac{t}{T} - \frac{a}{T}\big)\\
        &\hspace{8cm} \times K_{h,1}\big(\frac{t}{T} - \frac{b}{T}\big) \E\big[Z_{a,t,T}Z_{b,t,T}\big].
\end{align*}
Observe that
\begin{align*}
    \E\big[\Pi_{t,T}^2\big]
        &= \frac{1}{(Th\phi(h))^2}\sum_{l=0}^{v_T-1} \sum_{a=l(r_T+s_T)+r_T+1}^{(l+1)(r_T+s_T)}  K_{h,1}^2\big(\frac{t}{T} - \frac{a}{T}\big) \E\big[Z_{a,t,T}^2\big]\\
        &\qquad + \frac{1}{(Th\phi(h))^2} \sum_{l=0}^{v_T-1} \sum_{\substack{a=l(r_T+s_T)+r_T+1 \\ \qquad \qquad \quad \mathrlap{a \neq b}}}^{(l+1)(r_T+s_T)} \sum_{b=l(r_T+s_T)+r_T+1}^{(l+1)(r_T+s_T)} K_{h,1}\big(\frac{t}{T} - \frac{a}{T}\big) \\
        &\hspace{8cm} \times K_{h,1}\big(\frac{t}{T} - \frac{b}{T}\big)  \E\big[Z_{a,t,T} Z_{b,t,T}\big]\\
        &\qquad \quad + \frac{1}{(Th\phi(h))^2} \sum_{\substack{l=0 \\ \quad \mathrlap{l\neq l'}}}^{v_T-1}\sum_{l'=0}^{v_T-1} \sum_{a=l(r_T+s_T)+r_T+1}^{(l+1)(r_T+s_T)} \sum_{b=l'(r_T+s_T)+r_T+1}^{(l'+1)(r_T+s_T)}  K_{h,1}\big(\frac{t}{T} - \frac{a}{T}\big) \\
        &\hspace{8cm} \times K_{h,1}\big(\frac{t}{T} - \frac{b}{T}\big) \E\big[Z_{a,t,T} Z_{b,t,T}\big]\\
        &=: {\mathsf \Sigma}_{1}^\Pi + {\mathsf \Sigma}_{2}^\Pi + {\mathsf \Sigma}_{3}^\Pi.
\end{align*}
\noindent {\textcolor{magenta}{\underline{\it Step 2.1. Control of $\mathlarger{\mathsf \Sigma}_{1}^\Pi$}}} First, let us consider $\mathlarger{\mathsf \Sigma}_{1}^\Pi$.
\begin{align*}
    {\mathsf \Sigma}_{1}^\Pi &= \frac{1}{(Th\phi(h))^2}\sum_{l=0}^{v_T-1} \sum_{a=l(r_T+s_T)+r_T+1}^{(l+1)(r_T+s_T)}  K_{h,1}^2\big(\frac{t}{T} - \frac{a}{T}\big) \E\big[ Z_{a,t,T}^2\big]\nonumber\\
        &= \frac{1}{(Th\phi(h))^2}\sum_{l=0}^{v_T-1} \sum_{a=l(r_T+s_T)+r_T+1}^{(l+1)(r_T+s_T)} K_{h,1}^2\big(\frac{t}{T} - \frac{a}{T}\big) \E\Big[ K_{h,2}^2(\mathsf{D} (x, X_{a,T}))(\mathds{1}_{Y_{a,T}\leq y} - F_t^\star (y|x))^2 \Big]\\
        &\leq \frac{2 C_2}{(Th\phi(h))^2}\sum_{l=0}^{v_T-1} \sum_{a=l(r_T+s_T)+r_T+1}^{(l+1)(r_T+s_T)} K_{h,1}^2\big(\frac{t}{T} - \frac{a}{T}\big) \E\Big[ K_{h,2}(\mathsf{D} (x, X_{a,T}))\big| \mathds{1}_{Y_{a,T}\leq y} - F_t^\star (y|x)\big| \Big].
\end{align*}
By Proposition \ref{Lemma: E of K2}.\textit{iii}, we get
\begin{align}\label{eqn: small block 1}
    {\mathsf \Sigma}_{1}^\Pi
        &\lesssim \frac{1}{T^2h \phi^2(h)} \Big(\frac{1}{T} + h\phi(h)\Big)  \sum_{l=0}^{v_T-1} \sum_{a=l(r_T+s_T)+r_T+1}^{(l+1)(r_T+s_T)}  K_{h,1}^2\big(\frac{t}{T} - \frac{a}{T}\big) \nonumber\\
        &\leq \frac{C_1}{T^2 h^2 \phi^2(h)} \Big(\frac{1}{T} + h\phi(h)\Big)  \sum_{l=0}^{v_T-1} \sum_{a=l(r_T+s_T)+r_T+1}^{(l+1)(r_T+s_T)}  K_{h,1}\big(\frac{t}{T} - \frac{a}{T}\big) \nonumber\\
        &\leq \frac{C_1}{Th\phi^2(h)} \Big(\frac{1}{T} + h\phi(h)\Big)  \underbrace{\frac{1}{Th} \sum_{a=1}^{T} K_{h,1}\big(\frac{t}{T} - \frac{a}{T}\big)}_{\bigO(1)} \nonumber\\
        &\lesssim \frac{1}{T h \phi^2(h)} \Big(\frac{1}{T} + h\phi(h)\Big)\nonumber\\
        &\lesssim \frac{1}{T^2 h \phi^2(h)} + \frac{1}{T \phi(h)}\nonumber\\
        &\lesssim \frac{1}{Th\phi^2(h)}.
\end{align}
\noindent {\textcolor{magenta}{\underline{\it Step 2.2. Control of $\mathlarger{\mathsf \Sigma}_{2}^\Pi.$}}}
On the other hand,

\begin{align*}
    {\mathsf \Sigma}_{2}^\Pi &= \frac{1}{(Th\phi(h))^2} \sum_{l=0}^{v_T-1} \sum_{\substack{a=l(r_T+s_T)+r_T+1 \\ \qquad \qquad \quad \mathrlap{a\neq b}}}^{(l+1)(r_T+s_T)} \sum_{b=l(r_T+s_T)+r_T+1}^{(l+1)(r_T+s_T)} K_{h,1}\big(\frac{t}{T} - \frac{a}{T}\big)K_{h,1}\big(\frac{t}{T} - \frac{b}{T}\big) \E[Z_{a,t,T}Z_{b,t,T}]\\
        &= \frac{1}{(Th\phi(h))^2} \sum_{l=0}^{v_T-1} \sum_{\substack{n_1=1 \\ \mathrlap{|n_1-n_2|>0}}}^{s_T} \sum_{n_2=1}^{s_T} K_{h,1}\big(\frac{t}{T} - \frac{\lambda+n_1}{T}\big)K_{h,1}\big(\frac{t}{T} - \frac{\lambda+n_2}{T}\big)\\
        &\quad \times \big\{ \text{$\C$ov}\big(Z_{\lambda+n_1,t,T},Z_{\lambda+n_2,t,T}\big) + \E[Z_{\lambda+n_1,t,T}]\E[Z_{\lambda+n_2,t,T}]  \big\},
\end{align*}
where $\lambda = l(r_T+s_T)+r_T$. So
\begin{align*}
    {\mathsf \Sigma}_{2}^\Pi
    &= \frac{1}{(Th\phi(h))^2} \sum_{l=0}^{v_T-1} \sum_{\substack{n_1=1 \\ \mathrlap{|n_1-n_2|>0}}}^{s_T} \sum_{n_2=1}^{s_T} K_{h,1}\big(\frac{t}{T} - \frac{\lambda+n_1}{T}\big) K_{h,1}\big(\frac{t}{T} - \frac{\lambda+n_2}{T}\big) \text{$\C$ov}\Big(Z_{\lambda+n_1,t,T}, Z_{\lambda+n_2,t,T}\Big)\\
        &\qquad + \frac{1}{(Th\phi(h))^2} \sum_{l=0}^{v_T-1} \sum_{\substack{n_1=1 \\ \mathrlap{|n_1-n_2|>0}}}^{s_T} \sum_{n_2=1}^{s_T} K_{h,1}\big(\frac{t}{T} - \frac{\lambda+n_1}{T}\big) K_{h,1}\big(\frac{t}{T} - \frac{\lambda+n_2}{T}\big) \\
        &\hspace{8cm} \times \E\Big[  Z_{\lambda+n_1,t,T}\Big]  \E\Big[Z_{\lambda+n_2,t,T}\Big]\\
        &=: {\mathsf \Sigma}_{21}^\Pi + {\mathsf \Sigma}_{22}^\Pi.
\end{align*}

\noindent {\textcolor{magenta}{\underline{\it Step 2.2.1. Control of $\mathlarger{\mathsf \Sigma}_{21}^\Pi$.}}} Taking $\mathlarger{\mathsf \Sigma}_{21}^\Pi$ into consideration, we have
\begin{align*}
    {\mathsf \Sigma}_{21}^\Pi &= \frac{1}{(Th\phi(h))^2} \sum_{l=0}^{v_T-1} \sum_{\substack{n_1=1 \\ \mathrlap{|n_1-n_2|>0}}}^{s_T} \sum_{n_2=1}^{s_T} K_{h,1}\big(\frac{t}{T} - \frac{\lambda+n_1}{T}\big) K_{h,1}\big(\frac{t}{T} - \frac{\lambda+n_2}{T}\big) \text{$\C$ov}\Big(Z_{\lambda+n_1,t,T}, Z_{\lambda+n_2,t,T}\Big).
\end{align*}
Using (\ref{eqn: cov within blocks}),
\begin{align*}
    \lefteqn{K_{h,1}\big(\frac{t}{T} - \frac{\lambda + n_1}{T}\big)K_{h,1}\big(\frac{t}{T} - \frac{\lambda + n_2}{T}\big)\big|\text{$\C$ov}\big(Z_{\lambda+n_1,t,T},Z_{\lambda+n_2,t,T}\big)\big|}\\
    &\lesssim K_{h,1}\big(\frac{t}{T} - \frac{\lambda + n_1}{T}\big)K_{h,1}\big(\frac{t}{T} - \frac{\lambda + n_2}{T}\big)\big(\frac{1}{T} + h\phi(h)\Big)^{\frac{2}{p}} \beta(|n_1 - n_2|)^{1-\frac{2}{p}}.
\end{align*} 
Thus
\begin{align*}
    {\mathsf \Sigma}_{21}^\Pi 
    &\lesssim \frac{1}{T^2 h^2 \phi^2(h)} \Big(\frac{1}{T} + h\phi(h)\Big)^{\frac{2}{p}} \sum_{l=0}^{v_T-1} \sum_{\substack{n_1=1 \\  \mathrlap{|n_1-n_2|>0}}}^{s_T} \sum_{n_2=1}^{s_T}  K_{h,1}\big(\frac{t}{T} - \frac{\lambda+n_1}{T}\big) K_{h,1}\big(\frac{t}{T} - \frac{\lambda+n_2}{T}\big) \beta(|n_1 - n_2|)^{1-\frac{2}{p}}\nonumber\\
        &\leq \frac{C_1^2}{T^2 h^2 \phi^2(h)} \Big(\frac{1}{T} + h\phi(h)\Big)^{\frac{2}{p}} \sum_{l=0}^{v_T-1} \sum_{\substack{n_1=1 \\  \mathrlap{|n_1-n_2|>0}}}^{s_T} \sum_{n_2=1}^{s_T} \beta(|n_1 - n_2|)^{1-\frac{2}{p}}.
\end{align*}
Using Assumption \ref{Assumption: mixing}, $\sum_{k=1}^{\infty} k^{\zeta} \beta(k)^{1-\frac{2}{p}}<\infty$, which can be expressed as $\sum_{k=1}^{s_T} k^{\zeta} \beta(k)^{1-\frac{2}{p}} + \sum_{k=s_T+1}^{\infty} k^{\zeta} \beta(k)^{1-\frac{2}{p}}$. In addition, letting $k = |n_1 -n_2|$ yields
\begin{align*}
    \sum_{\substack{n_1=1 \\  \mathrlap{|n_1-n_2|>0}}}^{s_T} \sum_{n_2=1}^{s_T} \beta(|n_1 - n_2|)^{1-\frac{2}{p}}
        &= \sum_{n_1=1}^{s_T} \Big( \sum_{n_2>n_1}^{s_T} \beta(n_2 - n_1)^{1-\frac{2}{p}} + \sum_{n_2<n_1}^{s_T} \beta(n_1 - n_2)^{1-\frac{2}{p}}\Big)\\
        &= \sum_{n_1=1}^{s_T}\sum_{k>0}^{s_T-n_1} \beta(k)^{1 - \frac{2}{p}} + \sum_{n_2=1}^{s_T}\sum_{k>0}^{s_T-n_2} \beta(k)^{1 - \frac{2}{p}}\\
        &= 2 \sum_{n=1}^{s_T}\sum_{k>0}^{s_T-n} \beta(k)^{1 - \frac{2}{p}}
                \leq 2s_T \sum_{k=1}^{s_T} \beta(k)^{1 - \frac{2}{p}}\\
        &\lesssim  s_T \sum_{k=1}^{s_T} k^{\zeta} \beta(k)^{1 - \frac{2}{p}}
        \leq s_T \sum_{k=1}^{\infty} k^{\zeta} \beta(k)^{1 - \frac{2}{p}},
\end{align*}
            since $\beta(k)\geq 0$ and $k^\zeta \geq 1$ for $\zeta > 1 - \frac{2}{p}$, where $p>2$. So
\begin{align}\label{eqn: small block 2_1}
    {\mathsf \Sigma}_{21}^\Pi 
    &\leq \frac{C_1^2 s_T}{T^2 h^2 \phi^2(h)} \Big(\frac{1}{T} + h\phi(h)\Big)^{\frac{2}{p}} \sum_{l=0}^{v_T-1} \sum_{k=1}^{\infty} k^{\zeta} \beta(k)^{1-\frac{2}{p}}\nonumber\\
        &\lesssim \frac{v_T s_T}{T^2 h^2 \phi^2(h)} \Big(\frac{1}{T} + h\phi(h)\Big)^{\frac{2}{p}} \sum_{k=1}^{\infty} k^{\zeta} \beta(k)^{1-\frac{2}{p}}\nonumber\\
        &\lesssim \frac{1}{T h^2 \phi^2(h)} \Big(\frac{1}{T} + h\phi(h)\Big)^{\frac{2}{p}}, \quad \text{since } 
    v_Ts_T \leq \frac{T}{s_T}s_T = T, \nonumber\\
        &= \Big(\frac{1}{T^{p} h^{2p} \phi^{2p}(h)} \Big(\frac{1}{T}  + h\phi(h)\Big)^2 \Big)^\frac{1}{p}\
        \lesssim \Big(\frac{1}{T^{p} h^{2p} \phi^{2p}(h)}   \Big(\frac{1}{T^2} + h^2 \phi^2(h) \Big) \Big)^\frac{1}{p}\nonumber\\
        &\lesssim \Big(\frac{1}{T^{2+p} h^{2p} \phi^{2p}(h) } + \frac{1}{T^p h^{2p-2} \phi^{2p-2}(h)}\Big)^\frac{1}{p}        \lesssim \Big(\frac{1}{T^{p} h^{2p} \phi^{2p}(h) }\Big)^\frac{1}{p}        \lesssim \frac{1}{Th^2\phi^2(h)}.
\end{align}

\noindent {\textcolor{magenta}{\underline{\it Step 2.2.2. Control of $\mathlarger{\mathsf \Sigma}_{22}^\Pi$.}}} Next, looking at $\mathlarger{\mathsf \Sigma}_{22}^\Pi$, we have
\begin{align*}
    {\mathsf \Sigma}_{22}^\Pi &= \frac{1}{(Th\phi(h))^2} \sum_{l=0}^{v_T-1} \sum_{\substack{a=l(r_T+s_T)+1\\ \qquad \quad \mathrlap{|a-b|>0}}}^{l(r_T+s_T)+r_T} \sum_{b=l(r_T+s_T)+1}^{l(r_T+s_T)+r_T} K_{h,1}\big(\frac{t}{T} - \frac{a}{T}\big) K_{h,1}\big(\frac{t}{T} - \frac{b}{T}\big) \E\big[Z_{a,t,T}\big] \E\big[Z_{b,t,T}\big]\\
        &= \frac{1}{(Th\phi(h))^2} \sum_{l=0}^{v_T-1} \sum_{\substack{n_1=1 \\ \mathrlap{|n_1-n_2|>0}}}^{s_T} \sum_{n_2=1}^{s_T} K_{h,1}\big(\frac{t}{T} - \frac{\lambda+n_1}{T}\big) K_{h,1}\big(\frac{t}{T} - \frac{\lambda+n_2}{T}\big) \nonumber\\
        &\hspace{6cm} \times \E\Big[K_{h,2}(\mathsf{D} (x, X_{\lambda+n_1,T}))(\mathds{1}_{Y_{\lambda + n_1,T}\leq y} - F_t^\star (y|x))\Big] \nonumber\\
        &\hspace{6cm} \times \E\Big[ K_{h,2}(\mathsf{D} (x, X_{\lambda+n_2,T}))(\mathds{1}_{Y_{\lambda + n_2,T}\leq y} - F_t^\star (y|x))\Big]. 
\end{align*}
Using Proposition \ref{Lemma: E of K2}.\textit{iii}, for $i=1,2$, $K_{h,1}\big(\frac{t}{T} - \frac{\lambda+n_i}{T}\big) \E\big[K_{h,2}(\mathsf{D} (x, X_{\lambda+n_i,T}))(\mathds{1}_{Y_{\lambda + n_i,T}\leq y} - F_t^\star (y|x))\big] \lesssim K_{h,1}\big(\frac{t}{T} - \frac{\lambda+n_i}{T}\big) \big(\frac{1}{T} + h\phi(h)\big)$, then 
\begin{align}\label{eqn: small block 2_2}
    {\mathsf \Sigma}_{22}^\Pi 
        &\lesssim \frac{1}{(Th\phi(h))^2}  \Big(\frac{1}{T} + h\phi(h)\Big)^2 \sum_{l=0}^{v_T-1} \sum_{\substack{n_1=1 \\ \mathrlap{|n_1-n_2|>0}}}^{s_T} \sum_{n_2=1}^{s_T} K_{h,1}\big(\frac{t}{T} - \frac{\lambda+n_1}{T}\big) K_{h,1}\big(\frac{t}{T} - \frac{\lambda+n_2}{T}\big)\nonumber\\
        &\leq \frac{C_1}{T h \phi^2(h)}  \Big(\frac{1}{T}  + h\phi(h)\Big)^2 \underbrace{\frac{1}{Th} \sum_{a=1}^{T} K_{h,1}\big(\frac{t}{T} - \frac{a}{T}\big)}_{\bigO(1)} \nonumber\\
        &\lesssim  \frac{1}{Th\phi^2(h)} \Big(\frac{1}{T}  + h\phi(h)\Big)^2\nonumber\\
        &\lesssim \frac{1}{Th\phi^2(h)} \Big( \frac{1}{T^2}  + h^2 \phi^2(h) \Big)\nonumber\\
        &\lesssim \frac{1}{T^3 h \phi^2(h)} + \frac{h}{T}\nonumber\\ 
        &\lesssim  \frac{1}{T h \phi^2(h)}.
\end{align}
\noindent {\textcolor{magenta}{\underline{\it Step 2.3. Control of $\mathlarger{\mathsf \Sigma}_{3}^\Pi$.}}} Now, let us deal with $\mathlarger{\mathsf \Sigma}_{3}^\Pi$.
\begin{align*}
    {\mathsf \Sigma}_{3}^\Pi &= \frac{1}{(Th\phi(h))^2} \sum_{\substack{l=0 \\ \quad \mathrlap{l\neq l'}}}^{v_T-1}\sum_{l'=0}^{v_T-1} \sum_{a=l(r_T+s_T)+r_T+1}^{(l+1)(r_T+s_T)} \sum_{b=l'(r_T+s_T)+r_T+1}^{(l'+1)(r_T+s_T)} K_{h,1}\big(\frac{t}{T} - \frac{a}{T}\big) K_{h,1}\big(\frac{t}{T} - \frac{b}{T}\big)\\
        &\hspace{8cm} \times  \text{$\C$ov}\big(Z_{a,t,T},Z_{b,t,T}\big)\\
        &\qquad + \frac{1}{(Th\phi(h))^2} \sum_{\substack{l=0 \\ \quad \mathrlap{l\neq l'}}}^{v_T-1}\sum_{l'=0}^{v_T-1} \sum_{a=l(r_T+s_T)+r_T+1}^{(l+1)(r_T+s_T)} \sum_{b=l'(r_T+s_T)+r_T+1}^{(l'+1)(r_T+s_T)} K_{h,1}\big(\frac{t}{T} - \frac{a}{T}\big) K_{h,1}\big(\frac{t}{T} - \frac{b}{T}\big)\\
        &\hspace{8cm} \times  \E\big[Z_{a,t,T}\big] \E\big[Z_{b,t,T}\big]\\
        &= {\mathsf \Sigma}_{31}^\Pi + {\mathsf \Sigma}_{32}^\Pi.
\end{align*}
\noindent {\textcolor{magenta}{\underline{\it Step 2.3.1 Control of $\mathlarger{\mathsf \Sigma}_{31}^\Pi$.}}} Looking at $\mathlarger{\mathsf \Sigma}_{31}^\Pi$, see that
\begin{align*}
    {\mathsf \Sigma}_{31}^\Pi &= \frac{1}{(Th\phi(h))^2} \sum_{\substack{l=0 \\ \quad \mathrlap{l\neq l'}}}^{v_T-1}\sum_{l'=0}^{v_T-1}\sum_{n_1 =1}^{s_T}\sum_{n_2=1}^{s_T} K_{h,1}\big(\frac{t}{T} - \frac{\lambda + n_1}{T}\big)K_{h,1}\big(\frac{t}{T} - \frac{\lambda'+n_2}{T}\big)\\
    &\hspace{8cm} \times \text{$\C$ov}\big(Z_{\lambda + n_1,t,T},Z_{\lambda'+n_2,t,T}\big),
\end{align*}
where $\lambda = l(r_T+s_T)+r_T$ and $\lambda' = l'(r_T+s_T)+r_T$, however, for $l\neq l'$, 
\begin{align*}
    |\lambda - \lambda' + n_1 - n_2| &\geq | l(r_T+s_T)+r_T - l'(r_T+s_T)-r_T + n_1 - n_2 |\\
        &\geq |(l-l')(r_T+s_T) + n_1 - n_2|\\
                &> r_T,
\end{align*}
since $n_1,n_2\in \{1, \ldots, s_T\}$. So if we let $q = \lambda+n_1$ and $q'=\lambda'+n_2$, we have
\begin{align*}
    {\mathsf \Sigma}_{31}^\Pi &= \frac{1}{(Th\phi(h))^2} \sum_{\substack{q=r_T+1 \\ \quad \mathrlap{|q-q'|> r_T}}}^{v_T(r_T+s_T)}\sum_{q'=r_T+1}^{v_T(r_T+s_T)} K_{h,1}\big(\frac{t}{T} - \frac{q}{T}\big) K_{h,1}\big(\frac{t}{T} - \frac{q'}{T}\big)  \text{$\C$ov}\big(Z_{q,t,T}, Z_{q',t,T}\big)\\
        &= \frac{1}{(Th\phi(h))^2} \sum_{\substack{m=1 \\ \qquad \mathrlap{|m-m'|> r_T}}}^{v_T(r_T+s_T)-r_T}\sum_{m'=1}^{v_T(r_T+s_T)-r_T} K_{h,1}\big(\frac{t}{T} - \frac{m}{T}\big) K_{h,1}\big(\frac{t}{T} - \frac{m'}{T}\big)\\
    &\quad \times  \text{$\C$ov}\big(Z_{m,t,T},Z_{m',t,T}\big)\\
        &\leq \frac{1}{(Th\phi(h))^2} \sum_{\substack{m=1 \\  \mathrlap{|m-m'|> r_T}}}^{T}\sum_{m'=1}^{T} K_{h,1}\big(\frac{t}{T} - \frac{m}{T}\big) K_{h,1}\big(\frac{t}{T} - \frac{m'}{T}\big) \big|\text{$\C$ov}\big(Z_{m,t,T},Z_{m',t,T}\big)\big|,
\end{align*}
where $m = q - r_T$ and $m' = q' - r_T$. Now, using (\ref{eqn: cov within blocks}), we have
\begin{align*}
    \lefteqn{ K_{h,1}\big(\frac{t}{T} - \frac{m}{T}\big) K_{h,1}\big(\frac{t}{T} - \frac{m'}{T}\big) \big|\text{$\C$ov}\big(Z_{m,t,T},Z_{m',t,T}\big)\big|}\\
        &\lesssim K_{h,1}\big(\frac{t}{T} - \frac{m}{T}\big) K_{h,1}\big(\frac{t}{T} - \frac{m'}{T}\big) \big(\frac{1}{T} + h\phi(h)\Big)^{\frac{2}{p}} \beta(|m - m'|)^{1-\frac{2}{p}}.
\end{align*}
Thus
\begin{align*}
    {\mathsf \Sigma}_{31}^\Pi &\lesssim \frac{1}{(Th\phi(h))^2} \Big(\frac{1}{T} + h\phi(h)\Big)^{\frac{2}{p}} \sum_{\substack{m=1 \\ \mathrlap{|m-m'|> r_T}}}^{T}\sum_{m'=1}^{T} K_{h,1}\big(\frac{t}{T} - \frac{m}{T}\big) K_{h,1}\big(\frac{t}{T} - \frac{m'}{T}\big) \beta(|m-m'|)^{1-\frac{2}{p}}\nonumber\\
        &\leq \frac{C_1^2}{(Th\phi(h))^2} \Big(\frac{1}{T} + h\phi(h)\Big)^{\frac{2}{p}} \sum_{\substack{m=1 \\ \mathrlap{|m-m'|> r_T}}}^{T}\sum_{m'=1}^{T} \beta(|m-m'|)^{1-\frac{2}{p}}.
\end{align*}
By Assumption \ref{Assumption: mixing}, $\sum_{k=1}^{\infty} k^{\zeta} \beta(k)^{1-\frac{2}{p}}<\infty$, which can be expressed as $\sum_{k=1}^{r_T} k^{\zeta} \beta(k)^{1-\frac{2}{p}} + \sum_{k=r_T+1}^{\infty} k^{\zeta} \beta(k)^{1-\frac{2}{p}}$. Additionally, observe that letting $k = |m - m'|$ yields
\begin{align*}
    \sum_{\substack{m=1 \\  \mathrlap{|m-m'|> r_T}}}^{T} \sum_{m'=1}^{T} \beta(|m - m'|)^{1-\frac{2}{p}}
        &\leq C \sum_{k = r_T+1}^{T} \beta(k)^{1-\frac{2}{p}}\\
        &\lesssim \frac{1}{k^\zeta} \sum_{k = r_T+1}^{T} k^\zeta \beta(k)^{1-\frac{2}{p}}\\
        &\leq  \frac{1}{r_T^\zeta} \sum_{k = r_T+1}^{T} k^\zeta \beta(k)^{1-\frac{2}{p}}, \quad \text{since } k > r_T,\\
        &\leq \frac{1}{r_T^\zeta} \sum_{k = r_T+1}^{\infty} k^\zeta \beta(k)^{1-\frac{2}{p}},
\end{align*}
            since $\beta(k)\geq 0$ and $\big(\frac{k}{r_T}\big)^\zeta \geq 1$ for $\zeta > 1 - \frac{2}{p}$, where $p>2$. So
\begin{align}\label{eqn: small block 3_1}
    {\mathsf \Sigma}_{31}^\Pi 
        &\lesssim \frac{1}{r_T^\zeta T^2 h^2 \phi^2(h)} \Big(\frac{1}{T} + h\phi(h)\Big)^{\frac{2}{p}} \sum_{k=r_T+1}^{\infty} k^{\zeta} \beta(k)^{1-\frac{2}{p}}\nonumber\\
        &\lesssim \frac{1}{ T^2 h^2 \phi^2(h)} \Big(\frac{1}{T} + h\phi(h)\Big)^{\frac{2}{p}}, \text{since } \frac{1}{r_T^\zeta} \leq 1,\nonumber\\
        &\lesssim \Big(\frac{1}{T^{2p} h^{2p} \phi^{2p}(h)}   \Big(\frac{1}{T} + h\phi(h)\Big)^2 \Big)^\frac{1}{p}
        \lesssim \Big(\frac{1}{T^{2p} h^{2p} \phi^{2p}(h)}   \Big(\frac{1}{T^2}   + h^2 \phi^2(h) \Big) \Big)^\frac{1}{p}\nonumber\\
        &\lesssim \Big(\frac{1}{T^{2p+2} h^{2p} \phi^{2p}(h)} + \frac{1}{T^{2p} h^{2p-2} \phi^{2p-2}(h)}\Big)^\frac{1}{p}
        \lesssim \Big(\frac{1}{T^{2p} h^{2p} \phi^{2p}(h)}\Big)^\frac{1}{p}\nonumber\\
        &\lesssim \frac{1}{T^2 h^2 \phi^2(h)}.
\end{align}

\noindent {\textcolor{magenta}{\underline{\it Step 2.3.2 Control of $\mathlarger{\mathsf \Sigma}_{32}^\Pi$.}}} In dealing with $\mathlarger{\mathsf \Sigma}_{32}^\Pi$, observe that 
\begin{align*}
    {\mathsf \Sigma}_{32}^\Pi &= \frac{1}{(Th\phi(h))^2} \sum_{\substack{l=0 \\ \quad \mathrlap{l\neq l'}}}^{v_T-1}\sum_{l'=0}^{v_T-1}\sum_{n_1 =1}^{s_T}\sum_{n_2=1}^{s_T} K_{h,1}\big(\frac{t}{T} - \frac{\lambda + n_1}{T}\big)K_{h,1}\big(\frac{t}{T} - \frac{\lambda'+n_2}{T}\big)\\
        &\hspace{8cm} \times \E\big[Z_{\lambda + n_1,t,T}\big] \E\big[Z_{\lambda'+n_2,t,T}\big]\\
        &= \frac{1}{(Th\phi(h))^2}  \sum_{\substack{l=0 \\ \quad \mathrlap{l\neq l'}}}^{v_T-1}\sum_{l'=0}^{v_T-1}\sum_{n_1 =1}^{s_T}\sum_{n_2=1}^{s_T} K_{h,1}\big(\frac{t}{T} - \frac{\lambda + n_1}{T}\big)K_{h,1}\big(\frac{t}{T} - \frac{\lambda'+n_2}{T}\big)\\
        &\hspace{6cm} \times \E\Big[K_{h,2}(\mathsf{D} (x, X_{\lambda+n_1,T}))(\mathds{1}_{Y_{\lambda + n_1,T}\leq y} - F_t^\star (y|x))\Big] \\
        &\hspace{6cm} \times \E\Big[K_{h,2}(\mathsf{D} (x, X_{\lambda+n_2,T}))(\mathds{1}_{Y_{\lambda' + n_2,T}\leq y} - F_t^\star (y|x))\Big].
\end{align*}
Using Proposition \ref{Lemma: E of K2}.\textit{iii}, $K_{h,1}\big(\frac{t}{T} - \frac{\lambda + n_1}{T}\big)\E\big[K_{h,2}(\mathsf{D} (x, X_{\lambda+n_1,T}))(\mathds{1}_{Y_{\lambda + n_1,T}\leq y} - F_t^\star (y|x))\big] \lesssim K_{h,1}\big(\frac{t}{T} - \frac{\lambda + n_1}{T}\big) \big( \frac{1}{T} + h\phi(h) \big)$, then
\begin{align*}
    {\mathsf \Sigma}_{32}^\Pi 
        &\lesssim \frac{1}{(Th\phi(h))^2} \Big(\frac{1}{T} + h\phi(h)\Big)^2\sum_{\substack{l=0 \\ \quad \mathrlap{l\neq l'}}}^{v_T-1}\sum_{l'=0}^{v_T-1}\sum_{n_1 =1}^{s_T}\sum_{n_2=1}^{s_T} K_{h,1}\big(\frac{t}{T} - \frac{\lambda + n_1}{T}\big) K_{h,1}\big(\frac{t}{T} - \frac{\lambda'+n_2}{T}\big).
\end{align*}
Similarly, for $l\neq l'$, $|\lambda - \lambda' + n_1 - n_2| > r_T$, then
\begin{align}\label{eqn: small block 3_2}
    {\mathsf \Sigma}_{32}^\Pi 
    &\lesssim \frac{1}{(Th\phi(h))^2} \Big(\frac{1}{T} + h\phi(h)\Big)^2 \sum_{\substack{m=1\\ \mathrlap{|m-m'|> r_T}}}^{T}\sum_{m'=1}^{T} K_{h,1}\big(\frac{t}{T} - \frac{m}{T}\big) K_{h,1}\big(\frac{t}{T} - \frac{m'}{T}\big)\nonumber\\
        &\leq \frac{1}{\phi^2(h)} \Big(\frac{1}{T} + h\phi(h)\Big)^2 \underbrace{\frac{1}{Th}\sum_{m=1}^{T} K_{h,1}\big(\frac{t}{T} - \frac{m}{T}\big)}_{\bigO(1)} \underbrace{\frac{1}{Th}\sum_{m'=1}^{T}  K_{h,1}\big(\frac{t}{T} - \frac{m'}{T}\big)}_{\bigO(1)}\nonumber\\
        &\lesssim \frac{1}{\phi^2(h)} \Big(\frac{1}{T} + h\phi(h)\Big)^2 
        \lesssim \frac{1}{\phi^2(h)} \Big(\frac{1}{T^2} + h^2 \phi^2(h)\Big)        \lesssim \frac{1}{T^2 \phi^2(h)} + h^2,  
\end{align}
which goes to zero as $T\rightarrow \infty$ using Assumption \ref{Assumption: bandwidth}. Now, comparing (\ref{eqn: small block 1}), (\ref{eqn: small block 2_1}), (\ref{eqn: small block 2_2}), (\ref{eqn: small block 3_1}), and (\ref{eqn: small block 3_2}), we get
\begin{align}\label{eqn: small blocks overall order}
    \E[\Pi_{t,T}^2] 
        &\lesssim  \frac{1}{Th^2\phi^2(h)} + h^2.
\end{align}

\noindent {\textcolor{magenta}{\underline{\it Step 3. Control of the remainder block.}}} 
Now, let us deal with $\E\big[\Xi_{t,T}^2\big]$. See that
\begin{align*}
    \E\big[\Xi_{t,T}^2\big]
        &= \E\Big[\Big(\frac{1}{Th\phi(h)} \sum_{a = v_T(r_T+s_T)+1}^T K_{h,1}\big(\frac{t}{T} - \frac{a}{T}\big) Z_{a,t,T}\Big)^2\Big]\\
        &= \frac{1}{(Th\phi(h))^2} \sum_{a = v_T(r_T+s_T)+1}^T K_{h,1}^2\big(\frac{t}{T} - \frac{a}{T}\big)\E\big[Z_{a,t,T}^2\big]\\
        &\qquad + \frac{1}{(Th\phi(h))^2}\sum_{\substack{a = v_T(r_T+s_T)+1\\ \qquad \qquad \quad \mathrlap{a\neq b}}}^T \sum_{b = v_T(r_T+s_T)+1}^T K_{h,1}\big(\frac{t}{T} - \frac{a}{T}\big)K_{h,1}\big(\frac{t}{T} - \frac{b}{T}\big) \E \big[Z_{a,t,T} Z_{b,t,T}\big]\\
                &= \frac{1}{(Th\phi(h))^2} \sum_{a = v_T(r_T+s_T)+1}^T K_{h,1}^2\big(\frac{t}{T} - \frac{a}{T}\big)\E\big[Z_{a,t,T}^2\big]\\
        &\qquad + \frac{1}{(Th\phi(h))^2}\sum_{\substack{a = v_T(r_T+s_T)+1\\ \qquad \qquad \quad \mathrlap{a\neq b}}}^T \sum_{b = v_T(r_T+s_T)+1}^T K_{h,1}\big(\frac{t}{T} - \frac{a}{T}\big)K_{h,1}\big(\frac{t}{T} - \frac{b}{T}\big) \text{$\C$ov} \big(Z_{a,t,T}, Z_{b,t,T}\big)\\
                &\qquad \quad + \frac{1}{(Th\phi(h))^2}\sum_{\substack{a = v_T(r_T+s_T)+1\\ \qquad \qquad \quad \mathrlap{a\neq b}}}^T \sum_{b = v_T(r_T+s_T)+1}^T K_{h,1}\big(\frac{t}{T} - \frac{a}{T}\big)K_{h,1}\big(\frac{t}{T} - \frac{b}{T}\big) \E\big[Z_{a,t,T}\big] \E\big[ Z_{b,t,T}\big]\\
                &=: {\mathsf \Sigma}_{1}^\Xi + {\mathsf \Sigma}_{2}^\Xi + {\mathsf \Sigma}_{3}^\Xi.
\end{align*}
\noindent {\textcolor{magenta}{\underline{\it Step 3.1. Control of $\mathlarger{\mathsf \Sigma}_{1}^\Xi$.}}} Considering $\mathlarger{\mathsf \Sigma}_{1}^\Xi$, we have
\begin{align*}
    {\mathsf \Sigma}_{1}^\Xi &= \frac{1}{(Th\phi(h))^2} \sum_{a = v_T(r_T+s_T)+1}^T K_{h,1}^2\big(\frac{t}{T} - \frac{a}{T}\big)\E\big[Z_{a,t,T}^2\big]\nonumber\\
        &= \frac{1}{(Th\phi(h))^2} \sum_{a = v_T(r_T+s_T)+1}^T K_{h,1}^2\big(\frac{t}{T} - \frac{a}{T}\big) \E\Big[K_{h,2}^2(\mathsf{D} (x, X_{a,T}))(\mathds{1}_{Y_{a,T}\leq y} - F_t^\star (y|x))^2 \Big]\nonumber\\
        &\leq \frac{2 C_2 }{(Th\phi(h))^2} \sum_{a = v_T(r_T+s_T)+1}^T K_{h,1}^2\big(\frac{t}{T} - \frac{a}{T}\big) \E\Big[K_{h,2}(\mathsf{D} (x, X_{a,T}))\big| \mathds{1}_{Y_{a,T}\leq y} - F_t^\star (y|x)\big| \Big].
\end{align*}
Using Proposition \ref{Lemma: E of K2}.\textit{iii}, we have
\begin{align*}    {\mathsf \Sigma}_{1}^\Xi 
        &\lesssim \frac{1}{T^2 h^2 \phi^2(h)} \Big(\frac{1}{T} + h\phi(h)\Big) \sum_{a = v_T(r_T+s_T)+1}^T K_{h,1}^2\big(\frac{t}{T} - \frac{a}{T}\big) \nonumber\\
        &\leq \frac{C_1}{Th\phi^2(h)} \Big(\frac{1}{T} + h\phi(h)\Big) \underbrace{\frac{1}{Th}\sum_{a = 1}^T K_{h,1}\big(\frac{t}{T} - \frac{a}{T}\big)}_{\bigO(1)},
\end{align*}
using \ref{eqn: O1 sum}. So
\begin{align}\label{eqn: remainder block 1}
    {\mathsf \Sigma}_{1}^\Xi 
        \lesssim \frac{1}{T h \phi^2(h)} \Big(\frac{1}{T} + h\phi(h)\Big)
        \lesssim \frac{1}{T^2 h \phi^2(h)} + \frac{1}{T \phi(h)}
        \lesssim \frac{1}{Th\phi^2(h)}.
\end{align}

\noindent {\textcolor{magenta}{\underline{\it Step 3.2. Control of $\mathlarger{\mathsf \Sigma}_{2}^\Xi$.}}} Taking $\mathlarger{\mathsf \Sigma}_{2}^\Xi$ into account, we have 
\begin{align*}
    {\mathsf \Sigma}_{2}^\Xi &= \frac{1}{(Th\phi(h))^2}\sum_{\substack{a = v_T(r_T+s_T)+1\\ \qquad \qquad \quad \mathrlap{a\neq b}}}^T \sum_{b = v_T(r_T+s_T)+1}^T K_{h,1}\big(\frac{t}{T} - \frac{a}{T}\big)K_{h,1}\big(\frac{t}{T} - \frac{b}{T}\big) \text{$\C$ov} \big(Z_{a,t,T}, Z_{b,t,T}\big)\\
        &= \frac{1}{(Th\phi(h))^2}\sum_{\substack{n_1 = 1\\ \qquad \quad \mathrlap{|n_1-n_2|>0}}}^{T-v_T(r_T+s_T)} \sum_{n_2 = 1}^{T-v_T(r_T+s_T)} K_{h,1}\big(\frac{t}{T} - \frac{\lambda+n_1}{T}\big)K_{h,1}\big(\frac{t}{T} - \frac{\lambda+n_2}{T}\big) \text{$\C$ov} \big(Z_{\lambda+n_1,t,T}, Z_{\lambda+n_2,t,T}\big),
\end{align*}
where $\lambda = v_T(r_T+s_T)$. Now, using (\ref{eqn: cov within blocks}), we have
\begin{align*}
    \lefteqn{K_{h,1}\big(\frac{t}{T} - \frac{\lambda+n_1}{T}\big)K_{h,1}\big(\frac{t}{T} - \frac{\lambda+n_2}{T}\big)\big|\text{$\C$ov}\big(Z_{\lambda+n_1,t,T},Z_{\lambda+n_2,t,T}\big)\big|}\\
        &\lesssim K_{h,1}\big(\frac{t}{T} - \frac{\lambda+n_1}{T}\big)K_{h,1}\big(\frac{t}{T} - \frac{\lambda+n_2}{T}\big) \Big(\frac{1}{T} + h\phi(h)\Big)^{\frac{2}{p}} \beta(|n_1 - n_2|)^{1-\frac{2}{p}}.
\end{align*}
Thus
\begin{align*}
    {\mathsf \Sigma}_{2}^\Xi &\lesssim \frac{1}{T^2 h^2 \phi^2(h)}\Big(\frac{1}{T} + h\phi(h)\Big)^{\frac{2}{p}}  \sum_{\substack{n_1 = 1\\ \qquad \mathrlap{|n_1 - n_2|>0}}}^{T-v_T(r_T+s_T)} \sum_{n_2 = 1}^{T-v_T(r_T+s_T)} K_{h,1}\big(\frac{t}{T} - \frac{\lambda+n_1}{T}\big)\nonumber\\
        &\hspace{8cm} \times K_{h,1}\big(\frac{t}{T} - \frac{\lambda+n_2}{T}\big) \beta(|n_1 - n_2|)^{1-\frac{2}{p}}\nonumber\\
        &\leq \frac{C_1^2}{T^2 h^2 \phi^2(h)}\Big(\frac{1}{T} + h\phi(h)\Big)^{\frac{2}{p}} \sum_{\substack{n_1 = 1\\ \qquad \mathrlap{|n_1 - n_2|>0}}}^{T-v_T(r_T+s_T)} \sum_{n_2 = 1}^{T-v_T(r_T+s_T)} \beta(|n_1 - n_2|)^{1-\frac{2}{p}}.
\end{align*}
Assumption \ref{Assumption: mixing} entails $\sum_{k=1}^{\infty} k^{\zeta} \beta(k)^{1-\frac{2}{p}}<\infty$. Moreover, letting $k = |n_1 - n_2|$ and $w_T = T - v_T(r_T + s_T)$ yields
\begin{align*}
    \sum_{\substack{n_1=1 \\  \mathrlap{|n_1-n_2|>0}}}^{w_T} \sum_{n_2=1}^{w_T} \beta(|n_1 - n_2|)^{1-\frac{2}{p}}
        &= \sum_{n_1=1}^{w_T} \Big( \sum_{n_2>n_1}^{w_T} \beta(n_2 - n_1)^{1-\frac{2}{p}} + \sum_{n_2<n_1}^{w_T} \beta(n_1 - n_2)^{1-\frac{2}{p}}\Big)\\
        &= \sum_{n_1=1}^{w_T}\sum_{k>0}^{w_T-n_1} \beta(k)^{1 - \frac{2}{p}} + \sum_{n_2=1}^{w_T}\sum_{k>0}^{w_T-n_2} \beta(k)^{1 - \frac{2}{p}}\\
        &= 2 \sum_{n=1}^{w_T}\sum_{k>0}^{w_T-n} \beta(k)^{1 - \frac{2}{p}}
        \leq 2w_T \sum_{k=1}^{w_T} \beta(k)^{1 - \frac{2}{p}}\\
        &\lesssim  w_T \sum_{k=1}^{w_T} k^{\zeta} \beta(k)^{1 - \frac{2}{p}}
        \leq w_T \sum_{k=1}^{\infty} k^{\zeta} \beta(k)^{1 - \frac{2}{p}},
\end{align*}
            since $\beta(k)\geq 0$ and $k^\zeta \geq 1$ for $\zeta > 1 - \frac{2}{p}$, where $p>2$. So
\begin{align}\label{eqn: remainder block 2}
    {\mathsf \Sigma}_{2}^\Xi 
        &\leq \frac{C_1^2 w_T }{T^2 h^2 \phi^2(h)}\Big(\frac{1}{T} + h\phi(h)\Big)^{\frac{2}{p}} \sum_{k=1}^{\infty} k^\zeta \beta(k)^{1-\frac{2}{p}}\nonumber\\
        &\lesssim \frac{1}{T h^2 \phi^2(h)}\Big(\frac{1}{T} + h\phi(h)\Big)^{\frac{2}{p}}, \quad \text{since } w_T \ll T,\nonumber\\
        &=\Big( \frac{1}{T^{p} h^{2p} \phi^{2p}(h)}   \Big( \frac{1}{T} + h\phi(h)\Big)^2 \Big)^\frac{1}{p}\nonumber\\
        &\lesssim \Big(\frac{1}{T^{p} h^{2p} \phi^{2p}(h)}   \Big(\frac{1}{T^2} + h^2 \phi^2(h) \Big) \Big)^\frac{1}{p}\nonumber\\
        &\lesssim \Big(\frac{1}{T^{2+p} h^{2p} \phi^{2p}(h) } + \frac{1}{T^p h^{2p-2} \phi^{2p-2}(h)}\Big)^\frac{1}{p}\nonumber\\
        &\lesssim \Big(\frac{1}{T^{p} h^{2p} \phi^{2p}(h) }\Big)^\frac{1}{p}\nonumber\\
        &\lesssim \frac{1}{Th^2\phi^2(h)}.
\end{align}

\noindent {\textcolor{magenta}{\underline{\it Step 3.3. Control of $\mathlarger{\mathsf \Sigma}_{3}^\Xi$.}}} Lastly, let us look at $\mathlarger{\mathsf \Sigma}_{3}^\Xi$. 
\begin{align*}
    {\mathsf \Sigma}_{3}^\Xi &= \frac{1}{(Th\phi(h))^2}\sum_{\substack{a = v_T(r_T+s_T)+1\\ \qquad \qquad \quad \mathrlap{a\neq b}}}^T \sum_{b = v_T(r_T+s_T)+1}^T K_{h,1}\big(\frac{t}{T} - \frac{a}{T}\big)K_{h,1}\big(\frac{t}{T} - \frac{b}{T}\big) \E\big[Z_{a,t,T}\big] \E\big[ Z_{b,t,T}\big]\nonumber\\
        &= \frac{1}{(Th\phi(h))^2}\sum_{\substack{n_1 = 1\\ \mathrlap{|n_1 - n_2|>0}}}^{w_T} \sum_{n_2 = 1}^{w_T} K_{h,1}\big(\frac{t}{T} - \frac{\lambda+n_1}{T}\big)K_{h,1}\big(\frac{t}{T} - \frac{\lambda+n_2}{T}\big)  \E\big[Z_{\lambda+n_1,t,T}\big] \E\big[Z_{\lambda+n_2,t,T}\big]\nonumber\\
        &= \frac{1}{(Th\phi(h))^2}\sum_{\substack{n_1 = 1\\ \mathrlap{|n_1 - n_2|>0}}}^{w_T} \sum_{n_2 = 1}^{w_T} K_{h,1}\big(\frac{t}{T} - \frac{\lambda+n_1}{T}\big)K_{h,1}\big(\frac{t}{T} - \frac{\lambda+n_2}{T}\big)\nonumber\\
        &\hspace{6cm} \times  \E\Big[\prod_{j=1}^d K_{h,2}(x^j-X_{\lambda+n_1,T}^j) (\mathds{1}_{Y_{\lambda+n_1,T}\leq y} - F_t^\star (y|x))\Big]\nonumber\\
        &\hspace{6cm} \times \E\Big[ \prod_{j=1}^d K_{h,2}(x^j-X_{\lambda+n_2,T}^j) (\mathds{1}_{Y_{\lambda+n_2,T}\leq y} - F_t^\star (y|x))\Big].
\end{align*}
Using Proposition \ref{Lemma: E of K2}.\textit{iii}, for $i=1,2$, $K_{h,1}\big(\frac{t}{T} - \frac{\lambda + n_i}{T}\big)\E\big[K_{h,2}(\mathsf{D} (x, X_{\lambda+n_i,T}))(\mathds{1}_{Y_{\lambda + n_i,T}\leq y} - F_t^\star (y|x))\big] \lesssim K_{h,1}\big(\frac{t}{T} - \frac{\lambda + n_i}{T}\big) \big( \frac{1}{T} + h\phi(h) \big)$, then

\begin{align*}    {\mathsf \Sigma}_{3}^\Xi 
        &\lesssim \frac{1}{(Th\phi(h))^2} \Big(\frac{1}{T} + h\phi(h)\Big)^2 \sum_{\substack{n_1 = 1\\ \mathrlap{|n_1 - n_2|>0}}}^{w_T} \sum_{n_2 = 1}^{w_T} K_{h,1}\big(\frac{t}{T} - \frac{\lambda+n_1}{T}\big)K_{h,1}\big(\frac{t}{T} - \frac{\lambda+n_2}{T}\big)\nonumber\\
        &\leq \frac{C_1}{T h \phi^2(h)} \Big(\frac{1}{T} + h\phi(h)\Big)^2 \underbrace{\frac{1}{Th}\sum_{a=1}^T K_{h,1}\big(\frac{t}{T} - \frac{a}{T}\big)}_{\bigO(1)},
\end{align*}
using \ref{eqn: O1 sum}. So
\begin{align}\label{eqn: remainder block 3}
    {\mathsf \Sigma}_{3}^\Xi 
        &\lesssim  \frac{1}{Th\phi^2(h)} \Big(\frac{1}{T}  + h\phi(h)\Big)^2
        \lesssim \frac{1}{Th\phi^2(h)} \Big( \frac{1}{T^2}  + h^2 \phi^2(h) \Big)\nonumber\\
        &\lesssim \frac{1}{T^3 h \phi^2(h)} + \frac{h}{T} 
        \lesssim  \frac{1}{T h \phi^2(h)}.
\end{align}
Now, comparing (\ref{eqn: remainder block 1}), (\ref{eqn: remainder block 2}), and (\ref{eqn: remainder block 3}), we have
\begin{align}\label{eqn: remainder blocks overall order}
    \E[\Xi_{t,T}^2] \lesssim \frac{1}{Th^2\phi^2(h)}.
\end{align}
Therefore, following (\ref{eqn: big blocks overall order}), (\ref{eqn: small blocks overall order}), and (\ref{eqn: remainder blocks overall order}),  we get
\begin{align}\label{eqn: bound E of S square}
    \E\big[Z_{t,T}^2\big]
        &=  \bigO\Big(\frac{1}{Th^2\phi^2(h)} + h^2\Big).
\end{align}

\bibliographystyle{plain}

\end{document}